\documentclass[reqno,a4paper]{amsart}
\usepackage{microtype}
\usepackage[dvipsnames]{xcolor}
\usepackage{amssymb}
\usepackage{stmaryrd} 
\usepackage[utf8]{inputenc}
\usepackage{a4wide}
\usepackage{paralist,float}
\usepackage{enumitem}
\usepackage[all]{xy}
\usepackage{ifthen}
\usepackage{graphicx}
\usepackage{slashed}
\usepackage{comment}
\usepackage[colorlinks, linkcolor=blue, citecolor=blue]{hyperref}

\newcommand\N{\mathbb{N}}
\newcommand\R{\mathbb{R}}
\newcommand\p{\mathbb{P}}

\newcommand\Z{\mathbb{Z}}

\newcommand\del{\partial}
\newcommand\tensor{\otimes}
\newcommand{\Aut}{\textup{Aut}}

\newcommand{\Tor}{\textup{Tor}}
\newcommand{\Hom}{\textup{Hom}}
\newcommand{\es}{\textup{es}}
\newcommand{\ES}{\textup{ES}}

\newcommand{\sm}{\textup{s}}
\newcommand{\eql}{\textup{eql}}
\newcommand{\EQL}{\textup{EQL}}
\newcommand{\ql}{\textup{i}}
\newcommand{\Ext}{\textup{Ext}}
\newcommand{\Coker}{\textup{Coker}}

\newcommand{\hght}{\textup{ht}}
\renewcommand{\Im}{\textup{Im}}
\newcommand{\Id}{\textup{Id}}
\renewcommand{\Im}{\textup{im}}

\newcommand{\hra}{\hookrightarrow}

\renewcommand{\mod}{\mathrm{mod}}
\newcommand{\ohat}{\mathop{\stackrel{{}_{\wedge}}{\otimes}}\nolimits}

\newcommand{\An}[1]{\bigl\langle{#1}\bigr\rangle}

\newcommand{\wh}{\widehat}
\newcommand{\ol}{\overline}
\newcommand{\ul}{\underline}
\newcommand{\mcal}{\mathcal}
\newcommand{\xra}{\xrightarrow}

\newcommand{\rk}{\textup{rk}}

\newcommand{\GW}{GW}
\newcommand{\HomZ}{\mathrm{Hom}_{\Z}}
\newcommand{\HomZns}{\HomZ^{\rm ns}}
\newcommand{\ICHomZns}{\mathrm{Hom}^{\rm ns}_{\Z}}
\newcommand{\FP}{\mathbf{FP}}

\newcommand{\Ab}{\mathbf{Ab}}
\newcommand{\AbF}{\mathbf{Ab}}
\newcommand{\rank}{\textup{rk}}
\renewcommand{\Im}{\mathrm{Im}}
\newcommand{\Ker}{\mathrm{Ker}}
\newcommand{\h}{\textup{h}}
\renewcommand{\p}{\textup{p}}
\newcommand{\ev}{\mathrm{ev}}
\newcommand{\CC}{\mathcal{C}}

\theoremstyle{plain}
\newtheorem{theorem}{Theorem}[section]
\newtheorem{thm}[theorem]{Theorem}
\newtheorem{lemma}[theorem]{Lemma}
\newtheorem{lem}[theorem]{Lemma}
\newtheorem{corollary}[theorem]{Corollary}
\newtheorem{proposition}[theorem]{Proposition}
\newtheorem{prop}[theorem]{Proposition}
\newtheorem{conjecture}[theorem]{Conjecture}

\theoremstyle{definition}
\newtheorem{definition}[theorem]{Definition}
\newtheorem{defin}[theorem]{Definition}
\newtheorem{example}[theorem]{Example}

\theoremstyle{remark}
\newtheorem{remark}[theorem]{Remark}
\newtheorem*{remark*}{Remark}

\advance \topmargin by -\headheight
\advance \topmargin by -\headsep
     
\evensidemargin \oddsidemargin
\title{The Witt groups of extended quadratic forms over $\Z$}
\author[D. Crowley]{Diarmuid Crowley}
\address{School of Mathematics and Statistics,
University of Melbourne,
Parkville, VIC, 3010, Australia}
\email{dcrowley@unimelb.edu.au}
\author[Cs. Nagy]{Csaba Nagy}
\address{Max-Planck-Institut f\"{u}r Mathematik, Bonn, Germany}
\email{nagy@mpim-bonn.mpg.de}
\date{\today}

\begin{document}
\maketitle

\begin{abstract}
%%%%%%%%
We study quadratic form parameters $Q$ over the integers and extended quadratic forms with values in $Q$, which we call {\em $Q$-forms}.
Certain form parameters $Q$ appeared in Wall's work on the classification of almost closed $(n{-}1)$-connected $2n$-manifolds via $Q$-forms. 
Baues, Ranicki and Schlichting independently developed definitions of extended quadratic forms in more general settings; when restricted to the ring $\Z$, each of those definitions is equivalent to those studied here.

In this paper we classify all quadratic form parameters $Q$ over the integers, 
determine the category of quadratic form parameters $\FP$ and 
compute the Witt group functor,
\[ W_0 \colon \FP \to \AbF, 
\quad
Q \mapsto W_0(Q),\]
where $\AbF$ is the category of finitely generated abelian groups
and $W_0(Q)$ is the Witt group of nonsingular $Q$-forms.
\end{abstract}
%%%%%%%%%%%%%%%%%%%%%%%%%%%%%%%%%%%%%%%%%%%%%%%%%%%

\section{Introduction} \label{s:intro}
%%%%%%%%%%%%%%%%%%%%%%%%%%%%%%%%%%%%%%%%%%%%%%%%%%%
Quadratic refinements of 
symmetric and anti-symmetric
bilinear forms are important objects in algebra and topology, 
and extended quadratic forms provide a general setting for their study.
Let $X$ be an abelian group 
(in this paper we assume that all groups are finitely generated)
and let
$\lambda \colon X \times X \to \Z$
be a bilinear form on $X$ such that $\lambda(x, y) = \epsilon_\lambda \lambda(y, x)$ for all $x, y \in X$ and some 
fixed $\epsilon_\lambda \in \{\pm 1\}$.
To define a quadratic refinement of $\lambda$
we require the additional structure of a {\em quadratic form parameter}.
This is a triple $Q = (Q_e, \h, \p)$, often written
\[   Q_e \xra{~\h~} \Z \xra{~\p~} Q_e,\]
where $Q_e$ is an abelian group and 
$\h$ and $\p$ are homomorphisms
such that 
$\h \p \h = 2 \h$ and $\p \h \p = 2 \p$.  
The equation
$\h \p = (1{+}\epsilon_Q) \Id_\Z$ defines $\epsilon_Q \in \{\pm 1\}$,
the symmetry of $Q$, and we require $\epsilon_Q = \epsilon_\lambda$.
A quadratic refinement of $\lambda$ with values in $Q$ is a function $\mu \colon X \to Q_e$,
which for all $x, y \in X$ satisfies the equations
\[ 
\mu(x+y) = \mu(x) + \mu(y) + \p(\lambda(x, y))
\qquad \text{and} \qquad
 \h(\mu(x)) = \lambda(x, x).
 \]
The triple $(X, \lambda, \mu)$ is a called an {\em extended quadratic form over $Q$}, or simply a {\em $Q$-form}.

Definitions of quadratic 
form parameters $Q$ and
$Q$-forms, equivalent to those above for the ring $\Z$, 
were developed independently by Baues \cite{Bau1}, Ranicki \cite{R1} and Schlichting \cite{S} in far broader contexts.
We also refer the reader to 
Section~\ref{ss:scm} for a review of the uses of $Q$-forms and their Witt groups
in geometric topology.  

\subsection{Form parameters and their Witt groups} \label{ss:FP+WG}
%%%%%%%%%%%%%%%%%%%%%%%%%%%%%%%%%%%%%%%%%%%%%%%%%%%%%%%%
We shall usually abbreviate ``quadratic form parameter'' to ``form parameter''.
For $k \in \Z^+ \cup \{\infty\}$ 
we have the following standard symmetric form parameters,
\begin{equation} \label{eq:ind+}
 Q_+ :=  \Bigl( \Z \xra{2} \Z \xra{1} \Z \Bigr),
\quad
\Z^P_k:= \Bigl(\Z \oplus \Z_{2^k} \xra{(1,0)} \Z \xra{(2,-1)^{\rm t}}
\Z \oplus \Z_{2^k} \Bigr),
\quad \text{and} \quad
 Q^+ :=  \Bigl( \Z \xra{1} \Z \xra{2} \Z \Bigr),
\end{equation}
where $\Z_{2^\infty} := \Z$. For $k \geq 2$ an integer  we have the following standard anti-symmetric form parameters,
\begin{equation} \label{eq:ind-}
 Q_- :=  \Bigl( \Z_2 \xra{0} \Z \xra{1} \Z_2 \Bigr),
 \quad
\Z^\Lambda_k := \Bigl( \Z_{2^k} \xra{~~0~~} \Z \xra{2^{k-1}} \Z_{2^k} \Bigr)
\quad \text{and} \quad
 Q^- :=  \Bigl( 0 \xra{} \Z \xra{} 0 \Bigr).
\end{equation}
In addition, given a form parameter $Q = (Q_e, \h, \p)$ and an abelian group $G$, we can construct a new form parameter, the sum of $Q$ and $G$, which is called a {\em split form parameter} and is defined by
\begin{equation} \label{eq:split}
Q \oplus G := 
\Bigl( Q_e \oplus G \xra{(\h,0)} \Z \xra{(\p,0)^{\rm t}} Q_e \oplus G \Bigr).
\end{equation}

In fact, every form parameter is isomorphic to 
a split form parameter with $Q$ standard,
as we now explain.
A morphism 
$\alpha \colon P \to Q$ in $\FP$, the category of form parameters, is a homomorphism
$\alpha \colon P_e \to Q_e$ such that $\h_P = \h_Q \circ \alpha$ and $\p_Q = \alpha \circ \p_P$.  It is an isomorphism if the map $P_e \to Q_e$ is
an isomorphism.  
A {\em splitting} of a form parameter $P$ is an isomorphism $\alpha \colon P \to Q \oplus G$.
If $P$ admits a splitting with $G$ non-trivial, then $P$ is {\em decomposable}; if not $P$ is {\em indecomposable}.
A splitting $P \cong Q \oplus G$ is called {\em maximal} if $Q$ indecomposable.
In Theorem~\ref{t:qfp-class-intro} we prove that the standard form parameters are indecomposable and that they  form the complete set of indecomposable form parameters up to isomorphism. 
Theorem~\ref{t:qfp-class-intro} also proves that every form parameter $P$ has a maximal splitting $P \cong Q \oplus G$, in which $Q$ and $G$ are well-defined up to isomorphism.

The {\em Witt group} of a form parameter $Q$ is defined as follows.
A $Q$-form $(X, \lambda, \mu)$ is 
{\em nonsingular} if $(X, \lambda)$ is nonsingular
(in which case $X$ is free), and it is
{\em metabolic} if it admits a {\em lagrangian}; i.e.~a half rank summand $L < X$
such that $\lambda(L \times L) = 0$ and $\mu(L) = 0$.
There are natural notions of isometry and orthogonal 
sum of $Q$-forms (see Section \ref{ss:Q_forms}),
and two nonsingular $Q$-forms are called {\em Witt equivalent}
if they become isometric after the addition of nonsingular metabolic $Q$-forms.
The Witt group of $Q$ is the set of Witt equivalence class of nonsingular $Q$-forms under the operation of orthogonal sum:
\[ W_0(Q):= \left( \bigl\{ [X, \lambda, \mu] \, \big| \,
\text{$(X, \lambda, \mu)$ is a nonsingular $Q$-form} \bigr\}, 
\oplus \right)\]

For example, the Witt-groups of
$Q_+, Q^+, Q_-$ and $Q^-$ are equal
to the $L$-groups of $\Z$ as follows:
\[
\begin{tabular}{c||c|c|c|c}
$Q$ & $Q_+$ & $Q^+$ & $Q_-$ & $Q^-$ \\
\hline
$W_0(Q)$ & $ L_0(\Z) \cong 8 \Z$ & $L^0(\Z) \cong \Z$ & $L_2(\Z) \cong \Z_2$ & $L^2(\Z) = 0$ \\
\end{tabular}
\]

If $\alpha \colon P \to Q$ is a morphism  of form parameters and $(X, \lambda, \mu)$ is a
$P$-form, then $(X, \lambda, \alpha \circ \mu)$ is a $Q$-form.  It follows that
$\alpha$ induces a homomorphism $W_0(\alpha) \colon W_0(P) \to W_0(Q)$, 
and so
taking Witt groups defines a functor
\[ W_0 \colon \FP \to \AbF,
\quad Q \mapsto W_0(Q),  \]
where $\AbF$ is the category of finitely generated abelian groups.
The main aim of this paper is to compute these Witt-groups and the Witt-group functor, 
and we achieve this in Theorems~\ref{t:W_0-ind},~\ref{t:W_0-split} and \ref{t:ES+EQL}.
%

%%%%%%%%%%%%%%%%%%%%%%%%%%%%%%%%%%%%%%%%%%%%%%%%
%
%
%

\begin{theorem} \label{t:W_0-ind}
%%%%%%%%%%%%%%%%%
For $k \in \N \cup \{\infty\}$
and $l \in \N$,
the Witt groups of 
indecomposable 
form parameters are given by 
\[ \sigma \colon W_0(Q_+) \xra{\cong} 8 \Z,
\qquad \sigma \oplus \rho_k \colon W_0(\Z^P_{k+1}) \xra{\cong} \Z \oplus \Z_{2^k},
\qquad \sigma \colon W_0(Q^+) \xra{\cong} \Z ,
\]
\[ c \colon W_0(Q_-) \xra{\cong} \Z_2,
\qquad
W_0(\Z^\Lambda_{l+2}) = 0
\qquad \text{and} \qquad
W_0(Q^-) = 0. 
\]
Here $\sigma$ is the signature homomorphism, 
$\infty{+}1 := \infty$, 
$\rho_k \colon W_0(\Z^P_{k+1}) \to \Z_{2^{k}}$
is a homomorphism defined using the squares of certain characteristic elements
(see Definitions~\ref{def:rho-inf} and~\ref{def:rho-k}) and $c$ is the Arf invariant. 
\end{theorem}

\noindent
{\em Remark.}
In Theorem~\ref{t:wf-ind} we compute $W_0$ on the full subcategory of indecomposable form parameters.

%%%%%%%%%%%%%%%%%%%%%%%%%%%%%%%%%%%%%%%%%%%%%%%%
\smallskip
For a split form parameter $Q \oplus G$, there are obvious morphisms
$i \colon Q \to Q \oplus G$ and $r \colon Q \oplus G \to Q$, with $r i = \Id_Q$,
and this ensures that 
\begin{equation} \label{eq:W(split)}
W_0(Q \oplus G) = W_0(Q) \oplus W^Q_0(G),
\end{equation} 
where $W^Q_0(G) := \Ker(W_0(r))$ is the {\em reduced Witt group} of $Q \oplus G$.  
A morphism $\alpha \colon Q_1 \oplus G_1 \to Q_2 \oplus G_2$ between split form parameters is called {\em split}, if $\alpha = \alpha_q \oplus \alpha_l$ 
for a morphism $\alpha_q \colon Q_1 \to Q_2$
and a homomorphism $\alpha_l \colon G_1 \to G_2$,
and it is clear that the homomorphisms induced on Witt groups by split morphisms preserve the splitting in \eqref{eq:W(split)}.
To compute reduced Witt groups
we will use the quadratic tensor product,
which is a covariant functor,
\[ \otimes_\Z \colon \AbF \times \FP \to \AbF, 
\quad (G, Q) \mapsto G \otimes_\Z Q, \]
defined by Baues \cite[\S 4]{Bau1}, 
and which we recall in Section~\ref{ss:qtp}.  
For example, there are natural isomorphisms
\[ G \otimes_\Z Q^+ \to \Gamma(G) 
\qquad \text{and} \qquad 
G \otimes_\Z Q_- \to \Lambda_1(G),\]
where $\Gamma(G)$ is Whitehead's universal quadratic functor~\cite{Wh}, 
and $\Lambda_1(G) := (G \otimes G)/(g \otimes h + h \otimes g)$.
We review the computation of $G \otimes_\Z Q$ for all indecomposable form parameters
$Q$ in Section~\ref{ss:compWG} and
in 
Section~\ref{ss:RWG} 
we define a homomorphism, called the
{\em split norm}, $\tilde F \colon W_0(Q \oplus G) \to G \otimes_\Z Q$, 
which is natural in $Q$ and $G$.

\begin{theorem} \label{t:W_0-split}
%%%%%%%%%%%%%%%%%
%
For a split form parameter $Q \oplus G$,
the restriction of the split norm
to $W^Q_0(G)$,
\[\tilde F|_{W^Q_0(G)} \colon W^Q_0(G) \to G \otimes_\Z Q,\]
is an isomorphism.
Hence the map
$H := (W_0(r), \tilde F) \colon W_0(Q \oplus G) \to W_0(Q) \oplus (G \otimes_\Z Q)$
is an isomorphism, 
natural in $Q$ and $G$,
such that for any split morphism 
$\alpha = \alpha_q \oplus \alpha_l \colon Q_1 \oplus G_1 \to Q_2 \oplus G_2$ we have
\[ H \circ W_0(\alpha) \circ H^{-1} = W_0(\alpha_q) \oplus (\alpha_l \otimes_\Z \alpha_q)  \colon W_0(Q_1) \oplus (G_1 \otimes_\Z Q_1)
\to W_0(Q_2) \oplus (G_2 \otimes_\Z Q_2).
\]
\end{theorem}
\noindent
For an arbitrary form parameter $P$, Theorems~\ref{t:W_0-ind} and \ref{t:W_0-split} provide an effective procedure to compute the Witt-group $W_0(P)$ by choosing a maximal splitting $P \cong Q \oplus G$. The details of this procedure are given in  Section~\ref{ss:compWG}.

\smallskip
To understand 
the Witt group functor on all morphisms, 
we first classify the category $\FP$.
If there is a morphism $\alpha \colon P \to Q$,
then $\epsilon_P = \epsilon_Q$ and so $\FP$ splits as the disjoint union of
$\FP_+$ and $\FP_-$, the categories of symmetric and anti-symmetric form parameters.  
We will show that each of $\FP_+$ and $\FP_-$ is equivalent to a simpler category. The {\em linearisation} of a form parameter $Q$ is the abelian group $SQ := Q_e/\p(\Z)$, and when $Q$ is symmetric, we consider the homomorphism induced by $\h \colon Q_e \to \Z$,
\[ v_Q \colon SQ \to \Z_2,
\quad [q] \mapsto \h(q)~\mathrm{mod}~2.\]
When $Q$ is anti-symmetric, $\p$ induces a homomorphism
\[ v'_Q \colon \Z_2 \to Q_e,
\quad [n] \mapsto \p(n).\]

We call the homomorphism $v_Q$ or $v'_Q$ the {\em quasi-Wu class} of $Q$;
this terminology is justified by Proposition~\ref{p:Q_beta}.
Let $\AbF/\Z_2$ and $\Z_2/\AbF$ denote the categories of abelian groups with homomorphisms to and from $\Z_2$, respectively. A morphism $\alpha \colon P \to Q$ of symmetric form parameters induces a homomorphism $S\alpha \colon SP \to SQ$ such that $v_P = v_Q \circ S\alpha$, and if $\alpha \colon P \to Q$ is a morphism of anti-symmetric form parameters, then $v'_Q = \alpha \circ v'_P$, showing that the quasi-Wu classes are functorial. We have the following

\begin{theorem} \label{t:FP_class}
%%%%%%%%%%%%%%%%%%%
The functors 
\[ v_{(-)}  \colon \FP_+ \to \AbF/\Z_2, \quad Q \mapsto v_Q
\qquad \text{and} \qquad
v'_{(-)}  \colon \FP_- \to \Z_2 / \AbF,  \quad Q \mapsto v'_Q \]
are equivalences of categories.
\end{theorem}

By combining Theorem \ref{t:FP_class} with the classification of the categories $\AbF/\Z_2$ and $\Z_2 / \AbF$, we obtain an explicit classification of form parameters.

\begin{thm}
\label{t:qfp-class-intro}
a) The standard form parameters in \eqref{eq:ind+} and \eqref{eq:ind-} are indecomposable and pairwise non-isomorphic. Any indecomposable form parameter is isomorphic to one of them. 

b) Every form parameter $P$ has a maximal splitting $P \cong Q \oplus G$. Moreover, $Q$ and $G$ are each well-defined up to isomorphism. 
\end{thm}

%%%%%%%%%%%%%%%%%%%%%%%%%%%%%%%%%%%%%%%%%%%%%
To describe the functor $W_0$ completely for general form parameters and morphisms,
we define the {\em extended symmetrisation} and {\em extended quadratic lift} morphisms
\[ \es \colon P \to Q^+ \oplus SP
\qquad \text{and} \qquad
\eql \colon Q_- \oplus P_e \to P,\]
for symmetric and anti-symmetric form parameters $P$ respectively. 
Applying Theorem~\ref{t:FP_class},
we note that for $\epsilon = \pm 1$, 
the form parameter $Q_\epsilon$ is an initial object in $\FP_\epsilon$, while $Q^\epsilon$ is terminal. 
So for any $\epsilon$-symmetric form parameter $Q$ there are unique morphisms
\[ Q_\epsilon \xra{\ql} Q \xra{\sm} Q^\epsilon, \]
where $\sm \circ \ql \colon Q_\epsilon \to Q^\epsilon$ is the canonical morphism.
Then $\es \colon P \to Q^+ \oplus SP$ is defined as the sum of $\sm$ and the linearisation map $P_e \to SP$ and $\eql \colon Q_- \oplus P_e \to P$ is the sum of $\ql$ and the identity $P_e \to P_e$; see Definitions~\ref{d:es} and~\ref{d:eql}.

By Theorem \ref{t:W_0-split} there are natural isomorphisms $W_0(Q^+ \oplus SP) \cong \Z \oplus \Gamma(SP)$ and $W_0(Q_- \oplus P_e) \cong \Z_2 \oplus \Lambda_1(P_e)$, and we will fix these identifications. We will show that $W_0(\es)$ is always injective and $W_0(\eql)$ is always surjective 
(see Lemmas \ref{l:Q_to_Q+SQ} and \ref{l:gen-antis-surj}). We will then get a natural description of $W_0(P)$ by identifying the image of $W_0(\es)$ or the kernel of $W_0(\eql)$ as a natural subgroup of $\Z \oplus \Gamma(SP)$ or $\Z_2 \oplus \Lambda_1(P_e)$.  For this we require some additional definitions.
Given an abelian group $A$ and a homomorphism $v \colon A \to \Z_2$, we define (see Definition~\ref{d:Sigma(v)}) the subgroup
\[
\Sigma(v) := \left< (1,x \otimes 1),(0,[k_1, k_2] \otimes 1),(8,0) \bigm| x \in v^{-1}(1), k_1, k_2 \in \Ker(v) \right> \leq \Z \oplus \Gamma(A).
\]
Given a homomorphism $v' \colon \Z_2 \to A$, we define (see Definition~\ref{d:Lambda(v')}) the subgroup 
\[ 
K(v') := \An{ (1, v'(1) \ohat v'(1)), (0, x \ohat x + x \ohat v'(1)) \mid x \in A } \leq \Z_2 \oplus \Lambda_1(A) 
\]
and the quotient group
\[ 
\Lambda(v') := (\Z_2 \oplus \Lambda_1(A)) / K(v') .
\]

It is elementary to check that the assignments $v \mapsto \Sigma(v)$, $v' \mapsto K(v')$ and $v' \mapsto \Lambda(v')$ define functors $\Sigma \colon \AbF/\Z_2 \rightarrow \AbF$ and $\Lambda \colon \Z_2/\AbF \rightarrow \AbF$. We will also include more explicit descriptions of these groups, see the discussion after Theorems \ref{t:sigma-comp} and \ref{t:lambda-comp}.
Our next result gives a natural identification of the Witt group functor (see Corollaries \ref{cor:nat+} and \ref{cor:nat-}).

\begin{theorem} \label{t:ES+EQL}
There are natural isomorphisms, induced by $W_0(\es)$ and $W_0(\eql)$ respectively: 
\[
\begin{aligned}
W_0 \big| _{\FP_+} &\cong \Sigma \circ v_{(-)} \colon \FP_+ \to \AbF \\
W_0 \big| _{\FP_-} &\cong \Lambda \circ v'_{(-)} \colon \FP_- \to \AbF
\end{aligned}
\]
\end{theorem}

\noindent
Theorem \ref{t:ES+EQL} allows us to determine 
the induced map $W_0(\alpha)$ for 
any morphism of form parameters $\alpha$.
\begin{corollary} \label{c:ES+EQL}
%%%%%%%%%%%%%%%
a) For any symmetric form parameter $P$, 
$W_0(\es) \colon W_0(P) \to \Z \oplus \Gamma(SP)$ is injective with image $\Sigma(v_P)$,
and so for any morphism $\alpha \colon P_1 \to P_2$ of symmetric form parameters there is a commutative diagram
\[
\xymatrix{
~~~~W_0(P_1) \ar[rr]^\cong \ar@/^5ex/[rrrr]^{W_0(\es)} \ar[d]^{W_0(\alpha)} &&
\Sigma(v_{P_1}) \ar[d]^{\Sigma(S\alpha)} \ar@{^(->}[rr] &&
\Z \oplus \Gamma(SP_1) \ar[d]^{\Id \oplus \Gamma(S\alpha)} \\
~~~~W_0(P_2) \ar[rr]^\cong \ar@/_5ex/[rrrr]_{W_0(\es)} &&
\Sigma(v_{P_2}) \ar@{^(->}[rr] &&
\Z \oplus \Gamma(SP_2).
}
\]

b) For any anti-symmetric form parameter $P$, $W_0(\eql) \colon \Z_2 \oplus \Lambda_1(P_e) \to W_0(P)$
is surjective with kernel $K(v'_P)$,
hence there is an induced isomorphism $W_0(\eql) \colon \Lambda(v'_P) \to W_0(P)$.
So for any morphism 
$\alpha \colon P_1 \to P_2$ of anti-symmetric form parameters,
there is a commutative diagram
\[ \hskip -0.9cm
\xymatrix{
\Z_2 \oplus \Lambda_1(P_{1e}) 
\ar[d]^{\Id \oplus \Lambda_1(\alpha)} \ar@{>>}[rr]
\ar@/^5ex/[rrrr]^{W_0(\eql)} &&
\Lambda(v'_{P_1}) \ar[d]^{\Lambda(\alpha)} \ar[rr]^-\cong &&
W_0(P_1) \ar[d]^{W_0(\alpha)} \\
 \Z_2 \oplus \Lambda_1(P_{2e}) \ar@{>>}[rr]
 \ar@/_5ex/[rrrr]_{W_0(\eql)}&&
\Lambda(v'_{P_2}) \ar[rr]^-\cong &&
W_0(P_2).
}
\]
\end{corollary}

%%%%%%%
%
%

\subsection{Grothendieck-Witt groups and absorbing forms} \label{ss:GWG}
%%%%%%%%%%%%%%%%%%%%%%%%%%%%%%%%%%%%%%%%%%%%%%%%%%%%%
The study of Witt groups of extended quadratic forms has recently received impetus from the breakthrough work of Calm\'es, Dotto, Harpaz, Hebestreit, Land, Moi, Nardin, Nikolaus and Steimle on Grothendieck-Witt groups, which defines spaces and spectra in a very general setting whose group of connected components are Witt groups of certain form parameters \cite{9author}.

\subsubsection*{Grothendieck-Witt groups} \label{ss:GW-Intro}
%%%%%%%%%%%%%%%%%%%%%%%%%%%%%%%%%%%%%%%%%%%%%%%%%%%%%
The Grothendieck-Witt group $GW_0(Q)$ of a form parameter $Q$ can be defined prior to the Witt group $W_0(Q)$.
Let $\ul \mu = (X, \lambda, \mu)$ be a nonsingular $Q$-form,
$\{ \ul \mu \}$ denote the isometry class of $\ul \mu$ and $\ICHomZns(Q) := \big\{ \{ \ul \mu \} \big\}$ be the set of isometry classes of nonsingular $Q$-forms, which is a commutative monoid under the operation of orthogonal sum.
The Grothendieck-Witt group of $Q$ is the Grothendieck group of this monoid,
\[ GW_0(Q) := {\mcal Gr}(\ICHomZns(Q), \oplus),\]
and there is a natural surjective homomorphism
\[ \pi \colon GW_0(Q) \to W_0(Q),
\quad (\{\ul \mu_0\}, \{\ul \mu_1\}) \mapsto 
[\ul \mu_0] - [\ul \mu_1].
\]
Taking the rank of the underlying abelian group of a $Q$-form defines a homomorphism
$\rk \colon GW_0(Q) \to \Z$,
and we let $\sigma \colon W_0(Q) \to \Z$ denote signature homomorphism 
(by convention, $\sigma = 0$ if $Q$ is anti-symmetric), 
$\varrho_2 \colon \Z \to \Z_2$
denote reduction mod~$2$ and set $\sigma_2 := \varrho_2 \circ \sigma$.

\begin{theorem} \label{t:GW_split}
%%%%%%%%%%%%%%%%%%
For all form parameters $Q$, there is a short exact sequence
\[ 0 \to GW_0(Q) \xra{\rk \oplus \pi} \Z \oplus W_0(Q) \xra{(\varrho_2, \sigma_2)} \Z_2 \to 0,\]
the surjection $\pi \colon GW_0(Q) \to W_0(Q)$ is split and there is an isomorphism
\[ GW_0(Q) \cong 2\Z \oplus W_0(Q),\]
where $\rk|_{\ker(\pi)} \colon \ker(\pi) \to 2\Z$ is also an isomorphism.
Moreover, the splitting 
$GW_0(Q) \cong 2\Z \oplus W_0(Q)$ is natural
on $\FP_-$ and on the full sub-category of $\FP_+$ consisting of symmetric form
parameters $Q$ with $v_Q = 0$.
\end{theorem}

\subsubsection*{Absorbing forms}
%%%%%%%%%%%%%%%%%%%%%%%%%%%%%%%%%%%%%%%%%%%%%%%%%%%%%
Given an nonsingular $Q$-form $\ul \mu = (X, \lambda, \mu)$, the group of  automorphisms of $\ul \mu$,
\[ \Aut(\ul \mu) = \{ A \in \Aut(X) \mid \lambda = \lambda \circ (A \times A) ~\text{and}~
\mu = \mu \circ A \}, \]
is an important object of study in algebra and topology.
For $k \ul \mu := \ul \mu \oplus \dots \oplus \ul \mu$, the $k$-fold orthogonal sum of $\ul \mu$ with itself, we often know more about $\Aut(k \ul \mu)$ and
$\lim_{k \to \infty} \Aut(k \ul \mu)$ than we do about $\Aut(\ul \mu)$, and for this and other reasons, 
it can be useful
to know whether a nonsingular $Q$-form $\ul \mu$ is 
{\em absorbing}, which means that
for every nonsingular $Q$-form $\ul \mu'$, 
there is a primitive embedding 
$\ul \mu' \hookrightarrow k \ul \mu$
for some positive integer $k$.
In this case there is an isomorphism $k \ul \mu \cong \ul \mu' \oplus \ul \mu^{'\perp}$, 
where $\ul \mu^{'\perp}$ is the orthogonal complement to an embedding of $\ul \mu' \hookrightarrow k \ul \mu$,
and so an injection $\Aut(\ul \mu') \to \Aut(k \ul \mu)$.
To determine whether a $Q$-form $\ul \mu = (X, \lambda, \mu)$ is absorbing we need two properties of $Q$-forms.
We say that $\ul \mu$ is  {\em indefinite} if $(X, \lambda)$ is indefinite
and that $\ul \mu$ is {\em full} if $S\mu \colon X \to SQ$ is onto.

\begin{theorem} \label{t:absorbing}
%%%%%%%%%%%%%%%%%
For any form parameter $Q$, a nonsingular $Q$-form $\ul \mu$ is absorbing if and only if $\ul \mu$ is indefinite and full.
\end{theorem}

\begin{remark}
%%%%%%%%
If a $Q$-form $\ul \mu = (X, \lambda, \mu)$ is not full, then there is a unique form parameter $Q'\subset Q$ with
linearisation $SQ' = (S\mu)(X)$ and such
that $\ul \mu$ lifts uniquely to a full $Q'$-form.  
To define $Q'$, consider the surjection
$Q_e \to SQ$ and the subgroup $(S\mu)(X) \subset SQ$,
and let $Q'_e$ be the pre-image of
$(S\mu)(X)$. 
Then $\p_Q(\Z) \subseteq Q'_e$ and $Q' := (Q'_e, \h_Q|_{Q'_e}, \p_Q)$ is
the desired form parameter.
\end{remark}

\subsection{Algebraic Poincar\'e cobordism} \label{ss:cb}
%%%%%%%%%%%%%%%%%%%%%%%%%%%%%%%%%%%%%%%%%%%%%%%%%%
Witt groups of extended quadratic forms appear in Ranicki's theory of algebraic surgery \cite{R1}, as we now quickly recall.
A chain bundle $(B_*, \beta)$ over $\Z$ is a chain complex $B_*$ of free abelian groups, which is bounded below, together with a degree $0$ hyper-quadratic structure $\beta$ on the dual of $B_*$.
Chain bundles were introduced by Weiss \cite{We} as algebraic models of stable vector bundles, or more generally, stable spherical fibrations.  
Some symmetric Poincar\'e complexes can
be enriched with $(B_*, \beta)$-structures and the bordism group of such complexes
is the $\beta$-twisted $L$-group $L^n(B_*, \beta)$.  
For each natural number $q$, a chain bundle
$(B_*, \beta)$ defines a $(-1)^q$-symmetric form parameter $Q_\beta(q)$ with 
$SQ_\beta(q) = H_q(B_*)$ \cite[\S 10]{R1},
and by \cite[\S 9]{R1} there is a natural homomorphism
\begin{equation} \label{eq:LBbeta+Witt}
\eta_\beta \colon W_0(Q_\beta(q)) \to L^{2q}(B_*, \beta). 
\end{equation}
In Conjecture~\ref{c:LBbeta}, we propose that $\eta_\beta$ is split injective if $H_{q-1}(B_*) = 0$, and 
an isomorphism if $H_i(B_*) = 0$ for $i < q$.

In Proposition \ref{p:Q_beta} we show that the algebraic Wu class $\wh v_{2k}(\beta) \in H^{2k}(B_*; \Z_2)$ of $(B_*, \beta)$ (defined by Ranicki in
e.g.\ \cite[\S 4]{R1}), determines the quasi-Wu classes of both $Q_\beta(2k)$ and $Q_\beta(2k{-}1)$.
Applying Theorem~\ref{t:FP_class}, it follows that $\wh v_{2k}(\beta)$ determines the isomorphism classes of the form parameters $Q_\beta(2k)$ and $Q_\beta(2k{-}1)$,
and so the isomorphism classes of the Witt groups $W_0(Q_\beta(2k))$ and
$W_0(Q_\beta(2k{-}1))$.
Hence Conjecture~\ref{c:LBbeta} and Theorem \ref{t:ES+EQL} give a conjectural calculation of
$L^{2q}(B_*, \beta)$ when $H_i(B_*) = 0$ for all $i < q$.

\subsection{$Q$-forms and the classification of simply-connected manifolds} \label{ss:scm}
%%%%%%%%%%%%%%%%%%%%%%%%%%%%%%%%%%%%%%%%%%%%%%%%%%
$Q$-forms appeared at the very start of the modern classification of manifolds in Wall's paper classifying $(q{-}1)$-connected $2q$-dimensional manifolds with connected boundary for $q \geq 3$ \cite{Wa1}.
In this subsection we give an historical and selective overview
of the use of $Q$-forms in the classification of simply-connected oriented manifolds.

\subsubsection*{The classification of $(q{-}1)$-connected $2q$-manifolds for $q \geq 3$}
%%%%%%%%%%%%%%%%%%%%%%%%%%%%%%%%%%%%%%%%%%%%%%%%%%%%%%%%
We follow Baues \cite[\S 8]{Bau2} in stating 
Wall's main classification result from \cite{Wa1} by first defining the form parameter
\[ \pi_{q-1}\{SO_q\} := \bigl( \pi_{q-1}(SO_q) \xra{~\h~} \Z \xra{~\p~} \pi_{q-1}(SO_q) \bigr).\]
Here $\h(\gamma) \in \Z = H^q(S^q)$ is the Euler class of an oriented rank~$q$ vector bundle over $S^q$ 
whose clutching function is in the homotopy class $\gamma$, and $\p(1) = \tau_q$, the homotopy class of the clutching function 
of the tangent bundle of the $q$-sphere.
Now let $W$ be a compact, oriented, $(q{-}1)$-connected $2q$-manifold with connected
boundary $\del W$ and $q$th integral 
homology group $H_q(W)$.
The intersection form of $W$ is a $(-1)^q$-symmetric bilinear form
\[\lambda_W \colon H_q(W) \times H_q(W) \to \Z.\] 
Applying results of Haefliger \cite{H}, Wall \cite{Wa1} proved 
that 
there is a function
\[ \mu_W \colon H_q(W) \to \pi_{q-1}(SO_q)\]
such that the triple 
$\ul{\mu}_W := (H_q(W), \lambda_W, \mu_W)$ is a $\pi_{q-1}\{SO_q\}$-form, 
and that 
the diffeomorphism classes of $(q{-}1)$-connected $2q$-manifolds with connected boundaries are in bijection with the isomorphism classes of their $\pi_{q-1}\{SO_q\}$-forms.
(When $q = 2$ the topological situation is vastly more complicated,
but the algebra is as simple as possible since $\pi_1\{SO_2\} = Q^+$ and so $\pi_1\{SO_2\}$-forms are 
nothing but symmetric bilinear forms.)

The Witt groups and Grothendieck-Witt groups of $\pi_{q-1}\{SO_q\}$-forms appear when 
$\del W$ is a homotopy sphere, for this is equivalent to $\ul \mu_W$ being nonsingular.  
Wall \cite[Theorem 2]{Wa1} computed the Grothendieck-Witt groups of 
$\pi_{q-1}\{SO_q\}$-forms and implicitly also their Witt groups, bearing 
Theorem~\ref{t:GW_split} in mind.
Moreover, \cite[Theorem 3]{Wa1} states that the diffeomorphism class of $\del W$ only depends on the Witt class of $\ul \mu_W$ in $W_0\bigl(\pi_{q-1}\{SO_q\} \bigr)$, and Wall proved that there is a well-defined homomorphism
\[ \del_{2q} \colon 
W_0 \bigl( \pi_{q-1}\{SO_q\} \bigr) \to \Theta_{2q-1},\]
where $\Theta_{2q-1}$ denotes the group of homotopy $(2q{-}1)$-spheres and $\del_{2q}$
maps $[X, \lambda, \mu]$ the diffeomorphism class of  the homotopy sphere $\Sigma$ with $\Sigma = \del W$ and 
$[H_q(W), \lambda_W, \mu_W] = [X, \lambda, \mu]$ 
(the homomorphism $\del_{2q}$ appears explicitly in \cite[\S 17]{Wa2}).
Wall \cite[Theorem 4]{Wa1}
initiated the study of the homomorphism $\del_{2q}$, which is a very delicate problem with a long history.
For example, when $q = 2^{k}{-}1$ and $k \geq 2$, computing $\del_{2q}$ requires solving the Kervaire-Invariant-One problem in dimension 
$2^{k+1}{-}2$.
Following Wall, many authors worked on the homomorphism $\del_{2q}$, with this effort recently culminating in a papers of Burklund and Senger \cite{B-S}, which computes $\del_{2q}$, except for the case of $\del_{126}$,
and Lin, Wang and Xu \cite{L-W-X}, gives that
$\del_{126} = 0$.

If $\del_{2q}([\ul \mu_M]) = 0$, then there are diffeomorphisms $f \colon \del W \cong S^{2q-1}$, $W$ can be completed to a closed $(q{-}1)$-connected $2q$-manifold $M = W \cup_f D^{2q}$, and Wall \cite{Wa1} also proved that such manifolds are classified up to the action of $\Theta_{2q}$ via connected sum by their $\pi_{q-1}\{SO_q\}$-forms.  This reduced the classification of closed $(q{-}1)$-connected $2q$-manifolds to the determination of their inertia groups, a task which was recently finally completed by Senger and Zhang \cite{S-Z}.

\subsubsection*{Classical, modified and extended surgery}
%%%%%%%%%%%%%%%%%%%%%%%%%%%%%%%%%%%%%%%%%%%%%%
Wall's classification of $(q{-}1)$-connected $2q$-manifolds was one of the precursors to
the development of the surgery classification of manifolds in general.
A key idea in surgery is to study normal maps from manifolds, not just manifolds
on their own.  Let $\xi \colon B \to BO$ be a fibration, where $B$ is homotopy equivalent to a simply-connected $CW$-complex.  In what follows, we shall abuse notation and also use $B$ to denote one of
its $CW$-approximations.
For $M$ a compact smooth manifold, 
a {\em normal map} from $M$ over $(B, \xi)$,
denoted $\ol \nu \colon M \to B$, is a map of stable vector bundles from the stable normal bundle of $M$ to $\xi$.
An important example of the above comes from classical surgery, 
when $B = X$ is a $2q$-dimensional Poincar\'e complex and the normal map $\ol \nu$ has degree one and is $q$-connected.
Then for $\epsilon = (-1)^q$, 
Browder \cite[Ch.\ III]{Browd2} proved that the surgery kernel
\[ K_q(\ol \nu) := H_{q+1}(X, M) \]
supports a nonsingular $Q_\epsilon$-form $(K_q(\ol \nu), \lambda, \mu_{\ol \nu})$,
where we use $\ol \nu$ to regard $M$ as a subspace of $X$ and $\lambda$ is induced by the intersection form of $M$.
Moreover, by \cite[Theorem 2.1 Ch.\ IV]{Browd2} for $q \geq 3$ the normal map $(M, \ol \nu)$ is $(X, \xi)$-bordant to a normal homotopy equivalence if and only if
\[ [K_q(\ol \nu), \lambda, \mu_{\ol \nu}] = 0 \in W_0(Q_\epsilon) = L_{1-\epsilon}(\Z).\]

In the more general setting of modified surgery, 
one seeks to classify $q$-connected $2q$-dimensional normal maps
$\ol \nu \colon M \to B$, which are called {\em normal $(q{-}1)$-smoothings}. An important result of
Kreck~\cite[Corollary 3]{K} states that $(B, \xi)$-bordant normal 
$(q{-}1)$-smoothings are stably diffeomorphic, where a stable diffeomorphism between two $2q$-manifolds $M_0$ and $M_1$ with the same Euler characteristic is a diffeomorphism $M_0 \sharp_r (S^q \times S^q) \to M_1 \sharp_r (S^q \times S^q)$.
Kreck also proved that cancellation of any $(S^q \times S^q)$-summands is always possible when $q$ is odd and also when $q$ is even, 
provided the genus of either manifold is at least 1,
\cite[Theorem 5]{K}.
Here cancellation is the phenomenon where the existence of a diffeomorphism $M_0 \sharp_r (S^q \times S^q) \to M_1 \sharp_r (S^q \times S^q)$ implies the existence of a diffeomorphism $M_0 \to M_1$,
and the {\em genus} of an oriented $2q$-manifold $M$, as defined by Galatius and Randal-Williams \cite[Definition 3.1]{G-RW}, is the largest natural number $r$ such that there is a diffeomorphism $M \to N \sharp_r (S^q \times S^q)$, for some other oriented $2q$-manifold $N$.
The above demonstrates the importance of finding invariants which can distinguish stably-diffeomorphic $4k$-manifolds of genus zero.
While Kreck \cite[Proposition 6]{K} generalised Browder's result above by showing that a $Q$-form $(K_q(\ol \nu), \lambda, \mu_{\ol \nu})$ is defined on $K_q(\ol \nu)$, where $Q = Q_\epsilon$ if $v_{q+1}(\pi_{q+1}(B)) = 0$ and $Q = Q^-$ if $v_{q+1}(\pi_{q+1}(B)) \neq 0$,
when $B$ is not a $2q$-dimensional Poincar\'e complex, 
the $Q$-form 
$(K_q(\ol \nu), \lambda, \mu_{\ol \nu})$ may be singular and this can significantly complicate the algebra of surgery obstructions; see \cite[\S 5]{K}.

One approach to handling the more complicated algebra associated to normal $(q{-}1)$-smoothings $\ol \nu \colon M \to B$ is via the theory of chain bundles, which provides more general form parameters $Q$ and $Q$-forms defined on $H_q(M)$.
Specifically, if we 
assume that $M$ admits a $CW$-structure such that $\ol \nu \colon M \to B$ is cellular, then the stable vector bundle $(B, \xi)$ defines a chain bundle $(B_*, \beta)$ over $C_*(B)$, the cellular chain complex of $B$, 
and the normal map $\ol \nu \colon M \to B$ equips $C_*(M)$, the cellular chain complex of $M$, with a normal $(B_*, \beta)$-structure by \cite[Theorem 3.4]{We}.
Moreover, as we saw in Section~\ref{ss:cb}, the chain bundle $(B_*, \beta)$ defines a form parameter $Q_\beta(q)$ and Ranicki shows \cite[\S 10 p.\ 250]{R1} how the normal $(B_*, \beta)$-structure on $C_*(M)$
defines a $Q_\beta(q)$-form 
\[ \ul \mu_{\ol \nu} := (H_q(M), \lambda_M, \mu_{\ol \nu}) \] 
on $H_q(M)$
\footnote{Ranicki assumes that $H_q(M)$ is torsion free, but this is to ensure the $Q$-form is nonsingular, and is not used in the definition.},
which is thus a $(B, \xi)$-diffeomorphism invariant of $\ol \nu \colon M \to B$.
In \cite[Theorem 1.2]{N2} the second author proves that when $q = 2k \geq 4$ is even, any $(B, \xi)$-bordism between $4k$-dimensional $(2k{-}1)$-smoothings with isomorphic $Q_\beta(q)$-forms is $(B, \xi)$-bordant to a $(B, \xi)$-h-cobordism, 
and so for $k \geq 2$, $4k$-dimensional
$(B, \xi)$-bordant $(2k{-}1)$-smoothings with isomorphic $Q_\beta(q)$-forms are diffeomorphic.
In \cite{CCPS} for $k \geq 2$, $Q_\beta(q)$-forms are used to distinguish $4k$-manifolds with the same Euler characteristic, 
which are stably diffeomorphic but not diffeomorphic.

%%%%%%%%%%%%%%%%%%%%%%%%%%%%
%

%%%%%%%%%%%%%%%%%
%

\subsection{A conjecture on cancellation for $Q$-forms} \label{ss:fw}
%%%%%%%%%%%%%%%%%%%%%%%%%%%%%%%%%%%%%%%%%%%%%%%%%%
In this subsection we consider the role of the Witt group $W_0(Q)$ in the classification of full nonsingular $Q$-forms.  
Recall that a $Q$-form $\ul \mu = (X, \lambda, \mu)$ is full if its linearisation $S\mu \colon X \to SQ$ is surjective.
We regard full $Q$-forms as algebraic analogues of closed $2q$-dimensional normal $(q{-}1)$-smoothings and $W_0(Q)$ as an analogue of the corresponding bordism group of closed normal $(q{-}1)$-smoothings.  
Hence we define the genus of a nonsingular full $Q$-form $\ul \mu$ to be the largest natural number $r$ such that there is an isomorphism 
$\ul \mu \to \ul \mu' \oplus H_{\epsilon_Q}(\Z^r)$,
where $\ul \mu'$ is a nonsingular full $Q$-form
and $H_{\epsilon_Q}(\Z^r)$ is the standard hyperbolic form
on $\Z^{2r}$; see Example~\ref{ex:hyp}.
In the light of Kreck's stable diffeomorphism and cancellation
theorems \cite[Corollary 3 \& Theorem 5]{K} reviewed in Section~\ref{ss:scm} above,
we make the following

\begin{conjecture} \label{c:alg_can}
%%%%%%%%%%%%%%%%%%%%%%%%%%%%%%%%%%%%%%%
Let $\ul \mu_0$ and $\ul \mu_1$ be nonsingular full $Q$-forms of the same rank, such that $[\ul \mu_0] = [\ul \mu_1] \in W_0(Q)$. Then
\begin{enumerate}
\item[a)] For some integer non-negative integer $m$, 
there is an isomorphism
$\ul \mu_0 \oplus H_{\epsilon_Q}(\Z^m) \to 
\ul \mu_1 \oplus H_{\epsilon_Q}(\Z^m)$;
\item[b)] If $\epsilon_Q = 1$ and $\mathrm{g}(\ul \mu_0) \geq 1$,
or if $\epsilon_Q = -1$, then then $\ul \mu_0$ and $\ul \mu_1$ are isomorphic.
\end{enumerate}
\end{conjecture}

\begin{remark} \label{r:Q-form-Class}
%%%%%%%%%%%%%%%
Conjecture~\ref{c:alg_can} is well-known for the classical form parameters $Q_\pm$ and $Q^\pm$; see e.g.\ \cite[Theorem 5.3]{M-H} for the symmetric case and \cite[Theorem 8.104]{L-M} for the anti-symmetric case.  For metabolic forms over general form parameters, Part a) of Conjecture~\ref{c:alg_can} is proven in Lemma~\ref{l:metabolic}.
\end{remark}

\subsection{Form parameters over more general rings} \label{ss:or}
%%%%%%%%%%%%%%%%%%%%%%%%%%%%%%%%%%%%%%%%%%%%%%%%%%%%%%%%%%%%%%%%%%
We hope that the results of this paper may be a source of ideas for future work over more general rings. 
The theory of form parameters over group rings $\Z[\pi]$ and more generally rings with involution goes back to Bak \cite{Bak} in special cases, and general theories were developed by Baues \cite{Bau1} and Ranicki \cite{R1} and Schlichting \cite{S}.  
For applications to non-simply-connected manifolds, one would like to develop the theory of form parameters $Q$ over group rings $\Z[\pi]$ and of $Q$-forms defined on $\Z[\pi]$-modules.   In particular, it seems likely that ``quasi-Wu-classes" can be defined more generally, and that an appropriate analogue of Theorem~\ref{t:FP_class} holds.
In this direction, Wolfgang Steimle informs us that 
results in Sections 1 \& 4 of \cite{9author} may
give another proof of Theorem~\ref{t:FP_class}
and also analogous results for other rings.
We also expect that analysing form parameters over general rings via indecomposable and linear parts will be a useful point of view for computing their Witt groups, but it seems likely that splittings $P \cong Q \oplus G$ will be replaced by short exact sequences $Q \to P \to G$.

\subsection*{Organisation} \label{ss:org}
%%%%%%%%%%%%%%%%%%%%%%%%%%%%%%%%%%%%%%%%%%%%%%%%%%
The rest of this paper is organised as follows.
In Section~\ref{s:prelims} we present basic definitions and examples of form parameters, $Q$-forms and Witt groups.
In Section~\ref{s:FPclass} we classify 
form parameters via their quasi-Wu classes, proving Theorems~\ref{t:FP_class} and \ref{t:qfp-class-intro}.  In Section~\ref{s:WG} we determine the Witt groups of indecomposable and split form parameters, proving 
Theorems~\ref{t:W_0-ind} and \ref{t:W_0-split},
and in Section~\ref{s:WGN} we give a natural computation of Witt groups, proving Theorem~\ref{t:ES+EQL}.
In Section~\ref{s:GW} we turn our attention to Grothendieck-Witt groups of form parameters, and prove 
Theorems~\ref{t:GW_split} and~\ref{t:absorbing}.
Finally, in Section~\ref{s:APC} we discuss Witt groups of form parameters in the context of Algebraic surgery over chain bundles.

\setcounter{tocdepth}{1}
\tableofcontents

\subsection*{Acknowledgements} \label{ss:ack}
%%%%%%%%%%%%%%%%%%%%%%%%%%%%%%%%%%%%%%%%%%%%%%%%%%
We would like to thank Wolfgang Steimle for several helpful comments on Grothendieck-Witt theory
and Mark Powell for a helpful comment on Algebraic Surgery.
We also thank two referees for many helpful comments, which have clarified several points and lead to many improvements.
The second author was supported by the Melbourne Research Scholarship and the EPSRC New Investigator grant EP/T028335/2.

\section{Preliminaries} \label{s:prelims}
%%%%%%%%%%%%%%%%%%%%%%%%%%%%%%%%%%%%%%%%%%%%%%%%%%%
In this section we establish basic definitions and results for form parameters and extended quadratic forms and their Witt groups.  We also give 
examples to illustrate these ideas.  For ease of exposition, we will repeat some definitions from the Introduction.
We defer similar preliminaries for Grothendieck-Witt groups to Section~\ref{s:GW}, where Grothendieck-Witt groups are treated separately.

\subsection{Quadratic form parameters} \label{subsec:quadratic_form_parameters}
%%%%%%%%%%%%%%%%%%%%%%%%%%%%%%%%%%%%%%%%%%%%%%%%%%%%%%%%%%%%%%%%%%%%%%%%%%%%%%%%%%%%%%%%%%%%%%%%%%%%%%
Quadratic form parameters provide the structure needed to define the quadratic refinements 
in extended quadratic forms.
The material in this subsection is adapted from, and builds on, Baues \cite[\S 2 \& \S 5]{Bau1}, \cite[\S 8]{Bau2} and Ranicki \cite[\S 10]{R1}.  

\subsubsection{Definitions and examples}
%%%%%%%%%%%%%%%%%%%%%%%%%%%%%%%%%%%%%%
We begin directly with the definition of a quadratic form parameter, discuss it, define morphisms of form parameters and give some 
important examples.

\begin{definition}
An {\em oriented quadratic form parameter} $Q = (Q_e, {\rm h}, {\rm p})$ (over the ground ring $\Z$) is a (finitely generated) abelian group $Q_e$ together with two homomorphisms ${\rm h}$ and ${\rm p}$,
\[
Q_e \xra{~\rm h~} \Z \xra{~\rm p~} Q_e \text{\,,}
\]
such that 
\[ 
{\rm h p h} = 2 {\rm h} \quad \text{~and~} \quad {\rm p h p} =  2 {\rm p} .
\]
The {\em symmetry} of $Q$, $\epsilon_Q \in \{\pm 1\}$, is determined by the equation ${\rm hp} - {\rm Id}_{\Z}  =\epsilon_Q {\rm Id}_{\Z}$.  
(To see that $\h \p(1) = 0$ or $2$,
note that if $\h \p(1) = d \neq 0$, then $p(1) \in Q_e$ much have infinite order and the equation $\p \h \p(1) = 2 \p(1) = d \p(1)$ ensures that $d = 2$.)  If $\epsilon_Q = 1$ we call $Q$ {\em symmetric} and if $\epsilon_Q = -1$ we call $Q$ {\em anti-symmetric}.
\end{definition}

\begin{remark}
An oriented quadratic form parameter $Q$ is a special case of a quadratic module as defined by Baues \cite[\S 2]{Bau1}, with $Q_{ee} = \Z$: $Q = (Q_e, \Z, {\rm h}, {\rm p}) = (Q_e, Q_{ee}, {\rm h}, {\rm p})$.  Henceforth we shall not explicitly mention ${\rm h}$ and ${\rm p}$, the ground ring $\Z$ or the orientation and so we use the letter $Q$ to refer to an oriented quadratic form parameter over $\Z$, which we now simply call ``a form parameter''.   We have used the term ``oriented" since we use the based abelian group $\Z$ and not simply an infinite cyclic group. 
\end{remark}

\begin{remark}
An oriented quadratic form parameter $Q$ is a special case of a form parameter for the ring $\Z$ with the trivial involution, as defined by 
Schlichting \cite[Definition 3.3]{S}.
In Schlichting's notation $(I, \sigma) = (\Z, \epsilon_Q), \Lambda = Q, \tau = \p,
\rho = \h$, and we discuss the quadratic action of $\Z$ on $Q$ in Section~\ref{ss:Z-grp}.
\end{remark}

We gave a list of standard form parameters in \eqref{eq:ind+} and \eqref{eq:ind-} and we give a classification of all form parameters in Section \ref{s:FPclass}.  Next we mention some examples of form parameters from topology and algebraic surgery.

\begin{example}
%%%%%%%%
Let $SO_q$ be the special orthogonal group on $\R^q$, so that $\pi_{q-1}(SO_q)$ classifies isomorphism classes of oriented rank $q$ vector bundles over $S^q$. For $\xi \in \pi_{q-1}(SO_q)$, let $V_\xi \to S^q$ be a corresponding oriented vector bundle and let $\tau_{q} \in \pi_{q-1}(SO_q)$ denote the homotopy class of the clutching function of the  tangent bundle of the $q$-sphere. We fix an orientation on $S^q$ and following Baues \cite[(8.6)]{Bau2} we define the 
form parameter
\[ \pi_{q-1}\{SO_q\} := \bigl( \pi_{q-1}(SO_q), {\rm h}, {\rm p} \bigr) , \]
where ${\rm h}(\xi) = e(V_\xi) \in \Z$ is the Euler number of $V_\xi$ and ${\rm p}(a) = a \tau_{q}$.  Elementary properties of the Euler class and $\tau_q$ show that $\pi_{q-1}\{SO_q\}$ is a 
form parameter with symmetry $(-1)^q$.
\end{example}

\begin{example}
%%%%%%%
A chain bundle over  the integers, $(B_*, \beta)$ (see \cite[\S 10]{R1}), defines for each integer $q$ a 
form parameter $Q_\beta(q)$, with symmetry $(-1)^q$.
\end{example}

%%%%%%%%

\begin{defin}
A morphism of 
form parameters $\alpha \colon (P_e, {\rm h}_P, {\rm p}_P)  \to (Q_e, {\rm h}_Q, {\rm p}_Q)$ is a homomorphism $\alpha \colon P_e \to Q_e$  such that ${\rm h}_R \circ \alpha = {\rm h}_Q$ and $\alpha \circ {\rm p}_Q = {\rm p}_R$. 
Form parameters form a category, which we denote by $\FP$,
and we denote the full subcategories of symmetric and anti-symmetric form parameters will be denoted 
by $\FP_+$ and $\FP_-$ respectively.
\end{defin}

\begin{remark}
Since morphisms preserve the symmetry of 
form parameters, $\FP$ is the disjoint union of $\FP_+$ and $\FP_-$. 
\end{remark}

\subsubsection{Linearisation and quasi-Wu classes}
We next define some isomorphism invariants of form parameters.

\begin{definition} \label{def:lin}
Let $Q = (Q_e, {\rm h}, {\rm p})$ be a 
form parameter. The {\em linearisation of $Q$} is the abelian group $SQ := Q_e/\mathrm{Im}(\mathrm{p})$. The quotient map will be denoted $\pi \colon Q_e \to SQ$.
\end{definition}

\begin{defin} \label{d:q-Wu}
\quad
\begin{compactenum}[a)]
\item For a symmetric 
form parameter $Q$, let $v_Q \colon SQ \to \Z_2$ be the homomorphism induced by $\h \colon Q_e \to \Z$.
The homomorphism $v_Q$ is well-defined since ${\rm hp} = 2 {\rm Id}_{\Z}$.
\item If $Q$ is an anti-symmetric 
form parameter, then ${\rm p} \colon \Z \to Q_e$ factors though the quotient homomorphism 
$\Z \to \Z_2$ and we let $v'_Q \colon \Z_2 \to Q_e$ denote the homomorphism induced by ${\rm p}$.
\end{compactenum}
The homomorphism $v_Q$ or $v'_Q$ is called the {\em quasi-Wu class} of $Q$.
\end{defin}

If $Q = (Q_e, {\rm h}, {\rm p})$ is a symmetric 
form parameter, then ${\rm h} \circ {\rm p} = 2\Id_{\Z}$ (hence ${\rm p}$ is injective), so we immediately get the following 

\begin{lem} \label{l:Wu-pb}
If $Q$ is a symmetric 
form parameter, then there is a commutative diagram of short exact sequences
\[
\xymatrix{
0 \ar[r] & \Z \ar[r]^-{{\rm p}} \ar@{=}[d] & Q_e \ar[r]^-{\pi} \ar[d]^-{{\rm h}} & SQ \ar[d]^-{v_Q} \ar[r] & 0 \\
0 \ar[r] & \Z \ar[r]^-{2} & \Z \ar[r]^-{\varrho_2} & \Z_2 \ar[r] & 0 , 
}
\] 
where $\varrho_2$ denotes reduction mod $2$. In particular the 
right hand square of this diagram is a pullback square.
\qed
\end{lem}

\subsubsection{Splittings of form parameters and their morphisms}  
%%%%%%%%%%%%%%%%%%%%%%%%%%%%%%%%%%%%%%%%%%%%%%%%%%%%%%%%%%%%
Recall from \eqref{eq:split} that given a form parameter $Q = (Q_e, {\rm h}, {\rm p})$ and abelian group $G$, we write $Q \oplus G$ for the {\em split} form parameter $(Q_e \oplus G, {\rm h} \oplus 0, {\rm p} \oplus 0)$.

\begin{definition} \label{d:maxd}
%%%%%%%%%%%%%%%%
Let $P$ be a form parameter.
\begin{compactenum}[a)]
\item
A {\em splitting} of $P$ is an isomorphism $P \cong Q \oplus G$ for some form parameter $Q$ and abelian group $G$. Equivalently, it is an isomorphism $P_e \cong Q_e \oplus G$ for some abelian groups $Q_e$ and $G$ such that under this identification $\h \big| _{0 \oplus G} = 0$ and $\p(1) \in Q_e \oplus 0$.  
\item
We say that $P$ is {\em decomposable} if there is a splitting $P \cong Q \oplus G$ with $G \not\cong 0$
and {\em indecomposable} otherwise.
\item
A splitting $P \cong Q \oplus G$ of $P$ is called {\em maximal} if $Q$ is indecomposable.
\end{compactenum}
\end{definition}

\begin{definition} \label{d:split}
\quad
\begin{compactenum}[a)]
\item A morphism $\alpha \colon Q_1 \oplus G_1 \to Q_2 \oplus G_2$ between split form parameters is called {\em split} if $\alpha = \alpha_q \oplus \alpha_l \colon Q_{1e} \oplus G_1 \to Q_{2e} \oplus G_2$ for some homomorphisms $\alpha_q \colon Q_{1e} \to Q_{2e}$ and $\alpha_l \colon G_1 \to G_2$.
\item A morphism $\alpha \colon P_1 \to P_2$ is called {\em splittable} if there are splittings $\alpha_1 \colon P_1 \to Q_1 \oplus G_1$ and $\alpha_2 \colon P_2 \to Q_2 \oplus G_2$ such that the composition $\alpha_2 \circ \alpha \circ \alpha^{-1}_1 \colon Q_1 \oplus G_1 \to Q_2 \oplus G_2$ is split.
\end{compactenum}
\end{definition}

The following example shows that not all morphisms of form parameters are splittable (in a non-trivial way).

\begin{example} \label{ex:non-splittable}
%%%%%%%%%%%%%%%%%%%%%%%%%%%%%%%%%%%%%
Let $\alpha \colon \Z^P_\infty \to Q^+ \oplus \Z_3$ be the morphism of form parameters given by $\left(\begin{smallmatrix} 1 & 0 \\ 1 & 2 \end{smallmatrix}\right) \colon \Z \oplus \Z \rightarrow \Z \oplus \Z_3$ (so that $S\alpha = \left(\begin{smallmatrix} 1 \\ 1 \end{smallmatrix}\right) \colon \Z \to \Z_2 \oplus \Z_3$). This homomorphism is surjective, hence $\alpha$ is not splittable.

Let $\alpha \colon \Z^\Lambda_2 \to \Z^\Lambda_3 \oplus \Z_2$ be the morphism of form parameters given by 
$\alpha \colon \Z_4 \to \Z_8 \oplus \Z_2, [a] \mapsto ([2a], [a])$. 
Then $\alpha$ is not splittable since $\alpha(\Z_4)$ is not contained in any proper summand of $\Z_8 \oplus \Z_2$.
\end{example}

\subsubsection{Extended symmetrisation and the
extended quadratic lift}  
%%%%%%%%%%%%%%%%%%%%%%%%%%%%%%%%%%%%%%%%%%%%%%%%%%%%%%%%%%%%
We introduce the morphisms which will play a key role in the natural description of the Witt group functor in Section~\ref{s:WGN}.

\begin{definition} \label{d:es}
%%%%%%%%%%%%%%%%%%
Let $Q = (Q_e, {\rm h}, {\rm p})$ be a symmetric 
form parameter. The {\em extended symmetrisation morphism}, 
$\es_Q \colon Q \rightarrow Q^+ \oplus SQ$, is given by 
$(\h,\pi) \colon Q_e \rightarrow \Z \oplus SQ$.
Unless it is required, we will drop the subscript ``$Q$" from the notation, writing ``$\es$" in place of ``$\es_Q$".
\end{definition}

Note that by Lemma \ref{l:Wu-pb} $(\h,\pi)$ determines an isomorphism between $Q_e$ and the pullback of the diagram 
\[
\xymatrix{
 & SQ \ar[d]^-{v_Q} \\
\Z \ar[r]^-{\varrho_2} & \Z_2 
}
\] 
that is, $\{ (n,x) \in \Z \times SQ \mid \varrho_2(n) = v_Q(x) \}$.

\begin{definition} \label{d:eql}
%%%%%%%%%%%%%%%%
Let $Q = (Q_e, {\rm h}, {\rm p})$ be an anti-symmetric form parameter. 
The {\em extended quadratic lift} is the morphism 
$\eql_Q \colon Q_- \oplus Q_e \rightarrow Q$, which is given by 
$v'_Q \oplus \Id \colon \Z/2 \oplus Q_e \rightarrow Q_e$.
Unless it is required, we will drop the subscript ``$Q$" from the notation,
writing ``$\eql$" in place of ``$\eql_Q$".
\end{definition}

\subsubsection{$\Z$-groups, quadratic $\Z$-modules and Schlichting's definition of a form parameter} \label{ss:Z-grp}
%%%%%%%%%%%%%%%%%%%%%%%%%%%%%%%%%%%%%%%%%%
A form parameter defines a ``quadratic action" of $\Z$ on $Q_e$, which is both a particular case of an ``$R$-group action'' as defined by Ranicki \cite[\S 4]{R1}, and gives $Q_e$ the structure of a quadratic $\Z$-module, a notion which is an essential part of Schlichting's definition of a form parameter \cite[Definition 3.3]{S}.

\begin{definition}
%%%%%%%%%
A {\em $\Z$-group} $U$ is an abelian group together with an action
\[ \Z \times U \to U, \quad (a,u) \to a \cdot u,\]
such that for all $a, b \in \Z$ and $u, v  \in \Z$
\[ a \cdot (u+v) =  a \cdot u + a \cdot v, \qquad (ab) \cdot u = a \cdot (b \cdot u),
\qquad 1 \cdot u = u. \]
A {\em linear $\Z$-group}, or $\Z$-module, i.e.\ abelian group, is a $\Z$-group such that 
\[(a+b) \cdot u = a\cdot u + b\cdot u.\]
In the future we will write $au$ for the linear action of $a \in \Z$ on an abelian group $U \ni u$.
A map of $\Z$-groups is a homomorphism of the underlying abelian groups which commutes with the 
$\Z$-actions.
\end{definition}

A form parameter $Q$ gives $Q_e$ the structure of a $\Z$-group (c.f.\ \cite[Remark 2.3]{S}), where $1 \in \Z$ acts by the identity and the action of $a \in \Z$ is inductively defined by 
$$(a + b) \cdot q := a \cdot q + b \cdot q + {\rm p}(a{\rm h}(q) b) = a \cdot q + b\cdot q + ab {\rm h}(q) {\rm p}(1),$$
where $a, b \in \Z$. 
If $C^a_2 := \frac{a(a-1)}{2}$ for $a \in \Z$, then induction shows that for all $a \in \Z$ and $q \in Q_e$
\begin{equation} \label{eq:Zactn}
a \cdot q = a q + C^a_2{\rm h}(q){\rm p}(1).
\end{equation}
It follows that a form parameter $Q$ is equivalent to a form parameter $((I, \sigma), \Lambda)$, as defined, by Schlichting \cite[Definition 3.3]{S}, under the following identifications.  
The ambient ring with with involution is $\Z$ with the trivial involution, the $\Z$-bimodule with involution is $(I, \sigma) = (\Z, \epsilon_Q)$,
the abelian group is $\Lambda = Q_e$, $\tau = \p \colon \Z \to Q_e$, $\rho = \h \colon Q_e \to \Z$
and the multiplicative $\Z$-linear left action $Q$ on $\Lambda$ is given by \eqref{eq:Zactn}.
%
%%%%%%%%%%%%%%%%%%%%%%%%%%%%%%%%%%%%%%%%%%%%%%%%%%%%

\subsection{Extended quadratic forms over $Q$} \label{ss:Q_forms} 
%%%%%%%%%%%%%%%%%%%%%%%%%%%%%%%%%%%%%%%%%%%%%%%%%%%%%%%%%%%%
In this subsection we define extended quadratic forms over a 
form parameter $Q$.  For brevity, we call such forms {\em $Q$-forms.} 
A $Q$-form is an enrichment of a symmetric or anti-symmetric bilinear form and we first review these notions.

Let $X$ be a finitely generated abelian group.
We define the {\em dual} of $X$ to be the group  $X^* := {\rm Hom}(X, \Z)$.
The canonical map $X \to X^{**}$ is defined by $x \mapsto (\phi \mapsto \phi(x))$ for 
$x \in X$ and $\phi \in X^*$.  It is surjective with kernel $TX$, the torsion subgroup of $X$.
We will be largely interested in the torsion free case, i.e.\ when $TX = 0$.

\begin{definition}
%%%%%%%%%
For $\epsilon = \pm 1$,
an {\em $\epsilon$-symmetric bilinear form} 
$(X, \lambda)$ 
consists of an abelian group $X$ and a pairing   
$$\lambda\colon X \times X \to \Z,$$
which for $x, x', y \in X$ and $a,b \in \Z$, satisfies
\begin{equation} \label{eq:1}
\lambda(ax,by) = b\lambda(x,y){a}, \quad \lambda(x+x',y) = \lambda(x,y) 
+\lambda(x',y) \quad \text{and} \quad \lambda(x,y) = \epsilon {\lambda(y,x)}.
\end{equation}
\noindent
Equivalently, 
$\lambda$ defines and is defined by its {\em adjoint homomorphism},
\[\wh \lambda\colon X \to X^*, \quad x \mapsto \bigl( y \mapsto \lambda(x, y) \bigr).\]
When $X$ is torsion free, 
the conditions in \eqref{eq:1} are equivalent to 
$\wh \lambda$ being a homomorphism
such that $\wh \lambda^* = \epsilon \wh \lambda$,
where $\wh \lambda^* \colon X^{**} = X \to X^*$.
If $\wh \lambda$ is an isomorphism then $(X, \lambda)$ 
is called {\em nonsingular}.
\end{definition}
%%%%%%%%%%%%%%%%%%%%%%%%%%%%%%%%%%%%%%%%%%%%%%%%%%%%

\begin{definition} \label{def:eqf}
%%%%%%%%%
Let $Q = (Q_e, {\rm h}, {\rm p})$ be a 
form parameter with symmetry $\epsilon_Q$.  
{\em An extended quadratic form with values in $Q_{}$,} 
or simply a {\em $Q$-form}, is a triple $(X, \lambda, \mu)$ consisting of 
a bilinear form $(X, \lambda)$ and a function $\mu\colon X \to Q_e$ such that for all $x, y \in X$
\[ \mu(x + y) = \mu(x) + \mu(y) + {\rm p}(\lambda(x, y))  
\qquad \text{and} \qquad
\lambda(x, x) = {\rm h}(\mu (x)). \]
Note that these conditions imply that $\lambda$ is symmetric if $\epsilon_Q=1$, 
or anti-symmetric if $\epsilon_Q=-1$.
A $Q$-form $(X, \lambda, \mu)$ is {\em nonsingular} when the underlying form $(X, \lambda)$ is nonsingular.
\end{definition}

\begin{remark}
%%%%%%%%
For a $Q$-form $(X, \lambda, \mu)$ we have $\mu(ax) = a \cdot \mu(x) = a \mu(x) + C^a_2 \lambda(x, x){\rm p}(1)$ for all $a \in \Z$ and $x \in X$.  This is a consequence of
the discussion in Section~\ref{ss:Z-grp} above.
\end{remark}

\begin{remark} \label{rem:eqf-not}
%%%%%%%%%%%%%%%%%
In computations an extended quadratic form $(X, \lambda, \mu)$ will be written as 
\[
\left(
X, 
M_{\lambda},
\left( \begin{array}{c} \mu(x_1) \\ \vdots \\ \mu(x_k) \end{array} \right) 
\right),
\]
where $x_1, \ldots, x_k$ is a basis of $X$ and $M_{\lambda}$ is the matrix of $\lambda$ in this basis. Note that, as $\mu$ is not linear, the vector in the third component should not be interpreted as a homomorphism $X \rightarrow Q_e$ (but the value of $\mu$ on an arbitrary element of $X$ can be unambiguously computed from this data using the properties of $\mu$). Conversely, a triple 
\[
\left(
X, 
M,
\left( \begin{array}{c} q_1 \\ \vdots \\ q_k \end{array} \right) 
\right),
\]
where $X$ is a free abelian group with a fixed basis $x_1, \ldots, x_k$, $M \in \Z^{k \times k}$ and $q_1, \ldots, q_k \in Q_e$, denotes a $Q$-form on $X$ if and only if $M$ is $\epsilon_Q$-symmetric and $M_{ii} = {\rm h}(q_i)$ for every $i$.

If $X = 0$ is the trivial group, we write $\ul 0$ for the unique $Q$-form defined on $0$.
\end{remark}

\begin{example} \label{ex:hyp}
%%%%%%%%
For a natural number $m$,
the standard rank-$2m$ hyperbolic form over an $\epsilon$-symmetric 
form parameter $Q$ is 
\[
H_{\epsilon}(\Z^m) = 
\left(
\Z^m \oplus \Z^m, 
\left[ \begin{array}{cc} 0 & I_m \\ \epsilon I_m & 0 \end{array} \right], 
\left( \begin{array}{c} 0 \\ \vdots \\ 0 \end{array} \right) 
\right)
\]
where $I_m$ denotes the $m \times m$ identity matrix.
If $m = 0$, then $H_\epsilon(0) = \ul 0$.
\end{example}

We briefly review some essential concepts and terminology related to $Q$-forms.  
When it is convenient we shall denote the $Q$-form $(X, \lambda, \mu)$ simply by $\ul \mu$.  
The rank of $\ul \mu$, denoted by $\rank(\ul \mu)$, is the rank of the abelian group $X$. If $Q$ is symmetric, then the signature of $\ul \mu$, $\sigma(\ul \mu)$, is the signature of the symmetric bilinear form $(X,\lambda)$.

\begin{defin} \label{d:HomZ}
%%%%%%%%%%%%%
Following \cite[\S 5]{Bau1},
we denote the set of $Q$-forms on $X$
by $\HomZ(X,Q)$ 
and the subset of nonsingular $Q$-forms by $\HomZns(X, Q) \subset \HomZ(X, Q)$.
We recall from Section~\ref{ss:GWG} that $\HomZns(Q)$ denotes the set of isomorphism classes of all nonsingular $Q$-forms.
\end{defin}

Elementwise addition (given by $(\lambda+\lambda')(x,y) = \lambda(x,y) + \lambda'(x,y)$ and $(\mu+\mu')(x) = \mu(x) + \mu'(x)$) gives $\HomZ(X,Q)$ the structure of an abelian group with zero element the zero form, $(X, 0, 0)$, and inverse defined by $-(X, \lambda, \mu) := (X, -\lambda, -\mu)$, where $(-\lambda)(x,y) = -\lambda(x,y)$ and $(-\mu)(x) = -\mu(x)$.  

\begin{defin}
\quad
\begin{compactenum}[a)]
\item Given a $Q$-form $(Y, \lambda, \mu)$ on $Y$ and homomorphism $f \colon X  \to Y$, the {\em pullback} of $(Y, \lambda, \mu)$ along $f$ is the $Q$-form $(X, f^*\lambda, f^*\mu) \in \HomZ(X, Q)$, where $f^*\lambda(x,y) = \lambda(f(x),f(y))$ and $f^*\mu(x) = \mu(f(x))$. 
\item If $(X, \lambda, \mu) \in \HomZ(X, Q)$, $(Y, \phi, \psi) \in \HomZ(Y, Q)$ and $f\colon X \to Y$ is a homomorphism then $f$ defines a {\em morphism} of $Q$-forms if $(X, f^*\phi, f^*\psi) = (X, \lambda, \mu)$.
\item A morphism of $Q$-forms is called an {\em isometry} if it is an isomorphism of the underlying abelian groups. 
\end{compactenum} 
\end{defin}

\begin{remark}
The assignment $X \mapsto \HomZ(X, Q)$ defines a contravariant {\em quadratic functor} 
$$ \mcal{F}_Q \colon \AbF \to \AbF,$$
in the sense of Baues \cite[\S 1]{Bau1}.
\end{remark}

Similarly, for a fixed $X$ the assignment $Q \mapsto \HomZ(X, Q)$ defines a covariant functor $\FP \rightarrow \AbF$. 
Specifically, for a morphism $\alpha \colon Q \rightarrow Q'$ for form parameters, we have the induced the homomorphism
\begin{equation*} \label{eq:HomZalpha}
\HomZ(X, \alpha) \colon \HomZ(X, Q) \to \HomZ(X, Q'),
\quad (X, \lambda, \mu) \mapsto (X, \lambda, \alpha \circ \mu).
\end{equation*}

\begin{defin}
%%%%%%
The direct sum of the $Q$-forms $(X, \lambda, \mu)$ and $(X', \lambda', \mu')$ is the $Q$-form $(X \oplus X', \lambda \oplus \lambda', \mu \oplus \mu')$, where $(\lambda \oplus \lambda')((x,x'),(y,y')) = \lambda(x,y)+\lambda'(x',y')$ and $(\mu \oplus \mu')(x,x') = \mu(x)+\mu'(x')$. 
\end{defin}

\begin{definition} \label{def:full}
%%%%%%%%%%%%%%%%%
Let $\ul \mu = (X, \lambda, \mu)$ be a $Q$-form. 
\begin{compactenum}[a)]
\item The composition of $\mu \colon X \to Q_e$ with the quotient map $Q_e \to SQ$ is a homomorphism $S\mu \colon X \to SQ$, 
which is called the {\em linearisation of $\ul \mu$.} 
\item The $Q$-form $\ul \mu$ is called {\em full} if its linearisation $S\mu \colon X \to SQ$ is onto.
\end{compactenum}
\end{definition}

\begin{definition} \label{def:char}
%%%%%%%%%%%%%%%%%
Let $(X, \lambda)$ be a symmetric bilinear form, 
and $f \colon X \to C$ a homomorphism to a cyclic group
with mod~$2$ reduction $f_2 \colon X \to C/2C \subseteq \Z_2$.
Then $f$ is called {\em characteristic for $\lambda$}, if for all $x \in X$ we have
\[ \varrho_2 \circ \lambda(x, x) = f_2(x).\]
\end{definition}

We next discuss the homomorphism of groups of forms,
\begin{equation*}  \label{eq:HomZes}
\HomZ(X, \es) \colon \HomZ(X, Q) \to \HomZ(X, Q^+ \oplus SQ),
\end{equation*}
where $Q$ is a symmetric form parameter, $X$ is an abelian group
and $\es \colon Q \to Q^+ \oplus SQ$ is the extended symmetrisation morphism.
Note that if $(X,\lambda,(\mu_0,\mu_1))$ is a $(Q^+ \oplus G)$-form for some abelian group $G$, 
then $\mu_0(x) = \lambda(x,x)$ for every $x$ and $\mu_1$ is a homomorphism by Definition \ref{def:eqf}. 
Moreover, if a $(Q^+ \oplus G)$-form $(X,\lambda,(\mu_0,\mu_1))$ is the extended symmetrisation of a $Q$-form $(X, \lambda, \mu)$,
then $G = SQ$ and $\mu_1 = S\mu$.

\begin{definition}
%%%%%%%%%%%%%%%%%%%
Let $Q$ be a symmetric form parameter 
and define $\HomZ^{\rm ch}(X, Q^+ \oplus SQ) \subset \HomZ(X, Q^+ \oplus SQ)$ to be the subgroup of $(Q^+ \oplus SQ)$-forms
$(X, \lambda, (\mu_0, \mu_1))$ such that $v_Q \circ \mu_1$ is characteristic for $\lambda$. 
\end{definition}

\begin{lem} \label{l:gmtrc}
%%%%%%%%%%%%%
%
Let $Q$ be a symmetric form parameter and $X$ an abelian group. 

a) For every $Q$-form $(X, \lambda, \mu)$, the composition
$v_Q \circ S\mu \colon X \to \Z_2$ is characteristic for $\lambda$;
i.e.\  $\HomZ(X, \es)$ takes values in $\HomZ^{\rm ch}(X, Q^+ \oplus SQ)$. 

b)
The map $\HomZ(X, \es) \colon \HomZ(X, Q) \to \HomZ^{\rm ch}(X, Q^+ \oplus SQ)$
is an isomorphism.
\end{lem}

\begin{proof}
%%%%%%%%
a) Let $Q = (Q_e, {\rm h}, {\rm p})$.
Let $\pi \colon Q_e \rightarrow SQ$ denote the quotient map so that the extended symmetrisation map is given by $({\rm h},\pi)$. By Lemma \ref{l:Wu-pb} $({\rm h},\pi)$ is an isomorphism between $Q_e$ and the subgroup $\left\{ (z,y) \in \Z \oplus SQ \mid \varrho_2(z) = v_Q(y) \right\}$. This means that if $(X,\lambda,\mu) \in \HomZ(X,Q)$ and its image is $(X,\lambda,(\mu_0,\mu_1))$ (i.e.\ $(\mu_0,\mu_1) = ({\rm h},\pi) \circ \mu$), then $v_Q \circ \mu_1(x) = \varrho_2 \circ \mu_0(x) = \varrho_2 \circ \lambda(x,x)$.

b) Now suppose that $(X,\lambda,(\mu_0,\mu_1)) \in \HomZ(X,Q^+ \oplus SQ)$ is such that $\varrho_2 \circ \lambda(x,x) = v_Q \circ \mu_1(x)$ for every $x \in X$. By the above $(\mu_0,\mu_1)(x) \in \Im({\rm h},\pi)$ for every $x \in X$, and so there is a unique function $\mu \colon X \rightarrow Q_e$ such that $(\mu_0,\mu_1) = ({\rm h},\pi) \circ \mu$. We will prove that $(X,\lambda,\mu)$ is $Q$-form.
First, ${\rm h}(\mu(x)) = \mu_0(x) = \lambda(x,x)$. Second, 
\begin{align*}
\mu(x+y) & = ({\rm h},\pi)^{-1} \circ (\mu_0,\mu_1)(x+y) \\
& = ({\rm h},\pi)^{-1}(\lambda(x+y,x+y),\mu_1(x+y)) \\
& = ({\rm h},\pi)^{-1}((\lambda(x,x),\mu_1(x))+(\lambda(y,y),\mu_1(y))+(2\lambda(x,y),0)) \\
& = ({\rm h},\pi)^{-1}(\mu_0(x),\mu_1(x)) + ({\rm h},\pi)^{-1}(\mu_0(y),\mu_1(y)) + ({\rm h},\pi)^{-1}(2\lambda(x,y),0) \\
& = \mu(x) + \mu(y) + {\rm p}(\lambda(x,y)). 
\end{align*}
Therefore $(X,\lambda,\mu) \in \HomZ(X,Q)$. 
Since the assignment $(X,\lambda,(\mu_0,\mu_1)) \mapsto (X,\lambda,\mu)$ is the inverse of the map induced by the extended symmetrisation map, part b) follows.
\end{proof}

\begin{remark}
%%%%%%%%
We briefly discuss the relationship between the terminology and notation of this 
paper and those of 
\cite[Chapter 4]{N1} and \cite{N2}. 
In \cite[Chapter 4]{N1} and \cite{N2} 
the second author
introduces (symmetric) {\em extended quadratic forms over an abelian group $G$}. Such a form is a triple $(X, \lambda, f)$, where $(X, \lambda)$ is a symmetric bilinear form and $f \colon X  \to G$ is a homomorphism. It corresponds to the $(Q^+ \oplus G)$-form $(X, \lambda, (\mu_0, f))$ (where $\mu_0(x)=\lambda(x,x)$). Given a homomorphism $t \colon G \to \Z_2$, the form is called {\em geometric} with respect to $t$ if $t \circ f \colon X \to \Z_2$ is characteristic for $\lambda$. Thus for a symmetric form parameter $Q$, it follows from Lemma~\ref{l:gmtrc} that $Q$-forms are equivalent to extended quadratic forms over the abelian group $SQ$ that are geometric with respect to $v_Q$. 
\end{remark}

%%%%%%

\subsection{The Witt group $W_0(Q)$} \label{ss:Witt_def}
%%%%%%%%%%%%%%%%%%%%%%%%%%%%%%%%%%%%%%%%%%%%%%%%%%%
%
Let $Q$ be a form parameter with symmetry $\epsilon_Q$.
In this subsection, we assume that all $Q$-forms
$(X, \lambda, \mu)$ are nonsingular,
and recall that $\HomZns(Q)$ is the set of isomorphism classes of nonsingular $Q$-forms.
Of course, if $(X, \lambda, \mu)$ is nonsingular
then $X$ is torsion free.

A summand, $L \subset X$, is called a {\em lagrangian} of $(X, \lambda, \mu)$ if $\rank(X) = 2 \rank(L)$ and the inclusion $i\colon L \hra X$ satisfies $i^*(X, \lambda, \mu) = (L, 0, 0)$. 
For any $Q$-form $(X, \lambda, \mu)$, the $Q$-form $(X \oplus X, \lambda \oplus -\lambda, \mu \oplus -\mu)$ has the diagonal, 
\[ \Delta := \left\{ (x, x) \,\big|\, x \in X \right\} \subset X \oplus X, \]  
as a lagrangian. 
A $Q$-form which admits a lagrangian is called {\em metabolic}. We warn the reader that we depart from \cite[\S 10]{R1} where the term {\em hyperbolic} is used 
for metabolic $Q$-forms.

We say that $Q$-forms $\ul{\mu_0}$ and $\ul{\mu_1}$ are {\em Witt equivalent} if there are metabolic 
$Q$-forms $\ul{\phi_0}$ and $\ul{\phi_1}$ such that  
$\ul{\mu_0} \oplus \ul{\phi_0} \cong \ul{\mu_1} \oplus \ul{\phi_1}$, and write $\ul \mu_0 \sim_{\rm W} \ul \mu_1$.
We denote the Witt equivalence class of a $Q$-form $\ul \mu$ by $[\ul \mu]$, and recall the following definition
from the Introduction.

\begin{definition} \label{def:WG}
%%%%%%%%%%%%%%%%
The \textit{Witt group} $W_0(Q)$ is the abelian group formed by the set of Witt equivalence classes of nonsingular $Q$-forms, with addition given by direct sum:
\[ W_0(Q):= \HomZns(Q)/\!\! \sim_{\rm W} ~ =
\left( \bigl\{ [X, \lambda, \mu] \, \big| \,
\text{$(X, \lambda, \mu)$ is a nonsingular $Q$-form} \bigr\}, 
\oplus \right) \]
\noindent
The zero element of $W_0(Q)$ is the equivalence class of metabolic forms and
the inverse of $[\ul{\mu}]$ is $[-\ul{\mu}]$ since $\ul{\mu} \oplus (-\ul{\mu})$ is metabolic.
\end{definition}

\begin{remark} \label{r:zero-in-Witt}
By definition, a $Q$-form $\ul \mu$ represents the zero element of $W_0(Q)$ if and only if it is {\em stably metabolic}; i.e.\ there are metabolic $Q$-forms $\ul \mu_0$ and $\ul \mu_1$ 
such that $\ul \mu \oplus \ul \mu_0 \cong\ul \mu_1$.
In Remark~\ref{r:sm-nm}, we observe that for all $k \geq 2$, there are stably metabolic forms over $Q = \Z^\Lambda_k$, which are not metabolic, and also that there is a stably metabolic form over $Q = \Z^P_1$, which is not metabolic.
\end{remark}

Now suppose that $\alpha \colon Q_1 \to Q_2$ is a morphism of form parameters.
Then $\alpha$ induces the map
\[ \HomZns(X,\alpha) \colon \HomZns(X, Q_1) \to \HomZns(X, Q_2),
\quad [X, \lambda, \mu] \mapsto [X, \lambda, \alpha \circ \mu].
\]
Since $\HomZns(X, \alpha)$ preserves the zero form (on any subgroup of $X$), it also preserves metabolic forms,
and so $\alpha$ descends to define a homomorphism of Witt groups.

\begin{definition} \label{def:WG(alpha)}
Let $\alpha \colon Q_0 \to Q_1$ be a morphism of form parameters.  Define 
\[ W_0(\alpha) \colon W_0(Q_1) \to W_0(Q_2), \quad
[X, \lambda, \mu] \mapsto [X, \lambda, \alpha \circ \mu]. \]
\end{definition}

\noindent
Hence taking the Witt group of a form parameter defines a functor, 
\[ W_0 \colon \FP \to \AbF,
\quad Q \mapsto W_0(Q), \]
from the category of form parameters to the category of abelian groups.  As we stated in the Introduction, the main aim of this paper is to compute the Witt group functor $ W_0 \colon \FP \to \AbF$.
%
%
%%%%%%%%%%%%%%%%%%%%%%%%%%%%%%%%%%%%%%%%%%%%%%%%%%%%%%%%%%%%

\section{The classification of quadratic form parameters} \label{s:FPclass}
%%%%%%%%%%%%%%%%%%%%%%%%%%%%%%%%%%%%%%%%%%%%%%%%%%%%
In this section we classify the category oriented form parameters over $\Z$.
In particular, 
we prove that every form parameter has a maximal splitting, that is, it is isomorphic to the sum of an indecomposable form parameter and an abelian group (see Definition \ref{d:maxd}). We also prove that \eqref{eq:ind+} and \eqref{eq:ind-} give a complete list of (isomorphism classes of) indecomposable form parameters. Finally, we will explicitly describe the possible morphisms between indecomposable form parameters.

We begin by recalling the list of standard form parameters from \eqref{eq:ind+} and \eqref{eq:ind-}. 

\begin{defin} \label{def:indec}
%%%%%%%%%%%%%%%%%%%%%%%%%%%%%%%
We have the following standard form parameters (where $\Z^P$ was previously denoted $\Z^P_\infty$): 
\[
\begin{aligned}
&  \ul{\quad \epsilon = +1 \quad} & & \ul{\quad \epsilon = -1 \quad} \\
Q_{+} &= \Bigl( \Z \xra{2} \Z \xra{1} \Z \Bigr) & Q_{-} &= \Bigl( \Z_{2} \xra{0} \Z \xra{1} \Z_{2} \Bigr) \\
\Z^P &= \Bigl( \Z \oplus \Z \xra{(1,0)} \Z \xra{\left(\!\begin{smallmatrix} 2  \\ -1 \end{smallmatrix}\!\right)} \Z \oplus \Z \Bigr)
&  \\
\Z^P_k &= \Bigl(\Z \oplus \Z_{2^k} \xra{(1,0)} \Z \xra{\left(\!\begin{smallmatrix} 2  \\ -1 \end{smallmatrix}\!\right)} \Z \oplus \Z_{2^k} \Bigr) \quad (k \geq 1) \quad \quad  & \Z^{\Lambda}_k &= \Bigl( \Z_{2^k} \xra{~~0~~} \Z \xra{2^{k-1}} \Z_{2^{k}} \Bigr) \quad (k \geq 2) \\
Q^{+} &= \Bigl( \Z \xra{1} \Z \xra{2} \Z \Bigr) 
& Q^{-} &= \Bigl( 0 \xra{} \Z \xra{} 0 \Bigr)
\end{aligned}
\]
\end{defin}

\begin{remark}
%%%%%%%%
We warn the reader that our notation differs somewhat from Baues' notation \cite{Bau1,Bau2}, where $Q_{+}$ is called 
$\Z^S$, $Q^{+}$ is called $\Z^{\Gamma}$, $Q_{-}$ is called $\Z^{\Lambda}_1$ and $Q^{-}$ is called
$\Z^{\Lambda}$.  Baues' notation is based on the quadratic functors corresponding to
these quadratic modules; 
see Example~\ref{ex:qtp-funct}.  In particular the quadratic functors corresponding to $\Z^S, \Z^P$ and $\Z^\Lambda_1$ are respectively the symmetric tensor product,
the universal quadratic functor and the skewsymmetric tensor product.
Where there is an overlap in notation, our notation is designed to fit with Ranicki's notation in algebraic surgery.
\end{remark}

Recall that given a form parameter $Q = (Q_e, {\rm h}, {\rm p})$ and an abelian group $G$, 
the split form parameter $Q \oplus G$ is $(Q_e \oplus G, {\rm h} \oplus 0, {\rm p} \oplus 0)$, and
that a splitting of $P$ is an isomorphism $P \cong Q \oplus G$ for some form parameter $Q$ and abelian group $G$.   
A form parameter $P$ is decomposable if there is a splitting $P \cong Q \oplus G$ with $G \not\cong 0$
and indecomposable otherwise, 
a splitting $P \cong Q \oplus G$ of $P$ is called maximal if $Q$ is indecomposable.
The following result restates Theorem~\ref{t:qfp-class-intro} from the introduction.

\begin{thm}[The classification of quadratic form parameters] \label{t:qfp-class}
\quad

a) The form parameters in Definition \ref{def:indec} are indecomposable and pairwise non-isomorphic. Any indecomposable form parameter is isomorphic to one of them. 

b) Every form parameter $P$ has a maximal splitting $P \cong Q \oplus G$. Moreover, $Q$ and $G$ are each well-defined up to isomorphism. 
\end{thm}

Following Theorem~\ref{t:qfp-class} we can define the following invariants of a form parameter (cf.\ Theorems \ref{thm:asymm-isom} and \ref{thm:symm-equiv}).

\begin{defin} \label{d:hght}
Let $P = Q \oplus G$ be a maximal splitting of a form parameter.
\begin{compactenum}[a)]
\item The {\em height} of an indecomposable form parameter,
$\hght(Q) \in \N_0 \cup 
\{\infty\}$ is defined by
\[  \hght(Q_+) = \hght(Q^-) =0, \qquad
\hght(Q^+) = \hght(Q_-) = 1, \qquad
\hght(\Z^P_k) = \hght(\Z^\Lambda_{k+1}) = k{+}1
\quad \text{and} \quad
\hght(\Z^P) = \infty.\]
\item
The height of $P$ is defined by 
$\hght(P) = \hght(Q)$.
\item 
The {\em complement} of $P$, $[G_P]$, is defined to be the isomorphism class of $G$.
\end{compactenum}
\end{defin}

As a direct consequence of 
Theorem~\ref{t:qfp-class} we have

\begin{corollary} \label{c:qfp-class}
Form parameters are classified up to isomorphism by their symmetry, height and complement;
i.e.\ $P_1 \cong P_2$ if and only if
$(\epsilon_{P_1}, \hght(P_1), [G_{P_1}]) = 
(\epsilon_{P_2}, \hght(P_2), [G_{P_2}])$. \qed
\end{corollary}

To prove Theorem~\ref{t:qfp-class}, we will identify $\FP_+$ and $\FP_-$ with categories of simpler objects. We start by introducing these categories in Section \ref{ss:slice-coslice}, and then describe the identifications in Section \ref{ss:cat}. The proof of Theorem \ref{t:qfp-class} 
is given at the end of Section \ref{ss:cat}.

\subsection{The slice and coslice categories $\AbF/\Z_2$ and $\Z_2/\AbF$} \label{ss:slice-coslice}
%%%%%%%%%%%%%%%%%%%%%%%%%%%%%%%%%%%%%%%%%%%%%%%%%%%%%%%
In this subsection we consider the target categories in the equivalences of Theorem~\ref{t:FP_class}. We will prove the analogue of Theorem \ref{t:qfp-class} for these categories, see Theorem \ref{t:(co)slice-d}.

\begin{defin}
We recall the following notation from the Introduction: 
\begin{compactenum}[a)]
\item $\AbF$ is the category of finitely generated abelian groups.
\item $\AbF / \Z_2$ denotes the slice category whose objects are homomorphisms $v \colon A \rightarrow \Z_2$ for finitely generated abelian groups $A$, and whose morphisms from $v_1 \colon A_1 \rightarrow \Z_2$ to $v_2 \colon A_2 \rightarrow \Z_2$ are homomorphisms $f \colon A_1 \to A_2$ such that 
the following diagram commutes:
\[
\xymatrix{
A_1 \ar[dr]_{v_1} \ar[rr]^f & & 
A_2 \ar[dl]^{v_2}\\
& \Z_2
}
\]
\item $\Z_2 / \AbF$ denotes the coslice category whose objects are homomorphisms $v' \colon \Z_2 \rightarrow A$ for finitely generated abelian groups $A$, and whose morphisms from $v'_1 \colon \Z_2 \rightarrow A_1$ to $v'_2 \colon \Z_2 \rightarrow A_2$ are homomorphisms $f \colon A_1 \to A_2$ such that 
the following diagram commutes:
\[
\xymatrix{
& \Z_2 \ar[dl]_{v'_1} \ar[dr]^{v'_2} \\
A_1 \ar[rr]^f & & 
A_2 }
\]
\end{compactenum}
\end{defin}

\begin{definition}
Given a homomorphism $v_0 \colon B \rightarrow \Z_2$ and an abelian group $C$ we write $v_0 \oplus C \colon B \oplus C \rightarrow \Z_2$ for the composition of the projection $B \oplus C \rightarrow B$ and $v_0$. 

Given a homomorphism $v'_0 \colon \Z_2 \rightarrow B$ and an abelian group $C$ we write $v'_0 \oplus C \colon \Z_2 \rightarrow B \oplus C$ for the composition of $v'_0$ and the inclusion $B \rightarrow B \oplus C$.
\end{definition}

\begin{definition}
Let $v \colon A \rightarrow \Z_2$ (resp.\ $v' \colon \Z_2 \rightarrow A$) be a homomorphism.
\begin{compactenum}[a)]
\item
A {\em splitting} of $v$ (resp.\ $v'$) is an isomorphism $v \cong v_0 \oplus C$ (resp.\ $v' \cong v'_0 \oplus C$) for some homomorphism $v_0 \colon B \rightarrow \Z_2$ (resp.\ $v'_0 \colon \Z_2 \rightarrow B$) and abelian group $C$.  
\item
We say that $v$ (resp.\ $v'$) is {\em decomposable} if there is a splitting $v \cong v_0 \oplus C$ (resp.\ $v' \cong v'_0 \oplus C$) with $C \not\cong 0$
and {\em indecomposable} otherwise.
\item
A splitting $v \cong v_0 \oplus C$ (resp.\ $v' \cong v'_0 \oplus C$) is called {\em maximal} if $v_0$ (resp.\ $v'_0$) is indecomposable.
\end{compactenum}
\end{definition}

\begin{example} \label{ex:slice-coslice}
We have the following homomorphisms, as objects in $\AbF / \Z_2$ and $\Z_2 / \AbF$: 
\[
\begin{aligned}
& \ul{\quad \AbF / \Z_2 \quad} & &  \ul{\quad \Z_2 / \AbF \quad} \\
0 & \colon 0 \rightarrow \Z_2 &  \iota_1 & \colon \Z_2 \rightarrow \Z_2 \\
1_{\infty} & \colon \Z \xra{1} \Z_2 \quad\quad\quad\quad & & & \\
1_{k+1} & \colon \Z_{2^{k+1}} \xra{1} \Z_2 \quad (k \geq 1) \quad &
\iota_k & \colon \Z_2 \xra{2^{k-1}} \Z_{2^k}
\quad (k \geq 2)\\
1_1 & \colon \Z_2 \xra{1} \Z_2 & 0 & \colon \Z_2 \xra{} 0 
\end{aligned}
\]
In particular $1_1 = \iota_1 = \Id_{\Z_2}$.
\end{example}

\begin{prop}
\label{p:slice-coslice-decomp}
a) The homomorphisms listed in Example \ref{ex:slice-coslice} are indecomposable and pairwise non-isomorphic.

b) Every homomorphism $v \colon A \rightarrow \Z_2$ (resp.\ $v' \colon \Z_2 \rightarrow A$) is isomorphic to some $v_0 \oplus C \colon B \oplus C \rightarrow \Z_2$ (resp.\ $v'_0 \oplus C \colon \Z_2 \rightarrow B \oplus C$), where $v_0 \colon B \rightarrow \Z_2$ (resp.\ $v'_0 \colon \Z_2 \rightarrow B$) is one of the homomorphisms listed in Example \ref{ex:slice-coslice} and $C$ is an abelian group.  

c) Every indecomposable homomorphism is isomorphic to one of the homomorphisms listed in Example \ref{ex:slice-coslice}.

d) If $v \colon A \rightarrow \Z_2$ (resp.\ $v' \colon \Z_2 \rightarrow A$) is isomorphic to $v_0 \oplus C \colon B \oplus C \rightarrow \Z_2$ (resp.\ $v'_0 \oplus C \colon \Z_2 \rightarrow B \oplus C$) for an indecomposable $v_0 \colon B \rightarrow \Z_2$ (resp.\ $v'_0 \colon \Z_2 \rightarrow B$), then the isomorphism classes of $v_0$ (resp.\ $v'_0$) and $C$ are determined by $v$ (resp.\ $v'$).
\end{prop}

\begin{proof}
a) This is straightforward to check.

b) First we consider the slice category. If $v = 0$, then $v \cong 0 \oplus A$. If $v \big| _{\Tor(A)} \neq 0$, then by Lemma \ref{lem:hom-decomp} $A \cong \Z_{2^k} \oplus H$ for some $k \geq 1$ and $H \leq \Ker(v)$, and $v \cong 1_k \oplus H \colon \Z_{2^k} \oplus H \rightarrow \Z_2$. Finally, if $v \neq 0$ but $v \big| _{\Tor(A)} = 0$, then $v$ factors through a non-zero map $A / \Tor(A) \rightarrow \Z_2$, and by Lemma \ref{lem:free-decomp} that map is isomorphic to $1_{\infty} \oplus H \colon \Z \oplus H \rightarrow \Z_2$ for some $H \leq A / \Tor(A)$. Therefore $v \cong 1_{\infty} \oplus (H \oplus \Tor(A))$. 

Next we consider the coslice category. If $v' = 0$, then $v' \cong 0 \oplus A$.  If $v' \neq 0$, then $v'(1) \in A$ is an element of order $2$. Let $l = \max \left\{ m \mid \exists x \in A : v'(1) = 2^m x \right\} \geq 0$ and choose an $x$ such that $v'(1) = 2^l x$. By Lemma \ref{lem:dir-summand} the subgroup $(x) \subseteq A$ is a direct summand, so $A \cong (x) \oplus (A / (x))$, and hence $v' \cong \iota_{l+1} \oplus (A / (x))$.

c) Follows by applying Part b) to the indecomposable homomorphism. 

d) Suppose that $v \cong v_0 \oplus C$ for a $v_0$ listed in Example \ref{ex:slice-coslice} and an abelian group $C$. If $v_0 = 0$, then $v = 0$. If $v_0 = 1_{\infty}$, then $v \neq 0$, but $\Tor(A) \leq \Ker(v)$. If $v_0 = 1_k$, then $\min \left\{ \text{order of $x$} \mid x \in \Tor(A), v(x)=1 \right\} = 2^k$. 

Similarly suppose that $v' \cong v'_0 \oplus C$ for a $v'_0$ listed in Example \ref{ex:slice-coslice} and an abelian group $C$. If $v'_0 = 0$, then $v' = 0$. If $v'_0 = \iota_k$, then $\max \left\{ m \mid \exists x \in A : v'(1) = 2^m x \right\} = k-1$. 

The above observations imply that for a fixed $v \colon A \rightarrow \Z_2$ (resp.\ $v' \colon \Z_2 \rightarrow A$) there is at most one homomorphism $v_0$ (resp.\ $v'_0$) listed in Example \ref{ex:slice-coslice} such that there exists an isomorphism $v \cong v_0 \oplus C$ (resp.\ $v' \cong v'_0 \oplus C$). By Part c) this means that in a maximal splitting $v \cong v_0 \oplus C$ (resp.\ $v' \cong v'_0 \oplus C$) the indecomposable homomorphism $v_0 \colon B \rightarrow \Z_2$ (resp. $v'_0 \colon \Z_2 \rightarrow B$) is determined up to isomorphism. This also determines the isomorphism class of $C$, because $A \cong B \oplus C$
\end{proof}

We used the following lemmas in the proof of Proposition~\ref{p:slice-coslice-decomp}:

\begin{lem} \label{lem:hom-decomp}
Let $G$ be a finitely generated abelian group and $f \colon G \rightarrow \Z_2$ a homomorphism such that $f \big| _{\Tor(G)} \neq 0$. Let $a = \min \left\{ \text{order of $x$} \mid x \in \Tor(G), f(x)=1 \right\}$, and let $g \in G$ be an element such that $f(g)=1$ and $ag=0$. Then $a$ is a power of $2$ and there is a subgroup $H \leq \Ker(f)$ such that $G = (g) \oplus H$.
\end{lem}

\begin{proof}
If $a$ is not a power of $2$, then $a = 2^bc$ for some $c>1$ odd, so $f(cg)=1$ and $cg$ has order $2^b < a$ contradicting the minimality of $a$. Therefore $a = 2^k$ for some $k$. We get similarly that if $f(x)=1$ for some $x \in \Tor(G)$, then $2^k$ divides the order of $x$. 

Consider the short exact sequence
\[
0 \rightarrow (g) \rightarrow G \rightarrow G / (g) \rightarrow 0 .
\]
We will construct the subgroup $H$ as the image of a splitting map $G / (g) \rightarrow G$. We write $G / (g)$ as a direct sum of cyclic groups, then for each generator $x$ we need to find an element $y \in G$ in the inverse image of $x$ such that $f(y)=0$ and if $x$ has finite order $r$, then $ry=0$. Fix a $z \in G$ that gets mapped to $x$. If $x$ has infinite order then we can choose $y=z$ or $y=z{+}g$ depending on whether $f(z)=0$ or $f(z)=1$. If $x$ has finite order $r$, then $rz \in (g)$. If $rz = 0$ and $f(z) = 0$, then we can choose $y=z$. If $f(z)=1$, then $2^k$ divides $r$, so $r(z+g) = 0$, and since $f(z+g)=0$, we can choose $y = z{+}g$. If $rz \neq 0$, then $rz = sg$ for some non-zero $s \in \Z_{2^k}$. Write $r = 2^{b_1}c_1$ and $s = 2^{b_2}c_2$, where $c_1$ and $c_2$ are odd and $b_2 < k$. If $b_1 > b_2$, then $2^{b_2}(2^{b_1-b_2}c_1z-c_2g)=0$ and $f(2^{b_1-b_2}c_1z-c_2g) = 2^{b_1-b_2}f(c_1z)-c_2f(g) = 1$, which contradicts the minimality of $a$. If $b_1 \leq b_2$, then $s \equiv 2^{b_2-b_1}c_2c_3r$ mod $2^k$, where $c_3 \in \Z_{2^k}$ is such that $c_1c_3 \equiv 1$ mod $2^k$. Then $r(z-2^{b_2-b_1}c_2c_3g) = rz - sg = 0$. Moreover, $2^k$ does not divide $r$, so $f(z-2^{b_2-b_1}c_2c_3g) = 0$. Therefore we can set $y = z-2^{b_2-b_1}c_2c_3g$. 
\end{proof}

\begin{lem} \label{lem:free-decomp}
Let $G$ be a finitely generated free abelian group and $f \colon G \rightarrow \Z_2$ a non-trivial homomorphism. Then there is a $g \in G$ and a subgroup $H \leq \Ker(f)$ such that $f(g)=1$ and $G = (g) \oplus H$.
\end{lem}

\begin{proof}
Fix a basis of $G$, then there is a basis element $g$ such that $f(g)=1$ (otherwise $f=0$). For each of the other basis elements $h$ we define $h' \in G$ by $h'=h$ if $f(h)=0$ and $h'=h+g$ otherwise. Then the subgroup generated by these elements $h'$ is a suitable $H$.
\end{proof}

\begin{lem} \label{lem:dir-summand}
Let $G$ be a finitely generated abelian group, and $g \in G$ an element of order $p$ for some prime $p$. Let $a := \max \left\{ b \mid \exists h \in G : g = p^b h \right\}$ and choose an $h$ such that $g = p^a h$. Then the subgroup $(h) \subseteq G$ is a direct summand.
\end{lem}

\begin{proof}
We need to construct a splitting map $G / (h) \rightarrow G$ for the short exact sequence
\[
0 \rightarrow (h) \rightarrow G \rightarrow G / (h) \rightarrow 0 .
\]
Since $G / (h)$ is a direct sum of cyclic groups, it is enough to find for each generator $x$ an element $y \in G$ in the inverse image of $x$ such that $ry=0$, where $r$ is the order of $x$ (if it is infinite, then there is no condition on $y$). Fix a $z \in G$ that gets mapped to $x$, then $rz \in (h)$. If $rz = 0$, then we can set $y = z$. If $rz \neq 0$, then $g = srz$ for some non-zero $s \in \Z_{p^{a+1}}$ (because $g$ is an order-$p$ element in $(h) \cong \Z_{p^{a+1}}$). Write $r = p^{c_1}d_1$ and $s = p^{c_2}d_2$, where $d_1$ and $d_2$ are not divisible by $p$. Then $g = p^{c_1+c_2}d_1d_2z$, so $c_1 + c_2 \leq a$. Also, $p^ah = p^{c_2}d_2rz$ means that $rz = p^{a-c_2}d_3h$ for a $d_3 \in \Z_{p^{a+1}}$ with $d_2d_3 \equiv 1$ mod $p$. Therefore $rz = rv$, where $v = p^{a-c_1-c_2}d_3d_4h \in (h)$ for a $d_4 \in \Z_{p^{a+1}}$ with $d_1d_4 \equiv 1$ mod $p^{c_2+1}$. So we can choose $y = z - v$. 
\end{proof}

The statements of Proposition \ref{p:slice-coslice-decomp} are summarised in the following theorem: 

\begin{thm}[The classification of 
abelian groups over or under $\Z_2$]
\label{t:(co)slice-d}
%%%%%%%%%%%%%%%
\quad

a) Example \ref{ex:slice-coslice} gives a complete list of indecomposable homomorphisms in $\AbF / \Z_2$ and $\Z_2 / \AbF$ (up to isomorphism). 

b) Every homomorphism in $\AbF / \Z_2$ or $\Z_2 / \AbF$ has a maximal splitting, and the components of a maximal splitting are unique up to isomorphism.
\qed
\end{thm}

\subsection{The categories $\FP_+$ and $\FP_-$} \label{ss:cat}

In this section we will identify $\FP_+$ and $\FP_-$ with $\AbF / \Z_2$ and $\Z_2 / \AbF$ (up to equivalence and isomorphism, respectively). 

Recall that in Definition \ref{d:q-Wu} we defined the homomorphism $v_Q \colon SQ \to \Z_2$ (induced by ${\rm h}$) if $Q$ is a symmetric form parameter, and the homomorphism $v'_Q \colon \Z_2 \to Q_e$ (induced by ${\rm p}$) if $Q$ is an anti-symmetric form parameter.

\begin{example} \label{ex:wu}
The following is the list of homomorphisms $v_Q$ and $v'_Q$ for the form parameters listed in Definition \ref{def:indec}:
\[
\begin{aligned}
& \ul{\quad \epsilon = +1 \quad} & &  \ul{\quad \epsilon = -1 \quad} \\
 v_{Q_{+}} &\cong 0 \colon 0 \rightarrow \Z_2
&  v'_{Q_{-}} & = \iota_1 \colon \Z_2 \rightarrow \Z_2\\
v_{\Z^P} &\cong 1_{\infty} \colon \Z \rightarrow \Z_2 \\
v_{\Z^P_k} &\cong 1_{k+1} \colon \Z_{2^{k+1}} \rightarrow \Z_2 \quad\quad\quad\quad 
 & v'_{\Z^{\Lambda}_k} &= \iota_k \colon \Z_2 \rightarrow \Z_{2^k} \\
v_{Q^{+}} &\cong 1_1 \colon \Z_2 \rightarrow \Z_2 
 & v'_{Q^{-}} &= 0 \colon \Z_2 \rightarrow 0 
\end{aligned}
\]
\end{example}

\begin{remark} \label{rem:v-decomp}
Note that if $P \cong Q \oplus G$ is a splitting of a symmetric form parameter $P$, then we have $SP \cong SQ \oplus G$ and $v_P \cong v_Q \oplus G \colon SP \cong SQ \oplus G \rightarrow \Z_2$. Similarly, if $P \cong Q \oplus G$ is a splitting of an anti-symmetric form parameter $P$, then $v'_P \cong v'_Q \oplus G \colon \Z_2 \rightarrow P_e \cong Q_e \oplus G$.
\end{remark}

\begin{remark}
The height of a homomorphism to or from $\Z_2$ can be defined in such a way that the height of a form parameter from Definition~\ref{d:hght} is
equal to the height of its quasi-Wu class; we leave the details to the reader.
\end{remark}

\begin{prop}
The assignment $P \mapsto v_P$ induces a functor $v_{(-)} \colon \FP_+ \rightarrow \AbF / \Z_2$. Similarly, the assignment $P \mapsto v'_P$ induces a functor $v'_{(-)} \colon \FP_- \rightarrow \Z_2 / \AbF$. 
\end{prop}

\begin{proof}
A morphism $\alpha \colon P \rightarrow P'$ of 
symmetric 
form parameters induces a homomorphism $S\alpha \colon SP \rightarrow SP'$ such that $v_{P'} \circ S\alpha = v_P$, so we can set $v_{\alpha} = S\alpha$. Similarly, a morphism $\alpha \colon P \rightarrow P'$ of anti-symmetric 
form parameters is a homomorphism $\alpha \colon P_e \rightarrow P'_e$ such that $v'_{P'} = \alpha \circ v'_P$, so we can set $v'_{\alpha} = \alpha$.
\end{proof}

\begin{thm} \label{thm:asymm-isom}
The functor $v'_{(-)} \colon \FP_- \rightarrow \Z_2 / \AbF$ is an isomorphism of categories.
\end{thm}

\begin{proof}
To a finitely generated abelian group $G$ with homomorphism $f \colon \Z_2 \rightarrow G$ (an object in $\Z_2 / \AbF$) we assign the anti-symmetric form parameter $(G, 0, {\rm p}_f)$, where ${\rm p}_f(1) = f(1)$. A homomorphism $g \colon G \rightarrow G'$ which is a morphism between $f \colon \Z_2 \rightarrow G$ and $f' \colon \Z_2 \rightarrow G'$ (i.e.\ $g \circ f = f'$) is also a morphism between $(G, 0, {\rm p}_f)$ and $(G', 0, {\rm p}_{f'})$. The functor defined this way is inverse to $v'_{(-)}$.
\end{proof}

\begin{thm} \label{thm:symm-equiv}
The functor $v_{(-)} \colon \FP_+ \rightarrow \AbF / \Z_2$ is an equivalence of categories.
\end{thm}

\begin{proof}
By Proposition \ref{p:slice-coslice-decomp} b) every homomorphism $v \colon A \rightarrow \Z_2$ is isomorphic to some $v_0 \oplus C$ where $v_0$ is one of the homomorphisms listed in Example \ref{ex:slice-coslice} and $C$ is an abelian group. It follows from Example \ref{ex:wu} that there is a symmetric form parameter $Q$ such that $v_0 \cong v_Q$. Then by Remark \ref{rem:v-decomp} we have $v \cong v_0 \oplus C \cong v_Q \oplus C \cong v_{Q \oplus C}$. This shows that $v_{(-)}$ is essentially surjective. By Lemma \ref{l:wu-mor-bij} below $v_{(-)}$ is also full and faithful, hence it is an equivalence of categories. 
\end{proof}

\begin{lem} \label{l:wu-mor-bij}
%%%%%%%%%%%%%%%%%%
For any pair of symmetric form parameters, $P$ and $P'$, the map $\Hom(P,P') \rightarrow \Hom(v_P, v_{P'})$ given by the functor $v_{(-)}$ is a bijection. 
\end{lem}

\begin{proof}
We will describe an inverse. Let $f \colon SP \rightarrow SP'$ be a morphism between $v_P$ and $v_{P'}$. Let $\pi_P \colon P_e \rightarrow SP$ and $\pi_{P'} \colon P'_e \rightarrow SP'$ denote the quotient maps (see the diagram below). 

\[
\xymatrix{
0 \ar[rr] & & \Z \ar@{=}[dd] \ar[rr]^-{{\rm p}_P} \ar@{=}[dr] & & P_e \ar@{-->}[dd]_(.7){\alpha_f} \ar[rr]^-{\pi_P} \ar[dr]^(.55){{\rm h}_P} & & SP \ar[dd]_(.7){f} \ar[dr]^(.55){v_P} \ar[rr] & & 0 & \\
 & 0 \ar '[r] [rr] & \, & \Z \ar '[r]^-{2} [rr] & & \Z \ar '[r]^-{\varrho_2} [rr] & & \Z_2 \ar[rr] & & 0 \\
0 \ar[rr] & & \Z \ar[rr]_-{{\rm p}_{P'}} \ar@{=}[ur] & & P'_e \ar[rr]_-{\pi_{P'}} \ar[ur]_(.55){{\rm h}_{P'}} & & SP' \ar[ur]_(.55){v{\vphantom{h}}_{P'}} \ar[rr] & & 0 & 
}
\] 
\quad \\
Since $\varrho_2 \circ {\rm h}_P = v_P \circ \pi_P = v_{P'} \circ f \circ \pi_P$, it follows from Lemma \ref{l:Wu-pb} (applied to $P'$) that there is a unique homomorphism $\alpha_f \colon P_e \rightarrow P'_e$ such that ${\rm h}_{P'} \circ \alpha_f = {\rm h}_P$ and $\pi_{P'} \circ \alpha_f = f \circ \pi_P$. We have ${\rm h}_{P'} \circ \alpha_f \circ {\rm p}_P = {\rm h}_P \circ {\rm p}_P = 2\Id_{\Z} = {\rm h}_{P'} \circ {\rm p}_{P'}$ and $\pi_{P'} \circ \alpha_f \circ {\rm p}_P = f \circ \pi_P \circ {\rm p}_P = 0 = \pi_{P'} \circ {\rm p}_{P'}$, hence $\alpha_f \circ {\rm p}_P = {\rm p}_{P'}$. Therefore $\alpha_f$ is a morphism between $P$ and $P'$. 

Since $S\alpha_f = f$ for any morphism $f$ between $v_P$ and $v_{P'}$, and $\alpha_{S\beta} = \beta$ for any morphism $\beta \colon P \rightarrow P'$, the assignment $f \mapsto \alpha_f$ is the inverse of the map $\Hom(P,P') \rightarrow \Hom(v_P, v_{P'})$ given by the functor $v_{(-)}$, therefore both of them are bijections. 
\end{proof}

\begin{prop} \label{prop:indec-equiv}
A symmetric (resp.\ anti-symmetric) form parameter $P$ is indecomposable if and only if $v_P$ (resp.\ $v'_P$) is indecomposable.
\end{prop}

\begin{proof}
It follows from Remark \ref{rem:v-decomp} that if $P$ is decomposable, then $v_P$ (resp.\ $v'_P$) is decomposable too. On the other hand, if $v_P$ (resp.\ $v'_P$) has a splitting $v_P \cong v_0 \oplus C$ (resp.\ $v'_P \cong v'_0 \oplus C$), then it follows from Theorem \ref{thm:symm-equiv} (resp.\ \ref{thm:asymm-isom}) and Remark \ref{rem:v-decomp} that $v_0 \cong v_Q$ (resp.\ $v'_0 \cong v'_Q$) for some $Q$, so $v_P \cong v_Q \oplus C \cong v_{Q \oplus C}$ (resp.\ $v'_P \cong v'_Q \oplus C \cong v'_{Q \oplus C}$) and hence $P \cong Q \oplus C$.
\end{proof}

\begin{proof}[Proof of Theorem \ref{t:qfp-class}]
By Theorems \ref{thm:symm-equiv} and \ref{thm:asymm-isom} the functors $v_{(-)}$ and $v'_{(-)}$ are both equivalences of categories. By Example \ref{ex:wu} they send the form parameters of Definition \ref{def:indec} to the homomorphisms in Example \ref{ex:slice-coslice}, and they preserve splittings and indecomposability (see Proposition \ref{prop:indec-equiv}).  Therefore all statements follow immediately from the corresponding parts of Theorem \ref{t:(co)slice-d}.
\end{proof}

\subsection{The category of indecomposable form parameters} \label{ss:indec-mor}
%%%%%%%%%%%%%%%%%%%%%%%%%%%%%%%%%%%%%%%%%%%%%%%%%%%%%%%%
Now we consider morphisms between form parameters. We will give an explicit description of the morphisms between indecomposable form parameters, equivalently, of the full subcategory of indecomposable form parameters. In the symmetric case we will use Theorem \ref{thm:symm-equiv} to first consider the morphisms in $\AbF / \Z_2$, and then take the corresponding morphisms between form parameters.

\begin{defin}[Standard morphisms] \label{def:std-mor}
Consider the diagram \\
\[
\xymatrix{
\Z \, \ar@{^{(}->}[r]^-{\left(\!\begin{smallmatrix} 2 \\ -1 \end{smallmatrix}\!\right)} & \Z \oplus \Z \ar@{^{(}->}@(ur,ul)[]_-{\left(\!\begin{smallmatrix} 1 & 0 \\ n & (2n+1) \end{smallmatrix}\!\right)} \ar@/^11.5pt/@{->>}[rr] \ar@/^20pt/@{->>}[rrrr]^(.8){\left(\!\begin{smallmatrix} 1 & 0 \\ 0 & 1 \end{smallmatrix}\!\right)} \ar@/^30pt/@{->>}[rrrrrr]^(.8){\left(\!\begin{smallmatrix} 1 & 0 \\ 0 & 1 \end{smallmatrix}\!\right)} \ar@/^40pt/@{->>}[rrrrrrrr]^(.8){\left(\!\begin{smallmatrix} 1 & 0 \\ 0 & 1 \end{smallmatrix}\!\right)} & & \,\, \vphantom{\Z^2_2} \ldots \,\, \ar@{->>}[rr] & & \Z \oplus \Z_{2^3} \ar@{->>}[rr]^(.4){\left(\!\begin{smallmatrix} 1 & 0 \\ 0 & 1 \end{smallmatrix}\!\right)} & & \Z \oplus \Z_{2^2} \ar@{->>}[rr]^(.4){\left(\!\begin{smallmatrix} 1 & 0 \\ 0 & 1 \end{smallmatrix}\!\right)} & & \Z \oplus \Z_2 \ar@{->>}[r]^(.55){(1,0)} & \Z
}
\] 
\quad \\
(where we have a homomorphism $\left(\begin{smallmatrix} 1 & 0 \\ n & (2n+1) \end{smallmatrix}\right) \colon \Z \oplus \Z \rightarrow \Z \oplus \Z$ for every $n \in \Z$). It determines a diagram \\ \\ \\
\[
\xymatrix{
Q_{+} \ar[r] & \Z^P \ar@(ur,ul)[] \ar@/^8pt/[rr] \ar@/^12pt/[rrrr] \ar@/^18pt/[rrrrrr] \ar@/^24pt/[rrrrrrrr] & & {} \vphantom{\Z} \ldots \ar[rr] & & \Z^P_3 \ar[rr] & & \Z^P_2 \ar[rr] & & \Z^P_1 \ar[r] & Q^{+}
}
\]
\quad \\
of indecomposable symmetric form parameters and morphisms. Similarly, the diagram \\
\[
\xymatrix{
\Z_2 \, \ar@{^{(}->}[r]^{2} & \Z_{2^2} \, \ar@{^{(}->}[rr]^{2} \ar@/_24pt/@{-}[rrrrrrrr] & & \Z_{2^3} \, \ar@{^{(}->}[rr]^{2} \ar@/_18pt/@{-}[rrrrrr] & & \Z_{2^4} \, \ar@{^{(}->}[rr] \ar@/_12pt/@{->>}[rrrr] & & \ldots \vphantom{\Z} \ar@/_7pt/@{-}[rr] & & 0 
}
\]
\quad \\ \\ \\
determines a diagram
\[
\xymatrix{
Q_{-} \ar[r] & \Z^{\Lambda}_2 \ar[rr] \ar@/_24pt/@{-}[rrrrrrrr] & & \Z^{\Lambda}_3 \ar[rr] \ar@/_18pt/@{-}[rrrrrr] & & \Z^{\Lambda}_4 \ar[rr] \ar@/_12pt/[rrrr] & & \ldots \vphantom{\Z} \ar@/_7pt/@{-}[rr] & & Q^{-} 
}
\]
\quad \\ \\ \\
of indecomposable anti-symmetric form parameters and morphisms. The morphisms between the form parameters in these two diagrams will be called {\em standard morphisms}.
\end{defin}

\begin{thm} \label{thm:mor}
Every morphism $\alpha \colon P \rightarrow P'$ between the indecomposable form parameters $P = (P_e, {\rm h}, {\rm p})$ and $P' = (P'_e, {\rm h}', {\rm p}')$ is the composition of standard morphisms and isomorphisms.
\end{thm}

\begin{proof}
First we consider the symmetric case. By Theorem \ref{thm:symm-equiv} it is enough to prove the analogous statement in $\AbF / \Z_2$. Using the bijection from Lemma \ref{l:wu-mor-bij}, we need to consider the diagram \\
\[
\xymatrix{
0 \, \ar@{^{(}->}[r] & \Z \ar@{^{(}->}@(ur,ul)[]_-{2n+1} \ar@/^11.5pt/@{->>}[rr] \ar@/^20pt/@{->>}[rrrr]^(.8){1} \ar@/^30pt/@{->>}[rrrrrr]^(.8){1} \ar@/^40pt/@{->>}[rrrrrrrr]^(.8){1} & & \,\, \vphantom{\Z^2_2} \ldots \,\, \ar@{->>}[rr] & & \Z_{2^4} \ar@{->>}[rr]^-{1} & & \Z_{2^3} \ar@{->>}[rr]^-{1} & & \Z_{2^2} \ar@{->>}[r]^-{1} & \Z_2
}
\] 
\quad \\
(where we have a homomorphism $2n{+}1 \colon \Z \rightarrow \Z$ for every $n \in \Z$), which determines a diagram \\ \\ \\
\[
\xymatrix{
0 \ar[r] & 1_{\infty} \ar@(ur,ul)[] \ar@/^8pt/[rr] \ar@/^12pt/[rrrr] \ar@/^18pt/[rrrrrr] \ar@/^24pt/[rrrrrrrr] & & {} \vphantom{\Z} \ldots \ar[rr] & & 1_4 \ar[rr] & & 1_3 \ar[rr] & & 1_2 \ar[r] & 1_1
}
\]
\quad \\
between the standard indecomposable homomorphisms in $\AbF / \Z_2$. We will show that any morphism $v_1 \rightarrow v_2$ between indecomposable homomorphisms in $\AbF / \Z_2$ is the composition of some morphisms in the above diagram and isomorphisms.

This is obvious if $v_1$ is initial (isomorphic to $0 \colon 0 \rightarrow \Z_2$), so we will assume that $v_1 \not\cong 0$. Then $v_1$ is surjective, so $v_2$ is also surjective. Therefore both $v_1$ and $v_2$ are isomorphic to either $1_{\infty}$ or $1_k$ for some $k \geq 1$.

If $v_1$ and $v_2$ are both isomorphic to $1_{\infty}$, then we can fix isomorphisms $v_1 \cong 1_{\infty}$ and $v_2 \cong 1_{\infty}$. Under these identifications the morphism $v_1 \rightarrow v_2$ corresponds to a homomorphism $\Z \rightarrow \Z$ that is compatible with the surjective map $1_{\infty} \colon \Z \rightarrow \Z_2$, equivalently, it is given by multiplication with some odd integer $2n+1$. 

If $v_1 \cong 1_k$ for some $k \geq 1$ and $v_2 \cong 1_l$ for some $l > k$ (or $v_2 \cong 1_{\infty}$), then the morphism $v_1 \rightarrow v_2$ corresponds to a homomorphism $\Z_{2^k} \rightarrow \Z_{2^l}$ (or $\Z_{2^k} \rightarrow \Z$). The image of any such homomorphism is contained in $\Ker(v_2)$, but $v_1 \neq 0$, showing that no morphism $v_1 \rightarrow v_2$ can exist in this case.

If $v_1 \cong 1_k$ for some $k \geq 1$ (or $v_1 \cong 1_{\infty}$) and $v_2 \cong 1_l$ for some $l \leq k$, then the morphism $v_1 \rightarrow v_2$ corresponds to a homomorphism $\Z_{2^k} \rightarrow \Z_{2^l}$ (or $\Z \rightarrow \Z_{2^l}$) compatible with the maps $1_k$ (or $1_{\infty}$) and $1_l$. Equivalently, this homomorphism sends $1$ to some odd $p \in \Z_{2^l}$, so it is the composition of the standard map $1 \colon \Z_{2^k} \rightarrow \Z_{2^l}$ (or $1 \colon \Z \rightarrow \Z_{2^l}$) and the automorphism of $\Z_{2^l}$ given by multiplication with $p$.

Now consider the anti-symmetric case. Again, the statement is obvious if $P$ is initial (isomorphic to $Q_{-}$) or $P'$ is terminal (isomorphic to $Q^{-}$). So we assume that $P \not \cong Q_{-}$ (so ${\rm p}(1)$ is divisible by $2$, therefore ${\rm p}'(1)$ is divisible by $2$, hence $P' \not \cong Q_{-}$) and $P' \not \cong Q^{-}$ (which means that ${\rm p}' \neq 0$, so ${\rm p} \neq 0$, therefore $P \not \cong Q^{-}$).

If $P \cong \Z^{\Lambda}_k$ and $P' \cong \Z^{\Lambda}_l$ with $k > l$, then ${\rm p}(1)$ is divisible by $2^{k-1} \geq 2^l$, so the same is true for ${\rm p}'(1)$, so ${\rm p}'(1)=0$. This contradiction shows that $\Z^{\Lambda}_k \rightarrow \Z^{\Lambda}_l$ morphisms with $k > l$ do not exist.

If $P \cong \Z^{\Lambda}_k$ and $P' \cong \Z^{\Lambda}_l$ with $k \leq l$, then $\alpha$ corresponds to a homomorphism $\Z_{2^k} \rightarrow \Z_{2^l}$ that sends ${\rm p}(1)=2^{k-1}$ to ${\rm p}'(1) = 2^{l-1}$. This means that $1 \in \Z_{2^k}$ has to be sent into $2^{l-k}p$ for some odd $p \in \Z_{2^l}$. Since multiplication by $p$ is an isomorphism of $\Z^{\Lambda}_l$, we can describe $\alpha$ as the composition of the inclusion $2^{l-k} \colon \Z_{2^k} \rightarrow \Z_{2^l}$ (which is the composition of $(l-k)$ standard morphisms) and $3$ isomorphisms. 
\end{proof}

\begin{defin}
${\rm Aut}(Q) \subseteq {\rm Aut}(Q_e)$ is the group of automorphisms of the form parameter $Q$.
\end{defin}

\begin{thm} \label{thm:aut}
The symmetric form parameters $Q_{+}$ and $Q^{+}$ have no non-trivial automorphisms.
In contrast, ${\rm Aut}(\Z^P) = \Z_{2}(\beta)$ (i.e.\ it is a group isomorphic to $\Z_{2}$, generated by $\beta$) and 
\[
{\rm Aut}(\Z^P_k) \cong 
\begin{cases}
\Z_2(\beta_k) & \text{if $k = 1$} \\
\Z_2(\beta_k) \oplus \Z_{2^{k-1}}(\gamma_k) & \text{if $k \geq 2$,}
\end{cases}
\]
where the generators are
\[
\beta = 
\begin{pmatrix}
1 & 0 \\ 
-1 & -1
\end{pmatrix}
\in {\rm Aut}(\Z \oplus \Z)
\text{\,,} \quad
\beta_k = 
\begin{pmatrix}
1 & 0 \\ 
-1 & -1
\end{pmatrix}
\in {\rm Aut}(\Z \oplus \Z_{2^k})
\quad \text{and} \quad 
\gamma_k = 
\begin{pmatrix}
1 & 0 \\ 
1 & 3
\end{pmatrix} 
\in {\rm Aut}(\Z \oplus \Z_{2^k})
. 
\]
In the anti-symmetric case $Q_{-}$ and $Q^{-}$ have no non-trivial automorphisms and for $k \geq 2$
\[
{\rm Aut}(\Z^{\Lambda}_k) \cong {\rm Aut}(\Z_{2^k}) \cong 
\begin{cases}
\Z_2(-1) & \text{if $k = 2$,} \\
\Z_2(-1) \oplus \Z_{2^{k-2}}(3) & \text{if $k \geq 3$.}
\end{cases}
\]
\end{thm}

\begin{proof}
Initial and terminal objects do not have non-trivial automorphisms. 

In the remaining anti-symmetric case, an automorphism of $\Z^{\Lambda}_k$ is an automorphism of the group $\Z_{2^k}$ that preserves ${\rm p}(1)$, the order-$2$ element. Since this is true for all automorphisms of $\Z_{2^k}$, we get that ${\rm Aut}(\Z^{\Lambda}_k) \cong {\rm Aut}(\Z_{2^k})$. This is isomorphic to $\Z_{2^k}^{\times}$, the multiplicative group of odd elements of $\Z_{2^k}$, which is $\Z_2$ if $k=2$ and $\Z_2 \oplus \Z_{2^{k-2}}$ (generated by $-1$ and $3$) if $k \geq 3$. 

Now we consider the symmetric case. By Theorem \ref{thm:symm-equiv} we have ${\rm Aut}(P) \cong {\rm Aut}(v_P)$. We start with ${\rm Aut}(\Z^P) \cong {\rm Aut}(v_{\Z^P}) \cong {\rm Aut}(1_{\infty})$. From the proof of Theorem \ref{thm:mor} we see that the monoid of morphisms $1_{\infty} \rightarrow 1_{\infty}$ is isomorphic to the multiplicative monoid of odd integers. The group of automorphisms of $1_{\infty}$ is the group of invertible elements in this monoid, so it is isomorphic to $\Z_2$, generated by $-1$. The generator corresponds to the matrix $\beta$ under the bijection of Lemma \ref{l:wu-mor-bij}.

Next we consider ${\rm Aut}(\Z^P_k) \cong {\rm Aut}(1_{k+1})$. As we saw in the proof of Theorem \ref{thm:mor}, this group is isomorphic to ${\rm Aut}(\Z_{2^{k+1}}) \cong \Z_{2^{k+1}}^{\times}$. The generator $-1 \in \Z_{2^{k+1}}^{\times}$ corresponds to $\beta_k$, and (if $k \geq 2$) $3$ corresponds to $\gamma_k$.
\end{proof}

%%%%%%%%%%%%%%%%%%%%%%%%%%%%%%%%%%%%%%%%%%%%%%%%%%%%%%%%%%

\section{Witt groups of split form parameters}
\label{s:WG}
%%%%%%%%%%%%%%%%%%%%%%%%%%%%%%%%%%%%%%%%%%%%%%%%%%%%%%%%%%

In this section we study the Witt groups $W_0(P)$ of form parameters $P$. First we compute the Witt groups of indecomposable form parameters and the homomorphisms induced by morphisms of such form parameters, giving a complete description of the Witt group functor restricted to the subcategory of indecomposable form parameters. 
Then, given a split form parameter $P = Q \oplus G$, where $Q$ is a form parameter and $G$ is an abelian group, we use the fact that the Witt group $W_0(Q \oplus G)$ decomposes as the direct sum of $W_0(Q)$ and the reduced Witt group $W_0^Q(G)$. We prove that $W_0^Q(G)$ is naturally isomorphic to Baues' quadratic tensor product $G \otimes_{\Z} Q$ \cite[\S 4]{Bau1}, and by our earlier classification of form parameters this allows us to determine the Witt group $W_0(P)$ for any form parameter $P$.  
Finally, we give some examples of computations of Witt groups and 
induced homomorphisms between Witt groups.

\subsection{Witt groups of indecomposable form parameters} \label{ss:Witt-ind}
%%%%%%%%%%%%%%%%%%%%%%%%%%%%%%%%%%%%%%%%%%%%%%%%%%%%%%%%
Let $Q$ be an indecomposable form parameter. 
Recall that by Theorem \ref{t:qfp-class}, Definition \ref{def:indec} gives a complete list of indecomposable form parameters (up to isomorphism). First we will define the invariants that classify the elements of the groups $W_0(Q)$.

For a symmetric form parameter $Q$ the signature induces a homomorphism
\[
\sigma \colon W_0(Q) \to \Z.
\]
Recall that for any symmetric form parameter $Q$
there is an extended symmetrisation map 
$\es \colon Q \rightarrow Q^+ \oplus SQ$ (see Definition \ref{d:es}). It induces a homomorphism 
$\HomZ(X,\es) \colon \HomZ(X,Q) \rightarrow \HomZ(X,Q^+ \oplus SQ)$ for any abelian group $X$, which sends $(X, \lambda, \mu)$ to $(X, \lambda, (\mu_0,\omega_0))$, where $\mu_0(x) = \lambda(x,x)$ and $v_Q \circ \omega_0 \colon X \to \Z_2$ is characteristic for $\lambda$, i.e.\ $\varrho_2 \circ \lambda(x,x) = v_Q \circ \omega_0(x)$ (see Lemma \ref{l:gmtrc}). We will use this map $\omega_0$ (more precisely a lift $\omega \colon X \rightarrow \Z$ of $\omega_0$) to define an invariant of $Q$-forms for $Q = \Z^P$ and $Q = \Z^P_k$.

If $Q =\Z^P$, then $v_{\Z^P} = 1_{\infty} = \varrho_2 \colon S\Z^P \cong \Z \rightarrow \Z_2$, and we define $\omega = \omega_0$, it satisfies
\[ 
\lambda(x,x) \equiv \omega(x) \text{${\rm~mod}~2$.} 
\]
If $\lambda$ is nonsingular with adjoint $\wh \lambda$, then for the element $\widehat{\omega} := \wh \lambda^{-1}(\omega) \in X$ it is well known that $\widehat{\omega}^2 := \lambda(\widehat{\omega}, \widehat{\omega}) \in \Z$ is congruent to $\sigma(\lambda)$ mod $8$ (see e.g.\ \cite[II/{\S}5]{M-H}). 
Recall that $\HomZns(Q)$ denotes the set of isomorphism classes of nonsingular $Q$-forms.

\begin{defin} \label{def:rho-inf}
We define the map $\rho_{\infty} \colon \HomZns(\Z^P) \rightarrow \Z$ by 
\[
\rho_{\infty}(X, \lambda, \mu) = \frac{\sigma(\lambda) - \widehat{\omega}^2}{8} .
\]
It is easy to see that this is a monoid homomorphism, and that $\rho_{\infty}$ vanishes on metabolic forms, so it induces a homomorphism $\rho_{\infty} \colon W_0(\Z^P) \rightarrow \Z$.  
\end{defin}

If $Q = \Z^P_k$, then $S\Z^P_k \cong \Z_{2^{k+1}}$ and $v_{\Z^P_k} = 1_{k+1} \colon \Z_{2^{k+1}} \rightarrow \Z_2$, so we have $\lambda(x, x) \equiv \omega_0(x)$ mod $2$. Now lift $\omega_0$ to a homomorphism $\omega \colon X \to \Z$, then $\lambda(x, x) \equiv \omega(x)$ mod $2$, and consider $\widehat{\omega}^2$. Again we have $\widehat{\omega}^2 \equiv \sigma(\lambda)$ mod $8$, and it is easy to check that $\widehat{\omega}^2$ mod $2^{k+2}$ is independent of the choice of the lift $\omega$.

\begin{defin} \label{def:rho-k}
We define the map $\rho_k \colon \HomZns(\Z^P_k) \rightarrow \Z_{2^{k-1}}$ by 
\[
\rho_k(X, \lambda, \mu) = \frac{\sigma(\lambda) - \widehat{\omega}^2}{8} \text{ mod $2^{k-1}$} .
\]
It induces a homomorphism $\rho_k \colon W_0(\Z^P_k) \rightarrow \Z_{2^{k-1}}$. 
\end{defin}

In the anti-symmetric case we will denote the Arf invariant by 
\[ 
c \colon W_0(Q_{-}) \rightarrow \Z_2 .
\]
The following result entails and extends Theorem~\ref{t:W_0-ind}.

\begin{thm}  \label{t:wf-ind}
%%%%%%%%%%%%%%%%%%%%%%%%%%%%%%
Applying the functor $W_0$ to the diagrams in Definition \ref{def:std-mor} gives the following diagrams of abelian groups, where we use the duals of the invariants detecting each summand to name generators:

\[
\!\!\!
\xymatrix{
\Z(8\sigma^*) \, \ar@{^{(}->}[r]^-{\left(\!\begin{smallmatrix} 8 \\ 1 \end{smallmatrix}\!\right)} & \Z(\sigma^*) \oplus \Z(\rho_{\infty}^*) \ar@{^{(}->}@(ur,ul)[]_-{\left(\!\begin{smallmatrix} 1 & 0 \\ -n(n+1)/2 & (2n+1)^2 \end{smallmatrix}\!\right)} \ar@/^11.5pt/@{->>}[r] \ar@/^20pt/@{->>}[rr]^(.8){\left(\!\begin{smallmatrix} 1 & 0 \\ 0 & 1 \end{smallmatrix}\!\right)} \ar@/^30pt/@{->>}[rrr]^(.8){\left(\!\begin{smallmatrix} 1 & 0 \\ 0 & 1 \end{smallmatrix}\!\right)} \ar@/^40pt/@{->>}[rrrr]^(.8){(1,0)} & \ldots\vphantom{\Z^2_2} \ar@{->>}[r] & \Z(\sigma^*) \oplus \Z_{2^2}(\rho_3^*) \ar@{->>}[r]^(.45){\left(\!\begin{smallmatrix} 1 & 0 \\ 0 & 1 \end{smallmatrix}\!\right)} & \Z(\sigma^*) \oplus \Z_2(\rho_2^*) \ar@{->>}[r]^(.55){(1,0)} & \Z(\sigma^*) \ar[r]^-{1} & \Z(\sigma^*)
}
\] 

\[
\xymatrix{
\Z_2(c^*) \ar[r] & 0 \ar[r] \ar@/_12pt/@{-}[rrrr] & 0 \ar[r] \ar@/_9pt/[rrr] & 0 \ar[r] \ar@/_6pt/@{-}[rr] & \ldots\vphantom{\Z} \ar@/_4pt/@{-}[r] & 0 
}
\]
\quad \\

In the symmetric case, the generators of the automorphism groups of the form parameters (see Theorem \ref{thm:aut}) induce the following automorphisms of the Witt groups: $W_0(\beta) = {\rm Id}$, $W_0(\beta_k) = {\rm Id}$ and $W_0(\gamma_k) = \left( \begin{smallmatrix} 1 & 0 \\ -1 & 9 \end{smallmatrix} \right)$. 
\end{thm}

\begin{proof}
First we consider the Witt groups in the symmetric case. It is a classical result that taking the signature of a symmetric form defines an isomorphism $\sigma\colon W_0(Q^{+}) \cong \Z$ and that $W_0(Q_{+}) \to W_0(Q^{+})$ is injective with image the subgroup of index $8$; see, for example, \cite[II Corollary 4.4, Lemma 5.2 \& Lemma 6.2]{M-H}.
By \cite[Theorem 2, (2)]{Wa1} we have $W_0(\Z^P) \cong \Z \oplus \Z$, where the isomorphism is given by the signature $\sigma \colon W_0(\Z^P) \rightarrow \Z$ and $\rho_{\infty}$.

The standard morphism $\Z^P \rightarrow \Z^P_k$ given by $(\Id,\varrho_{2^k}) \colon \Z \oplus \Z \rightarrow \Z \oplus \Z_{2^k}$ induces a homomorphism $W_0(\Z^P) \rightarrow W_0(\Z^P_k)$ that commutes with the signature maps. It follows from the definitions of $\rho_{\infty}$ and $\rho_k$ that the diagram
\[
\xymatrix{
W_0(\Z^P) \ar[r] \ar[d]_-{\rho_{\infty}} & W_0(\Z^P_k) \ar[d]^-{\rho_k} \\
\Z \ar[r]^-{\varrho_{2^{k-1}}} & \Z_{2^{k-1}}
}
\]
commutes. 
This implies that $\Ker(W_0(\Z^P) \rightarrow W_0(\Z^P_k)) \leq \Ker((\Id_{\Z} \oplus \varrho_{2^{k-1}}) \circ (\sigma \oplus \rho_{\infty})) = (\sigma \oplus \rho_{\infty})^{-1}(0 \oplus 2^{k-1} \Z)$. 
Now consider the $\Z^P_k$-form $\ul \phi$ and $\Z^P$-form $\ul \phi'$:
\[  \ul \phi  = \left( \Z^2(x,y),
\left[ \begin{array}{cc}
1 & 1 \\
1 & 0
\end{array} \right],
\left( \begin{array}{c} (1,2^{k-1}) \\ (0,0) \end{array} \right)\right) 
\quad \quad
\ul \phi'  = \left( \Z^2(x,y),
\left[ \begin{array}{cc}
1 & 1 \\
1 & 0
\end{array} \right],
\left( \begin{array}{c} (1,2^{k-1}) \\ (0,2^k) \end{array} \right)\right).\]
For $\ul \phi'$ we have $\omega = (1+2^k,2^{k+1})$, $\widehat{\omega} = \bigl( \begin{smallmatrix} 2^{k+1} \\ 1-2^k \end{smallmatrix} \bigr)$ and $\widehat{\omega}^2 = 2^{k+2}$, so $(\sigma \oplus \rho_{\infty}) ([\ul \phi']) = (0, -2^{k-1})$.
The image of $[\ul \phi']$ in $W_0(\Z^P_k)$ is $[\ul \phi] = 0$, because $\ul \phi$ is metabolic, showing that $\Ker(W_0(\Z^P) \rightarrow W_0(\Z^P_k)) = (\sigma \oplus \rho_{\infty})^{-1}(0 \oplus 2^{k-1} \Z)$. Finally, by Lemma \ref{l:gmtrc} a $\Z^P_k$-form corresponds to some $(Q^+ \oplus \Z_{2^{k+1}})$-form $(X, \lambda, (\mu_0,\omega_0))$ such that $\varrho_2 \circ \omega_0 \colon X \to \Z_2$ is characteristic for $\lambda$. For any lift $\omega \colon X \to \Z$ of $\omega_0$, $\varrho_2 \circ \omega$ is also characteristic for $\lambda$, so $(X, \lambda, (\mu_0,\omega))$ corresponds to a $\Z^P$-form. Its image is the $\Z^P_k$-form we started with, showing that the map $W_0(\Z^P) \rightarrow W_0(\Z^P_k)$ is surjective. Therefore $W_0(\Z^P_k) \cong W_0(\Z^P) / \Ker(W_0(\Z^P) \rightarrow W_0(\Z^P_k)) \cong (\Z \oplus \Z) / (0 \oplus 2^{k-1} \Z) = \Z \oplus \Z_{2^{k-1}}$, and the isomorphism is given by $(\sigma \oplus \rho_k)$.

In the anti-symmetric case it is a classical result that the Arf invariant defines an isomorphism $W_0(Q_{-}) \cong \Z_2$ (see \cite{K-M}), and that $W_0(Q^{-}) \cong 0$ (as every $Q^{-}$-form is hyperbolic). It remains to show that $W_0(\Z^{\Lambda}_k) \cong 0$ if $k>1$. This is a matter of a small number of simple calculations.  Firstly, recall that every nonsingular anti-symmetric form $(X, \lambda)$ splits as a direct sum of hyperbolic forms and thus every nonsingular $\Z^{\Lambda}_k$-form splits as a sum of rank $2$ forms.  Thus it is enough to show that every rank $2$ form is stably metabolic.  So consider a $\Z^{\Lambda}_k$-form
\[ ( X, \lambda, \mu)  = \left( \Z^2(x,y),
\left[ \begin{array}{cc}
0 & 1 \\
-1 & 0
\end{array} \right], 
\left( \begin{array}{c} p \\ q \end{array} \right)\right),\]
where $p, q \in \Z_{2^k}$.  Suppose to begin that one of $p$ or $q$ is odd; given the symmetry of the form, say that $p$ is odd.  Since $p$ is odd and $k > 1$, $p + 2^{k-1}$ is again odd and so there is an integer $a$ such that $a(p + 2^{k-1}) \equiv -q~{\rm mod}~2^k$.  It follows that the summand generated by $ax + y$ is a lagrangian since
\[\mu(ax + y) = ap + q + a2^{k-1} = 0 \in \Z_{2^{k}}.\]
For the general case one checks that there is an isometry
\[
\left( \Z^2(x, y), 
\left[ \begin{array}{cc}
0 & 1 \\
-1 & 0
\end{array} \right], 
\left( \begin{array}{c} p \\ q \end{array} \right)\right) 
\oplus 
\left( \Z^2(w, z), 
\left[ \begin{array}{cc}
0 & 1 \\
-1 & 0
\end{array} \right], 
\left( \begin{array}{c} 1 \\ 0 \end{array} \right)\right)
\cong \]
\[
\left(\Z^2(x - pw, y), 
\left[ \begin{array}{cc}
0 & 1 \\
-1 & 0
\end{array} \right], 
\left( \begin{array}{c}  0 \\ q \end{array} \right)\right)
\oplus
 \left(\Z^2(w, z + py), 
\left[ \begin{array}{cc}
0 & 1 \\
-1 & 0
\end{array} \right], 
\left( \begin{array}{c} 1 \\ pq \end{array} \right)\right).
\]
By the previous case, the second summand in the bottom line is metabolic and so the form 
$(X, \lambda, \mu)$ above represents zero in $W_0(\Z^{\Lambda}_k)$.

Next we will determine the image of the morphisms under the functor $W_0$. First consider the morphism $\alpha \colon Q_+ \rightarrow \Z^P$ given by $\left(\begin{smallmatrix} 2 \\ -1 \end{smallmatrix}\right) \colon \Z \rightarrow \Z \oplus \Z$. Let $(X, \lambda, \mu)$ be a representative of the generator of $W_0(Q_+) \cong \Z$ with $\sigma(\lambda) = 8$. Its image under $W_0(\alpha)$ is represented by $(X, \lambda, \mu')$, where $\mu' = \alpha \circ \mu$. For this $\mu'$ we have $\omega' = \pi_{\Z^P} \circ \mu' = \pi_{\Z^P} \circ \alpha \circ \mu = S\alpha \circ \pi_{Q_+} \circ \mu = 0$, because $SQ_+ = 0$. Hence $\widehat{\omega}' = 0$ and $\rho_{\infty}([X, \lambda, \mu']) = \frac{\sigma(\lambda) - (\widehat{\omega}')^2}{8} = 1$. Therefore $W_0(\alpha) = \left(\begin{smallmatrix} 8 \\ 1 \end{smallmatrix}\right)$. 

Now let $\alpha \colon \Z^P_k \rightarrow \Z^P_l$ be a morphism, where we use the notation $\Z^P = \Z^P_{\infty}$. In Theorems \ref{thm:mor} and \ref{thm:aut} we saw that such a morphism exists only if $k \geq l$, and it is given by a matrix $\left(\begin{smallmatrix} 1 & 0 \\ n & (2n+1) \end{smallmatrix}\right)$ for some $n \in \Z_{2^l}$ or $n \in \Z$ (if $k = l = \infty$). Let $(X, \lambda, \mu)$ be a form over $\Z^P_k$, and define $\omega_0$, $\omega$ and $\widehat{\omega}$ as in Definition \ref{def:rho-k} (or \ref{def:rho-inf}). Its image under $W_0(\alpha)$ is represented by $(X, \lambda, \mu')$, where $\mu' = \alpha \circ \mu$. Then $\omega_0' = \pi_{\Z^P_l} \circ \mu' = S\alpha \circ \pi_{\Z^P_k} \circ \mu = S\alpha \circ \omega_0$, and the map $S\alpha \colon S\Z^P_k \rightarrow S\Z^P_l$ is given by $S\alpha(1) = 2n+1$. So we can choose the lift $\omega' = (2n+1)\omega$, and then $\widehat{\omega}' = (2n+1)\widehat{\omega}$. So 
\[
\rho_l(X, \lambda, \mu') = \frac{\sigma(\lambda) - (\widehat{\omega}')^2}{8} = \frac{\sigma(\lambda) - (2n+1)^2\widehat{\omega}^2}{8} = (2n+1)^2\frac{\sigma(\lambda) - \widehat{\omega}^2}{8} - \frac{n(n+1)}{2}\sigma(\lambda)
\]
Therefore $W_0(\alpha) = \left(\begin{smallmatrix} 1 & 0 \\ -\frac{n(n+1)}{2} & (2n+1)^2 \end{smallmatrix}\right)$. 
\end{proof}

\begin{remark} \label{r:sm-nm}
Note that for the morphisms $i_k\colon Q_{-} \hra \Z^{\Lambda}_k$ ($k \geq 2$) given by the injective homomorphisms $\iota_k \colon \Z_2 \rightarrow \Z_{2^k}$, the induced maps $W_0(i_k)$ are not injective. This corresponds to the following fact: if $(X, \lambda, \mu)$ is a $(Q_{-})$-form with Arf invariant $1$ then $(X, \lambda, \iota_k \circ \mu)$ is a stably metabolic $\Z^{\Lambda}_k$-form that is not metabolic.  

The same phenomenon occurs in the symmetric case,
and we are grateful to a referee for pointing out the following example.  The $\Z^P_1$-form given by
\[  \left( \Z^2,
\left[ \begin{array}{cc}
0 & 1 \\
1 & 0
\end{array} \right],
\left( \begin{array}{c} (0,1) \\ (0,1) \end{array} \right)\right) \]
has signature zero,
and so it is stably metabolic by Theorem~\ref{t:wf-ind}.
But, by inspection, it is not metabolic.
\end{remark}

\subsection{The quadratic tensor product} \label{ss:qtp}
%%%%%%%%%%%%%%%%%%%%%%%%%%%%%%%%%%%%%%%%%%%%%%%%%%%%
We briefly recall the definition and basic properties of the quadratic tensor product (see Baues \cite[\S 4]{Bau1} and \cite[\S 8]{Bau2}), which we will use in Section \ref{ss:RWG} to compute reduced Witt groups.

\subsubsection{Definitions}

We will introduce the quadratic tensor product in the more general context of quadratic modules.

\begin{defin}
A quadratic module (over $\Z$) is a pair of abelian groups $Q_{ee}$ and $Q_e$ with homomorphisms ${\rm h} \colon Q_e \rightarrow Q_{ee}$ and ${\rm p} \colon Q_{ee} \rightarrow Q_e$ such that ${\rm h p h} = 2 {\rm h}$ and ${\rm p h p} =  2 {\rm p}$.
\end{defin}

In particular, a form parameter is a quadratic module with $Q_{ee} = \Z$. The definition of an extended quadratic form over $Q$, the group $\HomZ(X,Q)$, and a morphism between form parameters can be generalised to quadratic modules $Q$ in the obvious way.

\begin{definition} \label{def:qtp}
%%%%%%%%%%%%%%%%
Let $G$ be an abelian group and $Q$ a quadratic module over $\Z$.
The {\em quadratic tensor product} of $G$ and $Q$ is the abelian group
$G \tensor_{\Z} Q$ which is generated by symbols of the form
\[ 
g \tensor q \quad \text{and} \quad [g, h] \tensor a, 
\]
where $g, h \in G$, $q \in Q_e$ and $a \in Q_{ee}$.
These symbols satisfy the following relations:
\begin{compactitem}
\item{$(g+h)\tensor q  = g \tensor q + h \tensor q + [g,h] \tensor {\rm h}(q)$,}
\item{$[g, g] \tensor a = g \tensor {\rm p}(a)$}
\end{compactitem}
and $g \tensor q$ is linear in $q$ and $[g, h] \tensor a$ is linear in each variable $g$, $h$ and $a$.  
\end{definition}

It immediately follows from the defining relations that $[g, h] \tensor a + [h, g] \tensor a = [g,h] \tensor {\rm h}({\rm p}(a))$, so if $Q$ is a form parameter (i.e.\ $Q_{ee} = \Z$), then $[h, g] \tensor a = \epsilon_Q[g,h] \tensor a$. 

The quadratic tensor product is a functor in both $G$ and $Q$, i.e.\ a homomorphism
$f\colon G \to H$ and a morphism $\alpha \colon Q \to P$ induce homomorphisms
$f_* \colon G \tensor_\Z Q \to H \tensor_\Z Q$ and $\alpha_*\colon G \tensor_\Z Q \to G \tensor_\Z P$.

\begin{example} \label{ex:qtp-funct}
By Baues \cite[(8.16)]{Bau2} for certain indecomposable form parameters $Q$ the quadratic tensor product $G \tensor _{\Z} Q$ can be expressed (up to natural isomorphism) as a well-known quadratic functor of $G$:
\begin{center}
\begin{tabular}{c|l}
$Q$ & $G \tensor _{\Z} Q$ \\ 
\hline
$Q_{+}$  & $S^2(G) = G \tensor G /\{g \tensor h - h \tensor g\} $ \\  
$\Z^P$ & $P^2(G) = \Delta(G)/\Delta(G)^3$ \\
$Q^{+}$  & $\Gamma(G)$ \\ 
$Q_{-}$ &  $\Lambda_1(G) = G \tensor G/\{g \tensor h + h \tensor g \}$ \\
$Q^{-}$  & $\Lambda(G) = G \tensor G/\{g \tensor g \}$ \\
\end{tabular}
\end{center}
Here $\Gamma(G)$ denotes Whitehead's universal quadratic functor and $\Delta(G)$ is the augmentation ideal in the group ring $\Z[G]$ and $\Delta(G)^3$ its third power. 
\end{example}

\subsubsection{Computing the quadratic tensor product and $\HomZ$}

Next we construct short exact sequences and pullback (or pushout) diagrams which can be used to determine the quadratic tensor product $G \otimes_{\Z} Q$ and $\HomZ(G,Q)$ for a symmetric (or anti-symmetric) form parameter $Q$.

\begin{remark} \label{rem:p+p}
We are going to make use of the following fact, which is easily verified by diagram chasing. Let 
\[
\xymatrix{
0 \ar[r] & A_1 \ar[r] \ar[d]_{\alpha} & B_1 \ar[r] \ar[d]_{\beta} & C_1 \ar[r] \ar[d]^{\gamma} & 0 \\
0 \ar[r] & A_2 \ar[r] & B_2 \ar[r] & C_2 \ar[r] & 0
}
\]
be a commutative diagram between short exact sequences of abelian groups. If $\alpha$ is an isomorphism, then the second square is both a pullback and pushout square. If $\gamma$ is an isomorphism, then the first square is both a pullback and pushout square. 
\end{remark}

\begin{defin}
For an abelian group $A$ we define the quadratic module 
\[
I(A) := (A \rightarrow 0 \rightarrow A).
\]
\end{defin}
\noindent
Note that $I$ defines a functor, and for any abelian group $G$ we have isomorphisms $G \otimes_{\Z} I(A) \cong G \otimes A$ and $\HomZ(G,I(A)) \cong \Hom(G,A)$. 

\begin{prop} \label{prop:qfp-seq-symm}
For every symmetric form parameter $Q = (Q_e, {\rm h}, {\rm p})$ there is a sequence of quadratic modules 
\[
\xymatrix{
Q_+ \ar[r] & Q \ar[r] & I(SQ),
}
\]
which is natural in $Q$.
\end{prop}

\begin{proof}
The first map is the unique morphism of form parameters $Q_+ \rightarrow Q$ (which is given by ${\rm p} \colon \Z \rightarrow Q_e$) and the second one is given by the quotient map $\pi \colon Q_e \rightarrow SQ$. For naturality we need to check that a morphism $\alpha \colon Q \rightarrow Q'$ induces a commutative diagram
\[
\xymatrix{
Q_+ \ar@{=}[d] \ar[r] & Q \ar[d]_-{\alpha} \ar[r] & I(SQ) \ar[d]^-{I(S\alpha)} \\
Q_+ \ar[r] & Q' \ar[r] & I(SQ').
}
\]
The first square commutes by the definition of a morphism of form parameters. The second square commutes because $S\alpha \colon SQ \rightarrow SQ'$ is induced by $\alpha \colon Q_e \rightarrow Q'_e$.
\end{proof}

\begin{prop} \label{p:qt-pback}
a) For every symmetric form parameter $Q = (Q_e, {\rm h}, {\rm p})$ there is a commutative diagram of quadratic modules
\[
\xymatrix{
Q_+ \ar@{=}[d] \ar[r] & Q \ar[d] \ar[r] & I(SQ) \ar[d] \\
Q_+ \ar[r] & Q^+ \ar[r] & I(\Z_2) .
}
\]

b) Given an abelian group $G$, if we apply the functor $G \otimes_{\Z} {-}$ to the above diagram, then we get a commutative diagram of abelian groups 
\[
\xymatrix{
0 \ar[r] & S^2(G) \ar@{=}[d] \ar[r] & G \otimes_{\Z} Q \ar[d] \ar[r] & G \otimes SQ \ar[d] \ar[r] & 0 \\
0 \ar[r] & S^2(G) \ar[r] & \Gamma(G) \ar[r] & G \otimes \Z_2 \ar[r] & 0
}
\]
consisting of short exact sequences. In particular the second square in this diagram is a pullback square.

c) Given a free abelian group $G$, if we apply the functor $\HomZ(G, {-})$ to the diagram in part a), then we get a commutative diagram of abelian groups 
\[
\xymatrix{
0 \ar[r] & \HomZ(G,Q_+) \ar@{=}[d] \ar[r] & \HomZ(G,Q) \ar[d] \ar[r] & \Hom(G,SQ) \ar[d] \ar[r] & 0 \\
0 \ar[r] & \HomZ(G,Q_+) \ar[r] & \HomZ(G,Q^+) \ar[r] & \Hom(G,\Z_2) \ar[r] & 0
}
\]
consisting of short exact sequences. In particular the second square in this diagram is a pullback square.
\end{prop}

\begin{remark}
Part b) of Proposition~\ref{p:qt-pback} was proven by Baues \cite[Remark 2.10 (4)]{Bau1} in the special case $Q = \Z^P$ (when $G \otimes_{\Z} \Z^P \cong P^2(G)$ and $G \otimes S\Z^P \cong G$).
\end{remark}

\begin{proof}
a) We get this by applying Proposition \ref{prop:qfp-seq-symm} to the unique morphism of form parameters $Q \rightarrow Q^+$, which is given by ${\rm h} \colon Q_e \rightarrow \Z$. Note that $SQ^+ = \Z_2$ and the map $I(SQ) \rightarrow I(\Z_2)$ is $I(v_Q)$, because $v_Q$ is induced by ${\rm h}$ (see also Lemma \ref{l:Wu-pb}).

b) The exactness of the bottom sequence was proved by Baues. It follows from part a) that the diagram commutes. This also implies that the map $S^2(G) \rightarrow G \otimes_{\Z} Q$ is injective (because its composition with $\Id_G \otimes_{\Z} {\rm h} \colon G \otimes_{\Z} Q \rightarrow \Gamma(G)$ is the injective map $S^2(G) \rightarrow \Gamma(G)$). 

The map $S^2(G) \rightarrow G \otimes_{\Z} Q$ sends $gh$ to $[g,h] \otimes 1$, so $G \otimes_{\Z} Q / S^2(G)$ is the group generated by the symbols $g \otimes q$, linear in both $g \in G$ and $q \in Q_e$, and subject to the relation $g \otimes {\rm p}(1) = 0$, therefore $G \otimes_{\Z} Q / S^2(G) \cong G \otimes SQ$. Moreover, this isomorphism is induced by the map $G \otimes_{\Z} Q \rightarrow G \otimes SQ$, because it sends $[g,h] \otimes 1$ to $0$ and $g \otimes q$ to $g \otimes \pi(q)$. Therefore the top sequence is exact too. 

c) First we show that the bottom sequence is exact. An element $(G,\lambda,\mu)$ of $\HomZ(G,Q^+)$ is determined by the symmetric bilinear form $(G,\lambda)$ (because $\mu$ satisfies $\mu(x) = \lambda(x,x)$). The image of $(G,\lambda)$ is $f \in \Hom(G,\Z_2)$ given by $f(x) = \lambda(x,x)$ mod $2$. For a free $G$ the map $\HomZ(G,Q^+) \rightarrow \Hom(G,\Z_2)$ is surjective, because for a homomorphism $f$ and a basis $x_1, \ldots , x_k$ of $G$ we can define a $\lambda$ in the preimage of $f$ by $\lambda(x_i,x_i) = 0$ if $f(x_i)=0$, $\lambda(x_i,x_i) = 1$ if $f(x_i)=1$ and $\lambda(x_i,x_j) = 0$ if $i \neq j$. The kernel of this map consists of even forms on $G$, which are exactly the elements of $\HomZ(G,Q_+)$ (because if $(G,\lambda,\mu') \in \HomZ(G,Q_+)$, then $2\mu'(x) = \lambda(x,x)$). Moreover, the map $\HomZ(G,Q_+) \rightarrow \HomZ(G,Q^+)$ sends $(G,\lambda,\mu')$ to $(G,\lambda,\mu)$ (where $2\mu'(x) = \lambda(x,x) = \mu(x)$), therefore the bottom sequence is exact.

By Lemma \ref{l:gmtrc} the second square is a pullback square (note that $\HomZ(G,Q^+ \oplus SQ) \cong \HomZ(G,Q^+) \oplus \Hom(G,SQ)$). It follows that the map $\HomZ(G,Q) \rightarrow \Hom(G,SQ)$ is surjective too, and its kernel is also $\HomZ(G,Q_+)$, therefore the top sequence is exact too.
\end{proof}

\begin{corollary} \label{cor:qfunctors-ses-symm}
For every symmetric form parameter $Q$ and abelian group $G$ there is a short exact sequence
\[
\xymatrix{
0 \ar[r] & S^2(G) \ar[r] & G \otimes_{\Z} Q \ar[r] & G \otimes SQ \ar[r] & 0 .
}
\]
For every symmetric form parameter $Q$ and free abelian group $G$ there is a short exact sequence
\[
\xymatrix{
0 \ar[r] & \HomZ(G,Q_+) \ar[r] & \HomZ(G,Q) \ar[r] & \Hom(G,SQ) \ar[r] & 0 .
}
\]
These sequences are natural in $Q$ and $G$.
\end{corollary}

\begin{proof}
These sequences are obtained from the sequence of Proposition \ref{prop:qfp-seq-symm} by applying $G \otimes_{\Z} {-}$ and $\HomZ(G, {-})$ respectively. So they are natural in $Q$ by Proposition \ref{prop:qfp-seq-symm} and they are natural in $G$, because $G \otimes_{\Z} {-}$ and $\HomZ(G, {-})$ are bifunctors. Exactness was proven in Proposition \ref{p:qt-pback}. 
\end{proof}

\begin{prop} \label{prop:qfp-seq-antisymm}
For every anti-symmetric form parameter $Q = (Q_e, {\rm h}, {\rm p})$ there is a sequence of quadratic modules 
\[
\xymatrix{
I(Q_e) \ar[r] & Q \ar[r] & Q^- ,
}
\]
which is natural in $Q$.
\end{prop}

\begin{proof}
The first map is given by $\Id_{Q_e}$ and the second one is the unique morphism of form parameters $Q \rightarrow Q^-$, both of these are natural.
\end{proof}

\begin{prop} \label{prop:qfunctors-ses-antisymm}
For every anti-symmetric form parameter $Q = (Q_e, {\rm h}, {\rm p})$ and abelian group $G$ there are short exact sequences
\[
\xymatrix{
0 \ar[r] & G \otimes Q_e \ar[r] & G \otimes_{\Z} Q \ar[r] & \Lambda(G) \ar[r] & 0
} 
\]
and
\[
\xymatrix{
0 \ar[r] & \Hom(G,Q_e) \ar[r] & \HomZ(G,Q) \ar[r] & \HomZ(G,Q^-) \ar[r] & 0 .
} 
\]
These sequences are natural in $Q$ and $G$.
\end{prop}

\begin{proof}
These sequences are obtained from the sequence of Proposition \ref{prop:qfp-seq-antisymm} by applying $G \otimes_{\Z} {-}$ and $\HomZ(G, {-})$ respectively, so naturality follows as in Corollary \ref{cor:qfunctors-ses-symm}. It remains to prove the exactness of the sequences:

a) The map $G \otimes_{\Z} Q \rightarrow \Lambda(G)$ sends $g \otimes q$ to $0$ and $[g,h] \otimes a$ to $ag \wedge h$, hence it is surjective. The map $G \otimes Q_e \rightarrow G \otimes_{\Z} Q$ sends $g \otimes q$ to $g \otimes q$, so its cokernel is generated by the symbols $[g,h] \otimes a$ (linear in all variables) subject to the relation $[g,g] \otimes a = 0$, hence it is isomorphic to $\Lambda(G)$. This shows that the sequence is exact at $G \otimes_{\Z} Q$. If $G$ is written as the direct sum of cyclic groups with generators $g_1, \ldots , g_k$, then we can define a map $G \otimes_{\Z} Q \rightarrow G \otimes Q_e$ by $g \otimes q \mapsto g \otimes q$ (for all $g$), $[g_i,g_j] \otimes a \mapsto 0$ if $i \neq j$ and $[g_i,g_i] \otimes a \mapsto g_i \otimes {\rm p}(a)$. 
This map is a left inverse of $G \otimes Q_e \rightarrow G \otimes_{\Z} Q$, so the sequence is also exact at $G \otimes Q_e$.  

b) The homomorphism $\HomZ(G,Q) \rightarrow \HomZ(G,Q^-)$ is the forgetful map $(G,\lambda,\mu) \mapsto (G,\lambda)$. It is surjective, because for every bilinear form $\lambda$ on $G$ with $\lambda(x,x)=0$ we can define a refinement $\mu \colon G \rightarrow Q_e$. Namely, if $G$ is written as the direct sum of cyclic groups with generators $g_1, \ldots , g_k$, then let $\mu(\sum_i a_ig_i) = {\rm p}(\sum_{i < j} a_ia_j\lambda(g_i,g_j))$. This is well-defined (because if $g_i$ has infinite order, then $a_i \in \Z$ is well-defined, and $\lambda(g_i,g_j)=0$ otherwise) and we can check that it satisfies $\mu(x+y)=\mu(x)+\mu(y)+{\rm p}(\lambda(x,y))$ (using that ${\rm p} \circ \lambda \colon G \times G \rightarrow Q_e$ is symmetric, because $2{\rm p}=0$). 

The kernel of the map $\HomZ(G,Q) \rightarrow \HomZ(G,Q^-)$ consists of forms $(G,\lambda,\mu)$ with $\lambda=0$ and $\mu \colon G \rightarrow Q_e$ a homomorphism. Therefore the sequence is exact.
\end{proof}

\begin{corollary} \label{cor:qt-pushout}
a) For every anti-symmetric form parameter $Q = (Q_e, {\rm h}, {\rm p})$ there is a commutative diagram of quadratic modules
\[
\xymatrix{
I(\Z_2) \ar[d] \ar[r] & Q_- \ar[d] \ar[r] & Q^- \ar@{=}[d] \\
I(Q_e) \ar[r] & Q \ar[r] & Q^-
}
\]

b) Given an abelian group $G$, if we apply the functor $G \otimes_{\Z} {-}$ to the above diagram, then we get a commutative diagram of abelian groups 
\[
\xymatrix{
0 \ar[r] & G \otimes \Z_2 \ar[d] \ar[r] & \Lambda_1(G) \ar[d] \ar[r] & \Lambda(G) \ar@{=}[d] \ar[r] & 0 \\
0 \ar[r] & G \otimes Q_e \ar[r] & G \otimes_{\Z} Q \ar[r] & \Lambda(G) \ar[r] & 0
}
\]
consisting of short exact sequences. In particular the first square in this diagram is a pushout square.

c) Given an abelian group $G$, if we apply the functor $\HomZ(G, {-})$ to the diagram in part a), then we get a commutative diagram of abelian groups 
\[
\xymatrix{
0 \ar[r] & \Hom(G,\Z_2) \ar[d] \ar[r] & \HomZ(G,Q_-) \ar[d] \ar[r] & \HomZ(G,Q^-) \ar@{=}[d] \ar[r] & 0 \\
0 \ar[r] & \Hom(G,Q_e) \ar[r] & \HomZ(G,Q) \ar[r] & \HomZ(G,Q^-) \ar[r] & 0
}
\]
consisting of short exact sequences. In particular the first square in this diagram is a pushout square.
\end{corollary}

\begin{proof}
This follows from Propositions \ref{prop:qfp-seq-antisymm} and \ref{prop:qfunctors-ses-antisymm}, applied to the unique morphism of form parameters $Q_- \rightarrow Q$ given by $v'_Q \colon \Z_2 \rightarrow Q_e$. 
\end{proof}

\subsubsection{An aside on the quadratic tensor product and $\HomZ$}

The quadratic tensor product and $\HomZ$ can be described in terms of a tensor product operation and morphisms of quadratic modules, respectively. These notions, in turn, are special cases of tensor products and morphisms of functors.

We start by observing that for an abelian group $A$, the group $P^2(A) = A \otimes_{\Z} \Z^P$ is generated by symbols $\overline{[a,b]}$ (bilinear in $a,b \in A$) and $\overline{a}$, satisfying $\overline{a+b} = \overline{a} + \overline{b} + \overline{[a,b]}$. Define the quadratic module 
\[
J(A) = (P^2(A) \xrightarrow{\overline{\h}} A \otimes A \xrightarrow{\overline{\p}} P^2(A))
\]
where $\overline{\h}(\overline{[a,b]}) = a \otimes b + b \otimes a$, $\overline{\h}(\overline{a}) = a \otimes a$ and $\overline{\p}(a \otimes b) = \overline{[a,b]}$. Then we have 
\[
\Hom(J(A),Q) \cong \HomZ(A,Q)
\]
where a morphism $(F \colon A \otimes A \to Q_{ee}, G \colon P^2(A) \to Q_e)$ corresponds to the form $(A,\lambda,\mu)$ with $\lambda(a,b) = F(a \otimes b)$ and $\mu(a) = G(\overline{a})$.

To get a similar isomorphism involving the quadratic tensor product, we define the tensor product $Q \otimes Q'$ of quadratic modules $Q = (Q_{ee},Q_e,\h,\p)$ and $Q' = (Q'_{ee},Q'_e,\h',\p')$ as follows. It is the abelian group generated by the (bilinear) symbols $q_1 \otimes q'_1$ and $q_2 \otimes q'_2$ (where $q_1 \in Q_e$, $q_2 \in Q_{ee}$, $q'_1 \in Q'_e$, $q'_2 \in Q'_{ee}$), subject to the relations $\h(q_1) \otimes q'_2 = q_1 \otimes \p'(q'_2)$ and $\p(q_2) \otimes q'_1 = q_2 \otimes \h'(q'_1)$. Then 
\[
J(A) \otimes Q \cong A \otimes_{\Z} Q
\]
with $\overline{a} \otimes q_1$, $\overline{[a,b]} \otimes q_1$ and $(a \otimes b) \otimes q_2$ corresponding to $a \otimes q_1$, $[a,b] \otimes \h(q_1)$ and $[a,b] \otimes q_2$, respectively. 

The notions of a tensor product and a morphism of quadratic modules are themselves special cases of corresponding notions for functors. To see this, we let $\CC$ be the category on two objects, $X$ and $Y$, generated by two morphisms, $\alpha : X \to Y$ and $\beta : Y \to X$. Then a quadratic module $Q$ determines a functor $F_Q : \CC \to \Ab$, $F_Q(X) = Q_{ee}$, $F_Q(Y) = Q_e$, $F_Q(\alpha) = \p$, $F_Q(\beta) = \h$. At the same time, any quadratic module $Q'$ determines a functor $G_{Q'} : \CC^{\text{op}} \to \Ab$, $G_{Q'}(X) = Q'_{ee}$, $G_{Q'}(Y) = Q'_e$, $G_{Q'}(\alpha) = \h'$, $G_{Q'}(\beta) = \p'$. The tensor product $Q \otimes Q'$ defined above coincides with the tensor product of the functors $F_Q$ and $G_{Q'}$. Of course, a morphism $Q \to Q'$ between quadratic modules is just a morphism $F_Q \to F_{Q'}$ between functors (i.e.\ a natural transformation).

\subsection{Reduced Witt groups and the split norm} \label{ss:RWG}
%%%%%%%%%%%%%%%%%%%%%%%%%%%%%%%%%%%%%%%%%%%%%%%%%%%%%%%%
Now we consider a split form parameter $Q \oplus G$, where $Q$ is a form parameter and $G$ is an abelian group.  This subsection is about the problem of computing $W_0(Q \oplus G)$, assuming that $W_0(Q)$ is known.  
We have seen that taking Witt groups defines a functor $W_0 \colon \FP \to \AbF$.
Now there are obvious morphisms $i\colon Q \to Q \oplus G$ and 
$r\colon Q \oplus G \to Q$ with $ri = {\rm Id}_Q$, 
and so $W_0(Q)$ is naturally identified with 
a summand of $W_0(Q \oplus G)$.

\begin{defin} \label{def:red-Witt}
We define the \emph{reduced Witt group} of $Q \oplus G$ by
\[
W_0^Q(G) := {\rm Ker}(W_0(r)) \leq W_0(Q \oplus G).
\]
\end{defin}
\noindent
Evidently $W_0(Q \oplus G) \cong W_0(Q) \oplus W_0^Q(G)$ and our problem reduces to computing $W_0^Q(G)$.  

A split morphism $\alpha = \alpha_Q \oplus \alpha_G \colon Q_1 \oplus G_1 \to Q_2 \oplus G_2$ induces a commutative diagram
\[
\xymatrix{
Q_1 \ar[d]_-{\alpha_Q} \ar[r]^-{i_1} & Q_1 \oplus G_1 \ar[d]_-{\alpha} \ar[r]^-{r_1} & Q_1 \ar[d]^-{\alpha_Q} \\
Q_2 \ar[r]^-{i_2} & Q_2 \oplus G_2 \ar[r]^-{r_2} & Q_2
}
\]
of form parameters. If we apply the functor $W_0$ to this, we get an analogous commutative diagram between the Witt groups. Therefore the map $W_0(\alpha)$ is of the form $W_0(\alpha_Q) \oplus f \colon W_0(Q_1) \oplus W_0^{Q_1}(G_1) \rightarrow W_0(Q_2) \oplus W_0^{Q_2}(G_2)$ for some homomorphism $f \colon W_0^{Q_1}(G_1) \rightarrow W_0^{Q_2}(G_2)$, showing that $W_0^Q(G)$ is a functor in both $Q$ and $G$.

To study reduced Witt groups, we introduce the \emph{split norm} of a nonsingular $(Q \oplus G)$-form, for any form parameter $Q$ and a finitely generated abelian group $G$, which is an element of $G \tensor_{\Z} Q$. 
The split norm is invariant under isomorphism of $(Q\oplus G)$-forms (and even Witt equivalence, see Proposition \ref{p:tilde-F}), and hence induces a function
\[
F \colon \HomZns(Q \oplus G) \rightarrow G \tensor_{\Z} Q.
\] 
We will define this invariant in two different ways, and then show that the two definitions are equivalent. Then, in Theorem \ref{thm:main}, we will use $F$ to identify $W_0^Q(G)$ with $G \tensor_{\Z} Q$.

\begin{defin}[Split norm, first definition] \label{def:F-1}
Suppose that $(X, \lambda, (\mu_Q,\mu_G))$ is a nonsingular $(Q \oplus G)$-form. Choose any basis $x_1, \ldots , x_k$ of $X$ and let $y_1, \ldots , y_k$ be the basis of $X$ which is dual to $x_1, \ldots , x_k$ with respect to $\lambda$ (i.e.\ $\lambda(x_i,y_i) = 1$ and $\lambda(x_i,y_j) = 0$ if $i \neq j$). We set
\[
F(X, \lambda, (\mu_Q,\mu_G)) := \sum_{i<j} \left[ \mu_G(x_i),\mu_G(x_j) \right] \tensor \lambda(y_i,y_j) + \sum_i \mu_G(x_i) \tensor \mu_Q(y_i) \in G \tensor_{\Z} Q.
\]
\end{defin}

\begin{prop} \label{prop:F-well-def}
$F(X, \lambda, (\mu_Q,\mu_G))$ is well-defined; i.e.\ it does not depend on the choice of the basis.
Moreover, if $(X, \lambda, (\mu_Q, \mu_G)) \cong (X', \lambda', (\mu'_Q, \mu'_G))$, then $F(X, \lambda, (\mu_Q, \mu_G)) = F(X', \lambda', (\mu'_Q, \mu'_G))$, so that $F$ induces a well-defined function $F \colon \HomZns(Q \oplus G) \rightarrow G \tensor_{\Z} Q$ on the set of isomorphism classes.
\end{prop}

\begin{proof}
If suffices to check that the value of $F$ does not change when we apply elementary transformations to the basis $x_1, \ldots , x_k$.

If we swap $x_1$ and $x_2$, then $y_1$ and $y_2$ also get swapped. In $F$ the only change is that $[\mu_G(x_1),\mu_G(x_2)] \tensor \lambda(y_1,y_2)$ is replaced with $[\mu_G(x_2),\mu_G(x_1)] \tensor \lambda(y_2,y_1)$, and these are equal, because $[\cdot , \cdot]$ and $\lambda$ are both $\epsilon_Q$-symmetric.

If we replace $x_1$ with $-x_1$, then $y_1$ is also replaced with $-y_1$. In the first part of $F$ the terms $[\mu_G(x_1),\mu_G(x_j)] \tensor \lambda(y_1,y_j)$ for $j \geq 2$ are replaced with 
\[[\mu_G(-x_1),\mu_G(x_j)] \tensor \lambda(-y_1,y_j) = [\mu_G(x_1),\mu_G(x_j)] \tensor \lambda(y_1,y_j).\]
In the second part $\mu_G(x_1) \tensor \mu_Q(y_1)$ is replaced with $\mu_G(-x_1) \tensor \mu_Q(-y_1)$, we will show that these are equal. We have $0 = \mu_Q(y_1 - y_1) = \mu_Q(y_1) + \mu_Q(-y_1) + {\rm p}(\lambda(y_1,-y_1))$ and also $\lambda(y_1,-y_1) = -\lambda(-y_1,-y_1) = -{\rm h}(\mu_Q(-y_1))$.
Hence
\[
0 = \mu_G(x_1) \tensor \mu_Q(y_1 - y_1) = \mu_G(x_1) \tensor \mu_Q(y_1) + \mu_G(x_1) \tensor \mu_Q(-y_1) - \mu_G(x_1) \tensor {\rm p}({\rm h}(\mu_Q(-y_1))).  
\]
We also have $[\mu_G(x_1), \mu_G(-x_1)] \tensor {\rm h}(\mu_Q(-y_1)) = -[\mu_G(x_1), \mu_G(x_1)] \tensor {\rm h}(\mu_Q(-y_1)) = -\mu_G(x_1) \tensor {\rm p}({\rm h}(\mu_Q(-y_1)))$, hence
\[
0 = (\mu_G(x_1) + \mu_G(-x_1)) \tensor \mu_Q(-y_1) = \mu_G(x_1) \tensor \mu_Q(-y_1) + \mu_G(-x_1) \tensor \mu_Q(-y_1) - \mu_G(x_1) \tensor {\rm p}({\rm h}(\mu_Q(-y_1))) .
\]
By comparing these two equalities we see that $\mu_G(x_1) \tensor \mu_Q(y_1) = \mu_G(-x_1) \tensor \mu_Q(-y_1)$.

Finally if we replace $x_1$ with $x_1 + x_2$, then $y_2$ is replaced with $y_2-y_1$. In $F$ the terms 
\[ [\mu_G(x_1),\mu_G(x_j)] \tensor \lambda(y_1,y_j) + [\mu_G(x_2),\mu_G(x_j)] \tensor \lambda(y_2,y_j)\] 
are replaced by
\begin{multline*}
[\mu_G(x_1 + x_2),\mu_G(x_j)] \tensor \lambda(y_1,y_j) + [\mu_G(x_2),\mu_G(x_j)] \tensor \lambda(y_2 - y_1,y_j)= \\ [\mu_G(x_1),\mu_G(x_j)] \tensor \lambda(y_1,y_j) + [\mu_G(x_2),\mu_G(x_j)] \tensor \lambda(y_2,y_j)
\end{multline*}
for every $j \geq 3$. Further, $[\mu_G(x_1),\mu_G(x_2)] \tensor \lambda(y_1,y_2) + \mu_G(x_1) \tensor \mu_Q(y_1) + \mu_G(x_2) \tensor \mu_Q(y_2)$ is replaced by $[\mu_G(x_1 + x_2),\mu_G(x_2)] \tensor \lambda(y_1,y_2-y_1) + \mu_G(x_1 + x_2) \tensor \mu_Q(y_1) + \mu_G(x_2) \tensor \mu_Q(y_2-y_1)$. The difference is 
\begin{multline*}
[\mu_G(x_2),\mu_G(x_2)] \tensor \lambda(y_1,y_2) + [\mu_G(x_1),\mu_G(x_2)] \tensor \lambda(y_1,-y_1) + [\mu_G(x_2),\mu_G(x_2)] \tensor \lambda(y_1,-y_1) + \\
+ \mu_G(x_2) \tensor \mu_Q(y_1) + [\mu_G(x_1),\mu_G(x_2)] \tensor {\rm h}(\mu_Q(y_1)) + \mu_G(x_2) \tensor \mu_Q(-y_1) + \mu_G(x_2) \tensor {\rm p}(\lambda(y_2,-y_1)) = \\
= [\mu_G(x_2),\mu_G(x_2)] \tensor \lambda(y_1,-y_1) + \mu_G(x_2) \tensor \mu_Q(y_1) + \mu_G(x_2) \tensor \mu_Q(-y_1) = \\
= \mu_G(x_2) \tensor {\rm p}(\lambda(y_1,-y_1)) + \mu_G(x_2) \tensor \mu_Q(y_1) + \mu_G(x_2) \tensor \mu_Q(-y_1) = \mu_G(x_2) \tensor \mu_Q(y_1-y_1) = 0,
\end{multline*}
where we used that ${\rm h}(\mu_Q(y_1)) = -\lambda(y_1,-y_1)$ and ${\rm p}(\lambda(y_2,-y_1)) = -{\rm p}(\lambda(y_1,y_2))$.

Hence the value of $F$ does not change when we change the basis. 

If $(X, \lambda, (\mu_Q, \mu_G)) \cong (X', \lambda', (\mu'_Q, \mu'_G))$, then we can see that $F$ takes the same value on both forms by choosing compatible bases in $X$ and $X'$. Therefore the split norm is an isomorphism invariant.
\end{proof}

We present a second definition of the split norm, which does not rely on a choice of basis.
We do this for completeness, added understanding of the split norm and possible future applications.
In order to present the second definition of the split norm, we need the following result of Baues
\cite[Lemma (8.3)]{Bau2} (see also \cite[Lemma 5.4]{Bau1}):  If $X$ is torsion free, then there is a natural isomorphism
\[ \chi \colon X \tensor_\Z Q \to \HomZ(X^*, Q), \]
which is defined as follows.
For $x,y \in X$, $q \in Q_e$ and $n \in \Z$
\[
\chi(x \tensor q) = (X^*, \lambda_{x \tensor q}, \mu_{x \tensor q})
\qquad \text{and} \qquad
\chi([x, y] \tensor n) = (X^*, \lambda_{[x, y] \tensor n}, \mu_{[x, y] \tensor n}),
\]
where for $a, b \in X^*$, $\lambda_{x \tensor q}(a, b) = a(x)\h(q)b(x)$ and $\mu_{x \tensor q}(a) = a(x) \cdot q$ using the $\Z$-action given in \eqref{eq:Zactn}, and $\lambda_{[x, y] \tensor n}(a, b) = a(x)nb(y) + \epsilon_Q a(y)nb(x)$ and $\mu_{[x, y] \tensor n}(a) = \p(a(x)na(y))$.

\begin{defin}[Split norm, second definition] \label{def:F-2}
Let $\ul \mu = (X, \lambda, (\mu_Q,\mu_G))$ be a nonsingular $(Q \oplus G)$-form.
We set 
\[ 
F\bigl( \ul \mu \bigr) := (\mu_G)_*(\chi^{-1}(\ul \mu_Q^{-1})) \in G \tensor_{\Z} Q,
\]
where
\begin{compactitem}[$\bullet$]
\item $\ul \mu_Q := r_*(\ul \mu) = (X, \lambda, \mu_Q) \in \HomZns(X, Q)$;
\item $\lambda^{-1} \colon X^* \times X^* \rightarrow \Z$ denotes the adjoint of $\hat \lambda^{-1} \colon X^* \rightarrow X = (X^*)^*$, where $\hat \lambda \colon X \to X^*$ is the adjoint of $\lambda$;
\item $\ul \mu_Q^{-1} := (\hat \lambda^{-1})^*(\ul \mu_Q) = (X^*, \lambda^{-1}, \mu_Q \circ \hat \lambda^{-1}) \in \HomZns(X^*, Q)$;
\item $(\mu_G)_* \colon X \tensor_{\Z} Q \rightarrow G \tensor_{\Z} Q$ is the map induced by $\mu_G$.
\end{compactitem}
That is, $F(\ul \mu) = C_{\ul \mu}(\ul \mu)$, where $C_{\ul \mu}$ is the composition 
\[
C_{\ul \mu} \colon \HomZns(X, Q \oplus G) \xra{r_*} \HomZns(X, Q) \xra{(\hat \lambda^{-1})^*} \HomZns(X^*, Q) \xra{\chi^{-1}} X \tensor_\Z Q \xra{(\mu_G)_*} G \tensor_{\Z} Q .
\]
\end{defin}

\begin{lem} \label{lem:F-same}
The two definitions of the split norm $F$ coincide.
\end{lem}

\begin{proof}
Suppose that $\ul \mu = (X, \lambda, (\mu_Q,\mu_G))$ is a nonsingular $(Q \oplus G)$-form. Let $\ul \mu_X = (X, \lambda, (\mu_Q,\Id_X)) \in \HomZns(X, Q \oplus X)$, so that $\ul \mu = (\Id_Q \oplus \mu_G)_*(\ul \mu_X)$. Then it is enough to show that the two definitions of $F(\ul \mu_X) \in X \tensor_\Z Q$ coincide, because both versions of $F(\ul \mu)$ can be obtained from that by applying $(\mu_G)_*$.

Let $x_1, \ldots , x_k$ be a basis of $X$, let $y_1, \ldots , y_k$ be the dual basis of $X$ with respect to $\lambda$, and let $x_1^*, \ldots, x_k^*$ be the dual basis of $X^*$. We need to show that 
\[
\chi \left( \sum_{i<j} \left[ x_i, x_j \right] \tensor \lambda(y_i,y_j) + \sum_i x_i \tensor \mu_Q(y_i) \right) = (X^*, \lambda^{-1}, \mu_Q \circ \hat \lambda^{-1}) .
\]

For brevity, we will denote the left hand side by $(X^*, \ol \lambda, \ol \mu)$.
By the definition of $\chi$ above we have
\[
\begin{aligned}
\ol \lambda &= \sum_{i<j} \lambda_{[x_i, x_j] \tensor \lambda(y_i,y_j)} + \sum_i \lambda_{x_i \tensor \mu_Q(y_i)} \\
\ol \mu &= \sum_{i<j} \mu_{[x_i, x_j] \tensor \lambda(y_i,y_j)} + \sum_i \mu_{x_i \tensor \mu_Q(y_i)} 
\end{aligned}
\]
Therefore $\ol \lambda(x_i^*, x_j^*)$ is $\lambda(y_i,y_j)$ if $i < j$, $\h(\mu_Q(y_i)) = \lambda(y_i, y_i)$ if $i = j$ and $\epsilon_Q \lambda(y_j,y_i) = \lambda(y_i,y_j)$ if $i > j$. That is, $\ol \lambda(x_i^*, x_j^*) = \lambda(y_i,y_j)$ in all cases. We also have $\ol \mu(x_i^*) = 1 \cdot \mu_Q(y_i) = \mu_Q(y_i)$ for all $i$. 

On the right hand side, since $\hat \lambda(y_i) = x_i^*$, we have $\lambda^{-1}(x_i^*, x_j^*) = \lambda(y_i, y_j)$ for all $i, j$, therefore $\lambda^{-1} = \ol \lambda$. Also $\mu_Q \circ \hat \lambda^{-1}(x_i^*) = \mu_Q(y_i) = \ol \mu(x_i^*)$ for all $i$, and since the quadratic refinement $\ol \mu$ is uniquely determined by its values on a basis, this means that $\mu_Q \circ \hat \lambda^{-1} = \ol \mu$, as required. 
\end{proof}

Since the second definition does not depend on the choice of a basis, and in the proof of Lemma \ref{lem:F-same} the basis $x_1, \ldots , x_k$ can be chosen arbitrarily, this argument also serves as an alternative proof for Proposition \ref{prop:F-well-def}.  From now on we will use the first definition.

\begin{prop}
$F$ is a monoid homomorphism. 
\end{prop}

\begin{proof}
Given nonsingular $(Q\oplus G)$-forms $(X, \lambda, (\mu_Q,\mu_G))$ and  $(X', \lambda', (\mu'_Q,\mu'_G))$,
let $x_1, \ldots , x_k$ be a basis of $X$ and $x'_1, \ldots , x'_l$ be a basis of $X'$ (with dual bases $y_1, \ldots , y_k$ and $y'_1, \ldots , y'_l$, respectively). Then 
\begin{multline*}
F((X, \lambda, (\mu_Q,\mu_G)) \oplus (X', \lambda', (\mu'_Q,\mu'_G))) =\\ 
F(X, \lambda, (\mu_Q,\mu_G)) + F(X', \lambda', (\mu'_Q,\mu'_G)) + \sum_{i \leq k, j \leq l} [\mu_G(x_i),\mu'_G(x'_j)] \tensor (\lambda \oplus \lambda')(y_i,y'_j) = \\ F(X, \lambda, (\mu_Q,\mu_G)) + F(X', \lambda', (\mu'_Q,\mu'_G)),
\end{multline*}
because $(\lambda \oplus \lambda')(y,y') = \lambda(y,0) + \lambda'(0,y') = 0$ for every $y \in X$ and $y' \in X'$. 
\end{proof}

\begin{prop} \label{p:tilde-F}
The map $F$ induces a homomorphism $\Tilde{F} \colon W_0(Q \oplus G) \rightarrow G \tensor_{\Z} Q$.
\end{prop}

\begin{proof}
It is enough to prove that $F$ vanishes on metabolic forms. So suppose that $(X, \lambda, (\mu_Q,\mu_G))$ is a metabolic form over $Q \oplus G$, let $L < X$ be a lagrangian and fix an identification $X \cong L \oplus L^*$. 
Let $x_1, \ldots , x_k$ be a basis of $L$, and let $\alpha_1, \ldots , \alpha_k$ be the dual basis of $L^*$. 
Then $x_1, \ldots , x_k, \alpha_1, \ldots , \alpha_k$ is a basis of $X$, and the dual basis of $X$ with respect to $\lambda$ is $y_1, \ldots , y_k, \beta_1, \ldots , \beta_k$, where $y_i = \alpha_i + z_i$ for some $z_i \in L$ and $\beta_i = \epsilon_Qx_i$. Since $L$ is a lagrangian, $\mu_G(x_i) = 0$, $\lambda(\beta_i,\beta_j) = 0$ and $\mu_Q(\beta_i) = 0$ for every $i,j$, therefore in $F(X, \lambda, (\mu_Q,\mu_G))$ every term vanishes.
\end{proof}

%%%%%%%%%%%%%%

\begin{defin} \label{d:gamma'}
We define a homomorphism $\gamma \colon G \tensor_{\Z} Q \rightarrow W_0(Q \oplus G)$ by 
\[ 
\gamma([g_1,g_2] \tensor 1) := 
\left(
\Z^2, 
\left[ \begin{array}{cc} 0 & 1 \\ \epsilon_Q & 0 \end{array} \right],
\left( \begin{array}{c} g_1 \\ g_2 \end{array} \right) 
\right)
\quad \text{and} \quad
\gamma(g \tensor q) := 
\left(
\Z^2, 
\left[ \begin{array}{cc} 0 & 1 \\ \epsilon_Q & -{\rm h}(q) \end{array} \right],
\left( \begin{array}{c} g \\ R(q) \end{array} \right) 
\right) ,
\]
where $R = \Id - {\rm p}{\rm h} \colon Q_e \rightarrow Q_e$. 
\end{defin}

\begin{prop}
The map $\gamma$ is a well-defined homomorphism.
\end{prop}

\begin{proof}
These triples denote $(Q \oplus G)$-forms, because ${\rm h}(g_1) = {\rm h}(g_2) = {\rm h}(g) = 0$ and ${\rm h}(R(q)) = -{\rm h}(q)$; 
see Remark \ref{rem:eqf-not}. 
The matrices are nonsingular, so these $(Q \oplus G)$-forms represent elements of $W_0(Q \oplus G)$. It remains to check that the values of $\gamma$ satisfy the relations coming from the relations in $G \tensor_{\Z} Q$.

The first relation we consider is $[g_1,g_2] \tensor 1 + [g_1',g_2] \tensor 1 = [g_1+g_1',g_2] \tensor 1$. Since 
\begin{multline*}
\left(
\Z^2 \left< e_1,f_1 \right>, 
\left[ \begin{array}{cc} 0 & 1 \\ \epsilon_Q & 0 \end{array} \right],
\left( \begin{array}{c} g_1 \\ g_2 \end{array} \right) 
\right)
\oplus
\left(
\Z^2 \left< e_2,f_2 \right>, 
\left[ \begin{array}{cc} 0 & 1 \\ \epsilon_Q & 0 \end{array} \right],
\left( \begin{array}{c} g_1' \\ g_2 \end{array} \right) 
\right)
=
\\
\left(
\Z^2 \left< e_1+e_2,f_1 \right>, 
\left[ \begin{array}{cc} 0 & 1 \\ \epsilon_Q & 0 \end{array} \right],
\left( \begin{array}{c} g_1+g_1' \\ g_2 \end{array} \right) 
\right)
\oplus
\left(
\Z^2 \left< e_2,f_2-f_1 \right>, 
\left[ \begin{array}{cc} 0 & 1 \\ \epsilon_Q & 0 \end{array} \right],
\left( \begin{array}{c} g_1' \\ 0 \end{array} \right) 
\right)
\sim
\\ 
\left(
\Z^2 \left< e_1+e_2,f_1 \right>, 
\left[ \begin{array}{cc} 0 & 1 \\ \epsilon_Q & 0 \end{array} \right],
\left( \begin{array}{c} g_1+g_1' \\ g_2 \end{array} \right) 
\right)
\end{multline*}
(where $\Z^2 \left< e_1,f_1 \right>$ denotes a copy of the group $\Z^2$ with standard basis $e_1,f_1$), we get that 
\[ \gamma([g_1,g_2] \tensor 1) + \gamma([g_1',g_2] \tensor 1) = \gamma([g_1+g_1',g_2] \tensor 1),\]
i.e.\ the corresponding relation is satisfied by the values of $\gamma$. 
Similarly, we get that 
\[ \gamma([g_1,g_2] \tensor 1) + \gamma([g_1,g_2'] \tensor 1) = \gamma([g_1,g_2+g_2'] \tensor 1).\]
The above calculation also shows that for any $n \in \Z$ we have
\[
\gamma([g_1,g_2] \tensor n) = 
\left(
\Z^2, 
\left[ \begin{array}{cc} 0 & 1 \\ \epsilon_Q & 0 \end{array} \right],
\left( \begin{array}{c} ng_1 \\ g_2 \end{array} \right) 
\right).
\] 

We also have $\gamma(g \tensor q) + \gamma(g \tensor q') = \gamma(g \tensor (q+q'))$, because
\begin{multline*}
\left(
\Z^2 \left< e_1,f_1 \right>, 
\left[ \begin{array}{cc} 0 & 1 \\ \epsilon_Q & -{\rm h}(q) \end{array} \right],
\left( \begin{array}{c} g \\ R(q) \end{array} \right) 
\right)
\oplus
\left(
\Z^2 \left< e_2,f_2 \right>, 
\left[ \begin{array}{cc} 0 & 1 \\ \epsilon_Q & -{\rm h}(q') \end{array} \right],
\left( \begin{array}{c} g \\ R(q') \end{array} \right) 
\right)
=
\\
\left(
\Z^2 \left< e_1,f_1+f_2 \right>, 
\left[ \begin{array}{cc} 0 & 1 \\ \epsilon_Q & -{\rm h}(q+q') \end{array} \right],
\left( \begin{array}{c} g \\ R(q+q') \end{array} \right) 
\right)
\oplus
\\
\left(
\Z^2 \left< e_2-e_1,f_2+{\rm h}(q')e_2 \right>, 
\left[ \begin{array}{cc} 0 & 1 \\ \epsilon_Q & \epsilon_Q{\rm h}(q') \end{array} \right],
\left( \begin{array}{c} 0 \\ q' + {\rm h}(q')g \end{array} \right) 
\right)
\sim
\\ 
\left(
\Z^2 \left< e_1,f_1+f_2 \right>, 
\left[ \begin{array}{cc} 0 & 1 \\ \epsilon_Q & -{\rm h}(q+q') \end{array} \right],
\left( \begin{array}{c} g \\ R(q+q') \end{array} \right) 
\right).
\end{multline*}
\noindent
Next, we have $\gamma([g,g] \tensor 1) = \gamma(g \tensor {\rm p}(1))$, because 
\[
\left(
\Z^2 \left< e,f \right>, 
\left[ \begin{array}{cc} 0 & 1 \\ \epsilon_Q & 0 \end{array} \right],
\left( \begin{array}{c} g \\ g \end{array} \right) 
\right) 
= 
\left(
\Z^2 \left< e,f-e \right>, 
\left[ \begin{array}{cc} 0 & 1 \\ \epsilon_Q & -(1+\epsilon_Q) \end{array} \right],
\left( \begin{array}{c} g \\ -{\rm p}(1) \end{array} \right) 
\right)
\]
and $-(1+\epsilon_Q) = -{\rm h}{\rm p}(1)$ and $-{\rm p}(1) = R({\rm p}(1))$. 

Finally $\gamma(g \tensor q) + \gamma(g' \tensor q) + \gamma([g,g'] \tensor {\rm h}(q)) = \gamma((g+g') \tensor q)$, because
\begin{multline*}
\left(
\Z^2 \left< e_1,f_1 \right>, 
\left[ \begin{array}{cc} 0 & 1 \\ \epsilon_Q & -{\rm h}(q) \end{array} \right],
\left( \begin{array}{c} g \\ R(q) \end{array} \right) 
\right)
\oplus
\left(
\Z^2 \left< e_2,f_2 \right>, 
\left[ \begin{array}{cc} 0 & 1 \\ \epsilon_Q & -{\rm h}(q) \end{array} \right],
\left( \begin{array}{c} g' \\ R(q) \end{array} \right) 
\right)
\oplus
\\
\left(
\Z^2 \left< e_3,f_3 \right>, 
\left[ \begin{array}{cc} 0 & 1 \\ \epsilon_Q & 0 \end{array} \right],
\left( \begin{array}{c} {\rm h}(q)g \\ g' \end{array} \right) 
\right)
=
\left(
\Z^2 \left< e_1+e_2,f_1 \right>, 
\left[ \begin{array}{cc} 0 & 1 \\ \epsilon_Q & -{\rm h}(q) \end{array} \right],
\left( \begin{array}{c} g+g' \\ R(q) \end{array} \right) 
\right)
\oplus
\\
\left(
\Z^4 \left< e_2,f_2-f_1-{\rm h}(q)e_1,e_3,f_3 \right>, 
\left[ \begin{array}{cccc} 0 & 1 & 0 & 0 \\ \epsilon_Q & 0 & 0 & 0 \\ 0 & 0 & 0 & 1 \\ 0 & 0 & \epsilon_Q & 0 \end{array} \right],
\left( \begin{array}{c} g' \\ -{\rm h}(q)g \\ {\rm h}(q)g \\ g' \end{array} \right) 
\right)
\sim
\\ 
\left(
\Z^2 \left< e_1+e_2,f_1 \right>, 
\left[ \begin{array}{cc} 0 & 1 \\ \epsilon_Q & -{\rm h}(q) \end{array} \right],
\left( \begin{array}{c} g+g' \\ R(q) \end{array} \right) 
\right),
\end{multline*}
where we used that $(\epsilon_Q-1){\rm h}(q) = 0$ (because ${\rm h} = 0$ if $\epsilon_Q=-1$) and that $\left< e_2-f_3,f_2-f_1-{\rm h}(q)e_1+e_3 \right>$ is a lagrangian in $\Z^4 \left< e_2,f_2-f_1-{\rm h}(q)e_1,e_3,f_3 \right>$.

So the values of $\gamma$ satisfy the appropriate relations and therefore we have a well-defined homomorphism $G \tensor_{\Z} Q \rightarrow W_0(Q \oplus G)$. 
\end{proof}

Having shown that $\gamma$ is a homomorphism, we next determine its image.

\begin{prop} \label{p:gamma'_wd}
The homomorphism $\gamma \colon G \tensor_{\Z} Q \rightarrow W_0(Q \oplus G)$ has image $\Im(\gamma) = W_0^Q(G)$. 
\end{prop}

\begin{proof}
The natural map $\HomZns(Q \oplus G) \rightarrow \HomZns(Q)$ sends both $\gamma([g_1,g_2] \tensor 1)$ and $\gamma(g \tensor q)$ to metabolic forms, so the composition of $\gamma$ with the natural map $W_0(Q \oplus G) \rightarrow W_0(Q)$ is trivial, therefore $\Im(\gamma) \subseteq W_0^Q(G)$.

For the other direction recall that every element of $W_0^Q(G)$ has a representative $(X, \lambda, (\mu_Q,\mu_G))$ such that $(X, \lambda, \mu_Q)$ is metabolic and that every metabolic $Q$-form can be written as a sum of metabolic rank-$2$ $Q$-forms, so $(X, \lambda, (\mu_Q,\mu_G))$ is the sum of $(Q \oplus G)$-forms of the form
\[
\left(
\Z^2, 
\left[ \begin{array}{cc} 0 & 1 \\ \epsilon_Q & {\rm h}(q) \end{array} \right],
\left( \begin{array}{c} g_1 \\ q+g_2 \end{array} \right) 
\right) .
\]
These elements are contained in $\Im(\gamma)$, because
\begin{multline*}
\gamma([g_1,g_2] \tensor 1) + \gamma(g_1 \tensor R(q)) = 
\\
\left(
\Z^2 \left< e_1,f_1 \right>, 
\left[ \begin{array}{cc} 0 & 1 \\ \epsilon_Q & 0 \end{array} \right],
\left( \begin{array}{c} g_1 \\ g_2 \end{array} \right) 
\right)
\oplus
\left(
\Z^2 \left< e_2,f_2 \right>, 
\left[ \begin{array}{cc} 0 & 1 \\ \epsilon_Q & {\rm h}(q) \end{array} \right],
\left( \begin{array}{c} g_1 \\ q \end{array} \right) 
\right)
=
\\
\left(
\Z^2 \left< e_1,f_1+f_2 \right>, 
\left[ \begin{array}{cc} 0 & 1 \\ \epsilon_Q & {\rm h}(q) \end{array} \right],
\left( \begin{array}{c} g_1 \\ q+g_2 \end{array} \right) 
\right)
\oplus
\\
\left(
\Z^2 \left< e_2-e_1,f_2-{\rm h}(q)e_2 \right>, 
\left[ \begin{array}{cc} 0 & 1 \\ \epsilon_Q & -\epsilon_Q {\rm h}(q) \end{array} \right],
\left( \begin{array}{c} 0 \\ R(q)-{\rm h}(q)g_1 \end{array} \right) 
\right)
\sim
\\ 
\left(
\Z^2 \left< e_1,f_1+f_2 \right>, 
\left[ \begin{array}{cc} 0 & 1 \\ \epsilon_Q & {\rm h}(q) \end{array} \right],
\left( \begin{array}{c} g_1 \\ q+g_2 \end{array} \right) 
\right) .
\end{multline*}
Therefore $\Im(\gamma) \supseteq W_0^Q(G)$.
\end{proof}

We now relate the homomorphisms 
$\tilde F$ and $\gamma$.

\begin{prop} \label{prop:sse}
The homomorphisms $\Tilde{F}$ and $\gamma$ fit into the following split short exact sequence: 
\[
\xymatrix{
0 \ar[r] & W_0(Q) \ar@<0,4ex>[r]^-{W_0(i)} & W_0(Q \oplus G) \ar@<0,4ex>[r]^-{\Tilde{F}} \ar@{.>}@<0,4ex>[l]^-{W_0(r)} & G \tensor_{\Z} Q \ar[r] \ar@{.>}@<0,4ex>[l]^-{\gamma} & 0
}
\]
\end{prop}

\begin{proof}
We already noted at the start of this subsection that the morphisms $i \colon Q \rightarrow Q \oplus G$ and $r \colon Q \oplus G \rightarrow Q$ of form parameters are such that $W_0(r) \colon W_0(Q \oplus G) \to W_0(Q)$ splits 
$W_0(i) \colon W_0(Q) \to W_0(Q \oplus G)$.
It follows immediately from the definition that for any nonsingular $Q$-form $(X, \lambda, \mu_Q)$ we have $F(X, \lambda, (\mu_Q,0)) = 0$, therefore the composition $W_0(Q) \rightarrow W_0(Q \oplus G) \rightarrow G \tensor_{\Z} Q$ is trivial. A straightforward calculation shows that $\Tilde{F} \circ \gamma = \Id_{G \tensor_{\Z} Q}$, i.e.\ $\gamma$ is a splitting map. This implies that $\Tilde{F}$ is surjective, and that the restriction of $\Tilde{F}$ to $\Im(\gamma) = W_0^Q(G)$ is injective, so $\Ker(\Tilde{F}) = W_0(Q)$, which means that the sequence is exact. 
\end{proof}

By observing that the map $F$ is natural in $Q$ and $G$,  Proposition \ref{prop:sse} immediately gives the main result of this section, which is Theorem~\ref{t:W_0-split} from the introduction (see also Proposition \ref{p:splittable}):

\begin{thm} \label{thm:main}
Let $Q$ be a form parameter and $G$ be an abelian group. Then the restriction of $\Tilde{F}$ and $\gamma$ define isomorphisms
\[ 
W_0^Q(G) \cong G \tensor_{\Z} Q, 
\]
which are natural in $G$ and $Q$.
\qed
\end{thm}

%%%%%%%%%%%%%%%%%%%%%%%%%%%%%%%%%%%%%%%%%%%%%%%%%%%%

\subsection{Computing Witt groups via splittings of form parameters} \label{ss:compWG}
%%%%%%%%%%%%%%%%%%%%%%%%%%%%%%%%%%%%%%%%%%%%%%%%%%%%
By combining the results of Sections \ref{ss:Witt-ind} and \ref{ss:RWG} we can calculate the Witt group $W_0(P)$ for any form parameter $P$. The necessary ingredients are summarised in the following table:

\begin{figure}[H]
\begin{center}
\ul{Table 1}
\[
\begin{array}{c|cccc}
~~~Q~~~ & \Z\tensor_\Z Q & \Z_n \tensor_\Z Q & ~~~~~G \tensor _{\Z} Q~~~~~ & ~~~~~W_0(Q)~~~~~\\
\hline \hline
Q_{+} & \Z & \Z_n & S^2(G) & \Z \\  
\Z^P & \Z \oplus \Z & \Z_{(n,2)n} \oplus \Z_{n/(2,n)} & P^2(G) & \Z \oplus \Z \\
\Z^P_k & ~~~\Z \oplus \Z_{2^k} & \Z_{(n,2)n} \oplus \Z_{(n/(2,n),2^k)} & - & \qquad \! \Z \oplus \Z_{2^{k-1}} \\
Q^{+} & \Z & \Z_{(2,n)n} & \Gamma(G) & \Z \\ 
\hline 
Q_{-} & \Z_2 & \Z_{(2,n)} & \Lambda_1(G) & \Z_2 \\
\Z^{\Lambda}_k & \Z_{2^k} & \Z_{(2^k,n)} & - & 0 \\
Q^{-}  & 0 & 0 & \Lambda(G) & 0 \\
\end{array}
\]
\end{center}
\end{figure}
\noindent
Table 1 lists the values of the quadratic tensor product $G \tensor_\Z Q$ for every cyclic group $G$ and indecomposable form parameter $Q$ (where $(a,b)$ denotes the greatest common divisor of $a$ and $b$), which can be computed directly from Definition \ref{def:qtp} (or using the results of Section \ref{ss:qtp}). In some cases we also have a simple description of the functor ${-} \otimes_{\Z} Q$ (see Example \ref{ex:qtp-funct}). Finally, we include the Witt groups $W_0(Q)$ which were computed in Theorem \ref{t:wf-ind}.

Now suppose that $P$ is an arbitrary form parameter. By Theorem \ref{t:qfp-class} it has a maximal splitting $P \cong Q \oplus G$ (and the isomorphism classes of $Q$ and $G$ are uniquely determined by that of $P$). Then we have the following: 

\begin{theorem} \label{t:main2}
a) A splitting $P \cong Q \oplus G$ of a form parameter $P$
defines an isomorphism
\[
W_0(P) \cong W_0(Q \oplus G) \cong W_0(Q) \oplus W_0^Q(G) \cong W_0(Q) \oplus (G \tensor_{\Z} Q) . 
\]
b) For any pair of abelian groups $G_1$ and $G_2$, and any form parameter $Q$ we have
\[ 
(G_1 \oplus G_2) \tensor_\Z Q \cong (G_1\tensor_\Z Q) \oplus (G_2\tensor_\Z Q) \oplus (G_1 \tensor G_2) .
\] 
\end{theorem}

\begin{proof}
Part a) follows from Theorem \ref{thm:main}, and Part b) follows from \cite[Proposition 4.2 (3)]{Bau1}, since a form parameter is a quadratic module $Q$ with $Q_{ee} = \Z$.
\end{proof}

Next we give some simple examples of computations using Theorem~\ref{t:main2}.

\begin{example}
%%%%%%%%
For a positive integer $l$, define $\bar \delta(l), \delta(l) \in \Z^+$ by
\[ \bar \delta(l) = 
\begin{cases}
1 & \text{if $l$ is odd,} \\
2 & \text{if $l$ is even,}
\end{cases}
\qquad \text{and} \qquad 
\delta(l) = \bar \delta(l) l.
\]
\begin{compactenum}[a)] 
\item
There is an isomorphism $W_0(Q_+ \oplus \Z_l) \cong \Z(8\sigma^*) \oplus \Z_l$, where the torsion summand is generated by the Witt class of the form
\[\left( \Z^2, 
\left[ \begin{array}{cc}
0 & 1 \\
1 & 0
\end{array} \right], 
\left( \begin{array}{c} (0, 1) \\ (0, 1) \end{array} \right)\right).
\]

\item
There is an isomorphism $W_0(Q^+ \oplus \Z_l) \cong \Z(\sigma^*) \oplus \Z_{\delta(l)}$, where the torsion summand is generated by the Witt class of the form
\[ \bigl( (\Z, [1], (1,0)) \oplus (\Z, [-1], (-1,1) \bigr) \cong \left( \Z^2, 
\left[ \begin{array}{cc}
0 & 1 \\
1 & 1
\end{array} \right], 
\left( \begin{array}{c} (0, 1) \\ (1, 0) \end{array} \right)\right).
\]

\item There is an isomorphism $W_0(Q_- \oplus \Z_l) \cong \Z_2(c^*) \oplus \Z_{\bar \delta(l)}$, where the second summand is generated by the Witt class of the form
\[\left( \Z^2, 
\left[ \begin{array}{cc}
0 & 1 \\
-1 & 0
\end{array} \right], 
\left( \begin{array}{c} (0, 1) \\ (0, 1) \end{array} \right)\right).
\]

\item We have $W_0(Q^- \oplus \Z_l) = 0$.
\end{compactenum}
\end{example}

To consider the naturality of the isomorphisms in Theorem~\ref{t:main2},
we recall from Definition~\ref{d:split} that a morphism $\alpha \colon Q_1 \oplus G_1 \to Q_2 \oplus G_2$ of split form parameters 
is split if $\alpha = \alpha_Q \oplus \alpha_G$ for a morphism $\alpha_Q \colon Q_1 \to Q_2$ and a homomorphism $\alpha_G \colon G_1 \to G_2$.  A morphism
$\alpha \colon P_1 \to P_2$ is splittable 
if there are splittings of $P_1$ and $P_2$ with respect to which $\alpha$ is split. By the naturality of the isomorphism in Theorem \ref{thm:main} (and the discussion after Definition \ref{def:red-Witt}), we get the following: 

\begin{prop} \label{p:splittable}
%%%%%%%%%%%%%%%%
Let $\alpha \colon P_1 \to P_2$ be a splittable morphism of form parameters and let
$\alpha = \alpha_Q \oplus \alpha_G$ be a splitting of $\alpha$
with respect to splittings $P_1 \cong Q_1 \oplus G_1$ and $P_2 = Q_2 \oplus G_2$.
Then 
\[ 
W_0(\alpha) = W_0(\alpha_Q) \oplus (\alpha_G \otimes_\Z \alpha_Q) \colon W_0(P_1) \cong W_0(Q_1) \oplus (G_1 \tensor_\Z Q_1) \to W_0(P_2) \cong W_0(Q_2) \oplus (G_2 \tensor_\Z Q_2). 
\qed
\]
\end{prop}

\begin{remark*}
If the splittings of $P_1$ and $P_2$ in Proposition~\ref{p:splittable} are maximal,
then $\alpha_Q$ is the composition of some standard morphisms and isomorphisms by Theorem \ref{thm:mor}, so $W_0(\alpha_Q)$ can be determined using Theorem \ref{t:wf-ind}.
\end{remark*}

\begin{example}
%%%%%%%%%
The morphism of form parameters $Q_{-} \oplus \Z_2 \xra{\ql \oplus {\rm Id}} \Z^\Lambda_2 \oplus \Z_2$ induces the zero map 
\[ W_{0}(Q_- \oplus \Z_2) \cong (\Z_2(c^*) \oplus \Z_2) \xra{~0~} (0 \oplus \Z_2) \cong  W_0(\Z^\Lambda_2 \oplus \Z_2). \]
\end{example}

\section{A natural description of Witt groups} \label{s:WGN}
%%%%%%%%%%%%%%%%%%%%%%%%%%%%%%%%%%%%%%%%%%%%%%%%%%%%%%%%%%%%%%%%%

By the results of the previous section we can determine the isomorphism type of $W_0(P)$ for any form parameter $P$, by choosing a maximal splitting. Proposition \ref{p:splittable} can also be used to determine the homomorphism between the Witt groups induced by a splittable morphism between form parameters. However, as Example~\ref{ex:non-splittable} shows, not all morphisms are splittable (with respect to some non-trivial splittings of the form parameters). 

In this section we give a natural description of the group $W_0(P)$, for any form parameter $P$, which, in particular, allows us to determine the homomorphism induced by any morphism of form parameters. 

We begin by outlining the strategy. 
Recall that for any symmetric form parameter $P$ there is an extended symmetrisation morphism,
$\es \colon P \to Q^+ \oplus SP$. 
In fact, the morphisms $\es$ define a natural transformation from the identity to the extended symmetrisation functor $\ES \colon \FP_+ \rightarrow \FP_+$, where $\ES(P) := Q^+ \oplus SP$ on objects and $\ES(\alpha) := \Id_{Q^+} \oplus S\alpha$ on morphisms. By Theorem~\ref{t:main2}, there is a natural isomorphism $W_0(Q^+\oplus SP) \cong \Z \oplus \Gamma(SP)$,
and we will show that the induced homomorphism 
\[
W_0(\es) \colon W_0(P) \to W_0(Q^+\oplus SP) \cong \Z \oplus \Gamma(SP)
\] 
is injective. Then we will determine $\Im(W_0(\es))$ as a natural subgroup of $\Z \oplus \Gamma(SP)$. 
Then, for any morphism $\alpha \colon P \rightarrow P'$, we can compute the induced map $W_0(\alpha) \colon W_0(P) \to W_0(P')$ as the restriction of $\Id_{\Z} \oplus \Gamma(S\alpha) \colon \Z \oplus \Gamma(SP) \rightarrow \Z \oplus \Gamma(SP)$, 
which is a homomorphism $\Im(W_0(\es_P)) \rightarrow \Im(W_0(\es_{P'}))$.

Dually, for any anti-symmetric form parameter $P$ there is an extended quadratic lift morphism, 
$\eql \colon Q_- \oplus P_e \to P$, 
and the morphisms $\eql$ define a natural transformation from the extended quadratic lift functor 
$\EQL$ to the identity, where $\EQL \colon \FP_- \rightarrow \FP_-$ has $\EQL(P) := Q_- \oplus P_e$ on objects and $\EQL(\alpha) := \Id_{Q_-} \oplus \alpha$ on morphisms. 
By Theorem~\ref{t:main2}, there is a natural isomorphism $W_0(Q_- \oplus P_e) \cong \Z_2 \oplus \Lambda_1(P_e)$. We will show that the induced homomorphism 
\[
W_0(\eql) \colon W_0(Q_- \oplus P_e) \cong \Z_2 \oplus \Lambda_1(P_e) \to W_0(P)
\] 
is surjective and determine $\Ker(W_0(\eql))$ as a natural subgroup of $\Z_2 \oplus \Lambda_1(P_e)$. 
Then, for any morphism $\alpha \colon P \rightarrow P'$, the homomorphism 
$W_0(\alpha) \colon W_0(P) \to W_0(P')$ can be computed by identifying it with the homomorphism $\Z_2 \oplus \Lambda_1(P_e) / \Ker(W_0(\eql_P)) \rightarrow \Z_2 \oplus \Lambda_1(P'_e) / \Ker(W_0(\eql_{P'}))$, 
which is induced by the homomorphism $\Id_{\Z_2} \oplus \Lambda_1(\alpha) \colon \Z_2 \oplus \Lambda_1(P_e) \rightarrow \Z_2 \oplus \Lambda_1(P'_e)$.

\subsection{Symmetric form parameters} 
%%%%%%%%%%%%%%%%%%%%%%%%%%%%%%%%%%%%%%%%%%%%%%%%%%%%%%
First we consider the restriction of the Witt group functor to the category of symmetric form parameters, $W_0 \colon \FP_+ \to \AbF$.

\begin{lem} \label{lem:indec-es-inj}
%%%%%%%%%%%%%%%%%%
If $Q$ is an indecomposable symmetric form parameter, then the induced homomorphism 
$W_0(\es) \colon W_0(Q) \to W_0(Q^+ \oplus SQ)$ is injective.
\end{lem}

\begin{proof}
Since $Q^+$ is the terminal object, the composition $Q \rightarrow Q^+ \oplus SQ \rightarrow Q^+$ is the unique morphism, 
and the induced map $W_0(Q) \to W_0(Q^+ \oplus SQ) \cong W_0(Q^+) \oplus W_0^{Q^+}(SQ) \rightarrow W_0(Q^+)$ 
is given by Theorem \ref{t:wf-ind}. To compute the other component $W_0(Q) \rightarrow W_0^{Q^+}(SQ) \cong SQ \tensor_\Z Q^+$ 
we will use the split norm $F$ from Definition \ref{def:F-1}, which induces the isomorphism $\Tilde{F}$ (Theorem \ref{thm:main}). 

If $Q = Q_+$, then $SQ_+ = 0$, and the map $W_0(Q_+) \to W_0(Q^+ \oplus SQ_+) \cong W_0(Q^+)$ is equal to 
the inclusion $\Z(8\sigma^*) \rightarrow \Z(\sigma^*)$. 

If $Q = Q^+$, then $SQ^+ = \Z_2$. The composition $W_0(Q^+) \to W_0(Q^+ \oplus SQ^+) \rightarrow W_0(Q^+)$ is of course the identity, hence the map $W_0(Q^+) \to W_0(Q^+ \oplus SQ^+)$ is injective.

If $Q = \Z^P$, then $S\Z^P = \Z$ and $W_0^{Q^+}(S\Z^P) \cong \Z \tensor_\Z Q^+ \cong \Z$ by Theorem \ref{t:main2}. The generators of $W_0(\Z^P) \cong \Z(\sigma^*) \oplus \Z(\rho_{\infty}^*)$ are (the equivalence classes of) $\sigma^* = (\Z, 1, (1,0))$ and 
$\rho_{\infty}^* = \left( \Z^2, \left[\begin{smallmatrix} 0 & 1 \\ 1 & 0 \end{smallmatrix}\right], \left(\begin{smallmatrix} (0,1) \\ (0,-1) \end{smallmatrix}\right) \right)$ 
(it is easy to verify that for these forms we have $\sigma(\sigma^*)=\rho_{\infty}(\rho_{\infty}^*)=1$ and $\rho_{\infty}(\sigma^*)=\sigma(\rho_{\infty}^*)=0$). The extended symmetrisation map is given by the homomorphism $\left(\begin{smallmatrix} 1 & 0 \\ 1 & 2 \end{smallmatrix}\right) \colon \Z^P_e = \Z \oplus \Z \rightarrow (Q^+ \oplus S\Z^P)_e = \Z \oplus \Z$, so the generators are mapped to $(\Z, 1, (1,1))$ and $(\Z^2, \left[\begin{smallmatrix} 0 & 1 \\ 1 & 0 \end{smallmatrix}\right], \left(\begin{smallmatrix} (0,2) \\ (0,-2) \end{smallmatrix}\right))$ respectively. Their images under $F$ are $1 \otimes 1$ and $[2,-2] \otimes 1$, and under the isomorphism $\Z \tensor_\Z Q^+ \cong \Z$ these correspond to $1$ and $-8$. Therefore the extended symmetrisation map induces the homomorphism $\left(\begin{smallmatrix} 1 & 0 \\ 1 & -8 \end{smallmatrix}\right) \colon W_0(\Z^P) = \Z \oplus \Z \rightarrow W_0(Q^+ \oplus S\Z^P) = \Z \oplus \Z$, and this is injective.

If $Q = \Z^P_k$, then $S\Z^P_k = \Z_{2^{k+1}}$ and $W_0^{Q^+}\!\!(S\Z^P_k) \cong \Z_{2^{k+1}} \tensor_\Z Q^+ \cong \Z_{2^{k+2}}$ by Theorem \ref{t:main2}. The generators of $W_0(\Z^P_k) \cong \Z(\sigma^*) \oplus \Z_{2^{k-1}}(\rho_k^*)$ are $\sigma^* = (\Z, 1, (1,0))$ and $\rho_k^* = (\Z^2, \left[\begin{smallmatrix} 0 & 1 \\ 1 & 0 \end{smallmatrix}\right], \left(\begin{smallmatrix} (0,1) \\ (0,-1) \end{smallmatrix}\right))$. The extended symmetrisation map is given by 
$\left(\begin{smallmatrix} 1 & 0 \\ 1 & 2 \end{smallmatrix}\right) \colon (\Z^P_k)_e = \Z \oplus \Z_{2^k} \rightarrow (Q^+ \oplus S\Z^P_k)_e = \Z \oplus \Z_{2^{k+1}}$, so the generators are mapped to $(\Z, 1, (1,1))$ and $(\Z^2, \left[\begin{smallmatrix} 0 & 1 \\ 1 & 0 \end{smallmatrix}\right], \left(\begin{smallmatrix} (0,2) \\ (0,-2) \end{smallmatrix}\right))$ respectively. Their images under $F$ are $1 \otimes 1$ and $[2,-2] \otimes 1$, and under the isomorphism $\Z_{2^{k+1}} \tensor_\Z Q^+ \cong \Z_{2^{k+2}}$ these correspond to $1$ and $-8$. Therefore the extended symmetrisation map induces 
$\left(\begin{smallmatrix} 1 & 0 \\ 1 & -8 \end{smallmatrix}\right) \colon W_0(\Z^P_k) = \Z \oplus \Z_{2^{k-1}} \rightarrow W_0(Q^+ \oplus S\Z^P_k) = \Z \oplus \Z_{2^{k+2}}$, and this map is also injective.
\end{proof}

\begin{lemma} \label{l:Q_to_Q+SQ}
%%%%%%%%%%%%%%%%%%%%%%
If $P$ is a symmetric form parameter, then
$W_0(\es) \colon W_0(P) \to W_0(Q^+ \oplus SP)$ is injective.
\end{lemma}

\begin{proof}
By Theorem \ref{t:qfp-class} it is enough to prove the statement for form parameters of the form $P = Q \oplus G$, where $Q$ is an indecomposable form parameter and $G$ is an abelian group.

We have $SP \cong SQ \oplus G$, and the 
map 
$\es \colon P = Q \oplus G \rightarrow Q^+ \oplus SP \cong (Q^+ \oplus SQ) \oplus G$ splits as the sum of $\es \colon Q \rightarrow Q^+ \oplus SQ$ and $\Id_G$. Therefore the induced map 
\[ W_0(P) = W_0(Q) \oplus W_0^Q(G) \rightarrow W_0(Q^+ \oplus SP) \cong W_0(Q^+ \oplus SQ) \oplus W_0^{Q^+ \oplus SQ}(G) \]
splits as the sum of $W_0(Q) \rightarrow W_0(Q^+ \oplus SQ)$ and the map
\[
W_0^Q(G) \cong G \otimes_{\Z} Q \rightarrow W_0^{Q^+ \oplus SQ}(G) \cong G \otimes_{\Z} (Q^+ \oplus SQ) \cong (G \otimes_{\Z} Q^+) \oplus (G \otimes SQ) \cong \Gamma(G) \oplus (G \otimes SQ) .
\]
The former is injective by Lemma \ref{lem:indec-es-inj}. 
By Proposition \ref{p:qt-pback} b) there is a pullback diagram
\[
\xymatrix{
G \otimes_{\Z} Q \ar[d] \ar[r] & G \otimes SQ \ar[d] \\
\Gamma(G) \ar[r] & G \otimes \Z_2 .
}
\]
The map $W_0^Q(G) \cong G \otimes_{\Z} Q \rightarrow W_0^{Q^+ \oplus SQ}(G) \cong \Gamma(G) \oplus (G \otimes SQ)$ is the direct sum of the maps $G \otimes_{\Z} Q \rightarrow \Gamma(G)$ and $G \otimes_{\Z} Q \rightarrow G \otimes SQ$ in this diagram, hence it is injective too. 
\end{proof}

Next we determine the image of $W_0(\es) \colon W_0(P) \to W_0(Q^+ \oplus SP) \cong \Z \oplus \Gamma(SP)$. Given an abelian group $A$, we will use the identification $\Gamma(A) \cong A \otimes_{\Z} Q^+$ to describe elements of $\Gamma(A)$ (in particular, $\Gamma(A)$ is generated by the symbols $a \otimes 1$ and $[a,b] \otimes 1$ for $a,b \in A$). Note that there is a short exact sequence
\[
\xymatrix{
0 \ar[r] & S^2(A) \ar[r] & \Gamma(A) \ar[r] & A \otimes \Z_2 \ar[r] & 0
}
\]
where the first map sends $ab$ to $[a,b] \otimes 1$, and the second map sends $[a,b] \otimes 1$ to $0$ and $a \otimes 1$ to $a \otimes 1$ (see Proposition \ref{p:qt-pback} b) or \cite[Remark 2.10]{Bau1}).

\begin{defin} \label{d:Sigma(v)}
Let $v \colon A \to \Z_2$ be a homomorphism. We define the subgroup 
\[
\Sigma(v) := \left< (1,x \otimes 1),(0,[k_1, k_2] \otimes 1),(8,0) \bigm| x \in v^{-1}(1), k_1, k_2 \in \Ker(v) \right> \leq \Z \oplus \Gamma(A).
\]
\end{defin}
\noindent
We note that if $\alpha \colon v_1 \rightarrow v_2$ is a morphism in $\AbF / \Z_2$, then the induced map $\Id_{\Z} \oplus \Gamma(\alpha) \colon \Z \oplus \Gamma(A_1) \rightarrow \Z \oplus \Gamma(A_2)$ sends the generators of $\Sigma(v_1)$ to generators of $\Sigma(v_2)$. This shows that $\Sigma$ defines a functor $\AbF / \Z_2 \to \AbF$.

\begin{lem} \label{lem:sigma-comp}
Let $v \colon A \to \Z_2$ be a homomorphism. 

a) If $v=0$, then $\Sigma(v) = \left< (0,[k_1, k_2] \otimes 1),(8,0) \bigm| k_1, k_2 \in \Ker(v) = A \right>$. 

b) If $v \neq 0$, then $\Sigma(v) = \left< (1,x \otimes 1) \bigm| x \in v^{-1}(1) \right>$.
\end{lem}

\begin{proof}
Part a) holds, because $v^{-1}(1) = \emptyset$. 

b) Let $\Sigma'(v)$ denote the right hand side, we need to show that it contains the remaining generators of $\Sigma(v)$. First note that if $x \in v^{-1}(1)$ and $k \in \Ker(v)$, then $[k, x] \otimes 1 + k \otimes 1 = (x+k) \otimes 1 - x \otimes 1$, hence $(0,[k, x] \otimes 1 + k \otimes 1) \in \Sigma'(v)$. We have $[k_1, k_2] \otimes 1 = ([k_1, x+k_2] \otimes 1 + k_1 \otimes 1) - ([k_1, x] \otimes 1 + k_1 \otimes 1)$ for any $k_1, k_2 \in \Ker(v)$ and $x \in v^{-1}(1)$, therefore $(0,[k_1, k_2] \otimes 1) \in \Sigma'(v)$. Finally, for any $y \in A$, we have $[2y,2y] \otimes 1 = 2y \otimes 2 = y \otimes 2 + y \otimes 2 + [y,y] \otimes 2 = y \otimes 8$. So, with any $x \in v^{-1}(1)$, we have $(8,0) = 8(1,x \otimes 1) - (0,x \otimes 8) \in \Sigma'(v)$.
\end{proof}

\begin{thm} \label{thm:im-sigma}
Let $P$ be a symmetric form parameter. Then under the identification $W_0(Q^+ \oplus SP) \cong \Z \oplus \Gamma(SP)$ we have $\Im(W_0(\es)) = \Sigma(v_P)$.
\end{thm}

\begin{proof}[Proof of Theorem~\ref{thm:im-sigma}]
First assume that $v_P = 0$, i.e.\ $P \cong Q_+ \oplus SP$. It follows from Lemmas \ref{lem:indec-es-inj} and \ref{l:Q_to_Q+SQ} that the map $W_0(\es) \colon W_0(P) \cong W_0(Q_+ \oplus SP) \cong W_0(Q_+) \oplus W_0^{Q_+}(SP) \rightarrow W_0(Q^+ \oplus SP) \cong W_0(Q^+) \oplus W_0^{Q^+}(SP)$ splits as the sum of the inclusions $W_0(Q_+) \cong \Z(8\sigma^*) \to W_0(Q^+) \cong \Z(\sigma^*)$ and $W_0^{Q_+}(SP) \cong S^2(SP) \rightarrow W_0^{Q^+}(SP) \cong \Gamma(SP)$. Therefore the image of $W_0(\es)$ is generated by $(8,0)$ and the elements $(0,[k_1, k_2] \otimes 1)$ for $k_1, k_2 \in SP = \Ker(v_P)$, so it is equal to $\Sigma(v_P)$ by Lemma \ref{lem:sigma-comp} a).

Now assume that $v_P \neq 0$. By adding the metabolic $P$-form
\[
\left(
\Z^2, 
\left[ \begin{array}{cc} 0 & 1 \\ 1 & {\rm h}(q) \end{array} \right],
\left( \begin{array}{c} 0 \\ q \end{array} \right) 
\right)
\]
(where $q \in v_P^{-1}(1)$) we get that every element of $W_0(P)$ has a representative $(X, \lambda, \mu)$ such that $\lambda$ is indefinite and odd. By the classification of nonsingular indefinite forms (see \cite[Ch.\ II.\ Theorem (5.3)]{M-H}), $(X, \lambda, \mu)$ is isomorphic to the sum of $P$-forms of the form $(\Z, [\pm 1], q)$ (where ${\rm h}(q) = \pm 1$). Therefore elements of this form generate $W_0(P)$. The image of such a form under extended symmetrisation is $(\Z, [\pm 1], (\pm 1, \pi_P(q)))$, which corresponds to $(\pm 1, \pi_P(q) \otimes (\pm 1)) = \pm (1, \pi_P(q) \otimes 1)$ under the isomorphism $W_0(Q_+ \oplus SP) \cong \Z \oplus \Gamma(SP)$ (see Definition \ref{def:F-1}). Since $\{ \pi_P(q) \mid {\rm h}(q) = \pm 1 \} = v_P^{-1}(1)$, this implies that $\Im(W_0(\es)) = \left< (1,x \otimes 1) \mid x \in v_P^{-1}(1) \right>$. This is equal to $\Sigma(v_P)$ by Lemma \ref{lem:sigma-comp} b).
\end{proof}

Using Lemma \ref{l:Q_to_Q+SQ}, Theorem \ref{thm:im-sigma} and the fact that 
extended symmetrisation 
is a natural transformation, we get

\begin{corollary} \label{cor:nat+}
The homomorphisms induced by extended symmetrisation define a natural isomorphism $W_0 \cong \Sigma \circ v_{(-)} : \FP_+ \to \AbF$. \qed
\end{corollary}

To give a more explicit description of the group $\Sigma(v)$, we include Theorem \ref{t:sigma-comp} below, which requires the following

\begin{defin} \label{d:Cv}
Let $v \colon A \to \Z_2$ be a homomorphism. We define the subgroups 
\begin{align*}
\Phi(v) & := \left< x \otimes 1, [k, x] \otimes 1 + k \otimes 1, [k_1, k_2] \otimes 1 \bigm| x \in v^{-1}(1), k, k_1, k_2 \in \Ker(v) \right>  \leq \Gamma(A), \\
\Psi(v) & := \left< [k, x] \otimes 1 + k \otimes 1, [k_1, k_2] \otimes 1 \bigm| x \in v^{-1}(1), k, k_1, k_2 \in \Ker(v) \right>  \leq \Phi(v).
\end{align*}
Note that, by the formula in the proof of Lemma \ref{lem:sigma-comp}, the generators $[k, x] \otimes 1 + k \otimes 1$ of $\Phi(v)$ could be omitted, and that if $v = 0$, then $\Phi(v) = \Psi(v)$. We also define the groups 
\[
\begin{aligned}
C(v) & := 
\begin{cases}
\Z_8 & \text{if } v = 0, \\
\Phi(v)/\Psi(v) & \text{if } v \neq 0,
\end{cases} 
\\
\Upsilon(v) & := 
\begin{cases}
\Z_8 \oplus \Gamma(A)/\Psi(v) & \text{if } v = 0, \\
\Gamma(A)/\Psi(v) & \text{if } v \neq 0.
\end{cases}
\end{aligned}
\]
\end{defin}

\begin{thm} \label{t:sigma-comp}
%%%%%%%%%%%%%%%%
Let $v \colon A \to \Z_2$ be a homomorphism. 
\begin{compactenum}[a)]
\item We have 
\[
C(v) \cong 
\begin{cases}
\Z_4 & \text{if $v \cong 1_1 \oplus G$ for some $G$,} \\
\Z_8 & \text{otherwise}.
\end{cases} 
\]
Moreover, if $v \neq 0$, then the elements of $C(v)$ of the form $[x \otimes 1]$ with $x \in v^{-1}(1)$ are all equal to each other, and this element generates $C(v)$.
\item There is a commutative diagram with short exact rows and columns
\[
\xymatrix{
 & & 0 \ar[d] & 0 \ar[d] & \\
0 \ar[r] & \Sigma(v) \ar@{=}[d] \ar[r] & \Z \oplus \Phi(v) \ar[d] \ar[r]^-{\iota - q} & C(v) \ar[d] \ar[r] & 0 \\
0 \ar[r] & \Sigma(v) \ar[r] & \Z \oplus \Gamma(A) \ar[d]^-{0+u_v} \ar[r]^-{\iota - \bar{q}} & \Upsilon(v) \ar[d]^-{\tilde{u}_v} \ar[r] & 0 \\
 & & \Ker(v_2) \ar@{=}[r] \ar[d] & \Ker(v_2) \ar[d] & \\
 & & 0 & 0, & 
}
\]
where the following hold:
\begin{compactitem}[$\bullet$]
\item $v_2 = v \otimes \Id_{\Z_2} \colon A \otimes \Z_2 \rightarrow \Z_2 \otimes \Z_2 = \Z_2$, hence $\Ker(v_2) \leq A \otimes \Z_2$;
\item The homomorphism $\iota \colon \Z \rightarrow C(v) \leq \Upsilon(v)$ is defined by
\[
\iota(1) = 
\begin{cases}
1 & \text{if } v = 0,\\
[x \otimes 1] \text{~for some $x \in v^{-1}(1)$} & \text{if } v \neq 0;
\end{cases}
\]
\item The homomorphism $q \colon \Phi(v) \rightarrow C(v)$ is the quotient map $\Phi(v) \rightarrow \Phi(v)/\Psi(v)$ if $v \neq 0$, and trivial if $v = 0$;
\item The homomorphism $\bar{q} \colon \Gamma(A) \rightarrow \Upsilon(v)$ is the quotient map $\Gamma(A) \rightarrow \Gamma(A)/\Psi(v)$;
\item The homomorphism $u_v \colon \Gamma(A) \rightarrow A \otimes \Z_2$ is given by $u_v([a_1,a_2] \otimes 1) = a_1 \otimes v(a_2) + a_2 \otimes v(a_1)$ and $u_v(a \otimes 1) = a \otimes 1 + a \otimes v(a)$;
\item The homomorphism $\tilde{u}_v \colon \Upsilon(v) \rightarrow A \otimes \Z_2$ is the map $\Gamma(A)/\Psi(v) \rightarrow A \otimes \Z_2$ induced by $u_v$ (precomposed with the projection $\Z_8 \oplus \Gamma(A)/\Psi(v) \rightarrow \Gamma(A)/\Psi(v)$ if $v = 0$);
\item The unlabelled maps are the obvious inclusions. 
\end{compactitem}
\end{compactenum}
\end{thm}

\begin{remark}
Theorem~\ref{t:sigma-comp} shows that $\Sigma(v)$ can be described as the kernel of a map $\Z \oplus \Gamma(A) \rightarrow \Upsilon(v)$, where $\Upsilon(v)$ is an extension of the $\Z_2$-vector space $\Ker(v_2)$ by the cyclic group $C(v)$. Alternatively, it is the subgroup of elements $(n, \alpha) \in \Z \oplus \Phi(v)$ such that $n \equiv q(\alpha)$ mod $4$ or $8$ (depending on $C(v)$), where $\Phi(v)$ is the kernel of the surjective map $u_v \colon \Gamma(A) \rightarrow \Ker(v_2)$. In particular, we see that $\Sigma(v) \cap (0 \oplus \Gamma(A)) = \Ker(\bar{q}) = \Psi(v)$, and that the image of the composition $\Sigma(v) \rightarrow \Z \oplus \Gamma(A) \rightarrow \Gamma(A)$ is $\Phi(v)$. 
\end{remark}

\begin{proof}[Proof of Theorem~\ref{t:sigma-comp}]
a) If $v = 0$, then the statement holds by definition. In the rest of this proof we restrict to the $v \neq 0$ case. Then the elements $[x \otimes 1]$ with $v(x)=1$ generate $\Phi(v) / \Psi(v)$. If $v(x)=v(y)=1$, then $x \otimes 1 - y \otimes 1 = (x-y) \otimes 1 + [y,x-y] \otimes 1 \in \Psi(v)$, so all of these generators coincide. Hence $\Phi(v) / \Psi(v)$ is cyclic. 

By Theorem \ref{t:(co)slice-d} $v$ is isomorphic to $v_0 \oplus G$ for some abelian group $G$, where $v_0 \colon A_0 \rightarrow \Z_2$ is either $1_k \colon \Z_{2^k} \rightarrow \Z_2$ for some positive integer $k$ or $1_{\infty} \colon \Z \rightarrow \Z_2$. Since $\Phi(v)$ and $\Psi(v)$ depend functorially on $v$, the morphism $v_0 \rightarrow v_0 \oplus G \cong v$ induces a homomorphism $C(v_0) \rightarrow C(v)$. This map is injective, because the composition $v_0 \rightarrow v_0 \oplus G \rightarrow v_0$ is the identity, and surjective, because it maps a generator of the cyclic group $C(v_0)$ (of the form $[x \otimes 1]$) into a generator of $C(v)$. Therefore $C(v) \cong C(v_0)$. 

We have $\Gamma(A_0) \cong \Z_{2^{k+1}}$ or $\Gamma(A_0) \cong \Z$, generated by $1 \otimes 1$. So from the definitions we see that in all cases $\Phi(v_0) = \Gamma(A_0)$ and $\Psi(v_0) = 8\Gamma(A_0)$. Moreover, $\Gamma(A_0) / 8\Gamma(A_0) \cong \Z_4$ if $k=1$, and $\Gamma(A_0) / 8\Gamma(A_0) \cong \Z_8$ otherwise. 

b) First consider the top horizontal sequence. It follows from the definitions that $\Sigma(v) \leq \Z \oplus \Phi(v)$. The map $\iota$ is well-defined, because if $v(x)=v(y)=1$, then $x \otimes 1 - y \otimes 1 = (x-y) \otimes 1 + [x-y,y] \otimes 1 \in \Psi(v)$. Moreover, $\iota$ is surjective, hence $\iota - q$ is surjective too. We need to prove that the sequence is exact at $\Z \oplus \Phi(v)$. If $v = 0$, then $\Sigma(v) = 8\Z \oplus \Phi(v) = \Ker(\iota - q)$ by Lemma \ref{lem:sigma-comp} a), so the sequence is exact. If $v \neq 0$, then $\Sigma(v)$ is generated by the elements $(1,x \otimes 1)$ with $x \in v^{-1}(1)$, so $\Sigma(v) \leq \Ker(\iota - q)$. Every element of $\Z \oplus \Phi(v)$ can be written as $n(1,x \otimes 1)+(0,\alpha)$ for some $n \in \Z$ and $\alpha \in \Phi(v)$, and such an element is in $\Ker(\iota - q)$ if and only if $\alpha \in \Ker(q) = \Psi(v)$. In the proof of Lemma \ref{lem:sigma-comp} b) we saw that $0 \oplus \Psi(v) \leq \Sigma(v)$, and this implies that $\Sigma(v) \geq \Ker(\iota - q)$. Therefore the sequence is exact. 

Next consider the middle horizontal sequence. The map $\iota - \bar{q}$ is surjective, because $\bar{q}$ is surjective onto $\Gamma(A)/\Psi(v)$, and if $v=0$, then $\iota$ is surjective onto $\Z_8$. If $v = 0$, then $\Sigma(v) = 8\Z \oplus \Phi(v) = 8\Z \oplus \Psi(v) = \Ker(\iota - \bar{q})$. If $v \neq 0$, then, using the observations from the top sequence, $\Sigma(v)$ and $\Ker(\iota - \bar{q})$ are both equal to the subgroup generated by the elements $(1,x \otimes 1)$ (where $x \in v^{-1}(1)$) and $0 \oplus \Psi(v)$. Therefore the middle sequence is exact too.

Next we consider the vertical sequence in the middle. The homomorphism $u_v$ is well-defined, because $u_v([a,a] \otimes 1) = 0 = u_v(a \otimes 2)$ and $u_v((a_1+a_2) \otimes 1) = (a_1+a_2) \otimes 1 + (a_1+a_2) \otimes v(a_1+a_2) = u_v(a_1 \otimes 1) + u_v(a_2 \otimes 1) + u_v([a_1,a_2] \otimes 1)$ for every $a,a_1,a_2 \in A$. Moreover, $\Im(u_v) \leq \Ker(v_2)$, because $v_2(a_1 \otimes v(a_2) + a_2 \otimes v(a_1)) = v(a_1) \otimes v(a_2) + v(a_2) \otimes v(a_1) = 0 \in \Z_2$ and $v_2(a \otimes 1 + a \otimes v(a)) = v(a) \otimes (1+v(a)) =0$.

For the exactness of the sequence we need to show that $\Phi(v) \leq \Ker(u_v)$ and that the induced homomorphism $\Gamma(A) / \Phi(v) \rightarrow \Ker(v_2)$ is an isomorphism. The former is verified by evaluating $u_v$ on the generators of $\Phi(v)$. For the latter, we show that the induced homomorphism has an inverse $\Ker(v_2) \rightarrow \Gamma(A) / \Phi(v)$. Namely, we have $\Ker(v_2) = \{ k \otimes 1 \mid k \in \Ker(v) \}$ and the inverse map is given by $k \otimes 1 \mapsto [k \otimes 1]$. This map is well-defined, because $k_1 \otimes 1 = k_2 \otimes 1 \in A \otimes \Z_2$ if and only if $k_2 = k_1 + 2a$ for some $a \in A$, and $(k_1 + 2a) \otimes 1 - k_1 \otimes 1 = 2a \otimes 1 + [k_1,2a] \otimes 1 \in \Phi(v)$. Here we used that $2a \in \Ker(v)$, so $[k_1,2a] \otimes 1 \in \Phi(v)$, and $2a \otimes 1 = 4(a \otimes 1) = 2([a,a] \otimes 1) \in \Phi(v)$ for both $a \in v^{-1}(1)$ and $a \in \Ker(v)$. This map is a homomorphism, because $[(k_1+k_2) \otimes 1] - [k_1 \otimes 1] - [k_2 \otimes 1] = [[k_1,k_2] \otimes 1] = 0 \in \Gamma(A) / \Phi(v)$. It is the inverse of the homomorphism induced by $u_v$, because $u_v(k \otimes 1) = k \otimes 1$ for $k \in \Ker(v)$, and the elements $[k \otimes 1]$ generate $\Gamma(A) / \Phi(v)$. Therefore $\Gamma(A) / \Phi(v) \cong \Ker(v_2)$ and the sequence is exact.
 
Finally consider the vertical sequence on the right. Since $u_v$ induces an isomorphism $\Gamma(A) / \Phi(v) \cong \Ker(v_2)$ and $\Psi(v) \leq \Phi(v)$, the map $\tilde{u}_v$ is well-defined. If $v=0$, then $\Psi(v) = \Phi(v)$, hence the restriction of $\tilde{u}_v$ to $\Gamma(A)/\Psi(v) \leq \Upsilon(v)$ is an isomorphism, so the sequence is exact. It is also exact if $v \neq 0$, because then it is isomorphic to the sequence $\Phi(v)/\Psi(v) \rightarrow \Gamma(A)/\Psi(v) \rightarrow \Gamma(A)/\Phi(v)$.
\end{proof}

\begin{remark}
We can give an alternative description of the map $u_v$ in terms of $(Q^+ \oplus A)$-forms, using that $\Gamma(A) \cong W_0^{Q^+}(A)$ by Theorem \ref{thm:main}. Let $(X, \lambda, (\mu_0,\mu_A))$ be a nonsingular $(Q^+ \oplus A)$-form and $\chi, \bar{\chi} \in A$ be characteristic elements for $\lambda$ and $v \circ \mu_A$ respectively; i.e.\ for all $x \in X$, $\lambda(x, x) \equiv \lambda(\chi, x)$ mod $2$ and $v(\mu_A(x)) = \varrho_2(\lambda(\bar{\chi}, x))$. Since $\lambda$ is nonsingular, $\chi$ and $\bar{\chi}$ are well-defined modulo $2X$. We have
\[ 
v(\mu_A(\chi - \bar{\chi})) = \varrho_2(\lambda(\bar{\chi}, \chi - \bar{\chi})) = \varrho_2(\lambda(\bar{\chi},\chi) - \lambda(\bar{\chi}, \bar{\chi})) = 0
\]
so $\mu_A(\chi - \bar{\chi}) \otimes 1 \in \Ker(v_2)$. By comparing the definitions we can see that the image of $[X, \lambda, (\mu_0,\mu_A)] \in W_0^{Q^+}(A) \cong \Gamma(A)$ is $u_v([X, \lambda, (\mu_0,\mu_A)]) = \mu_A(\chi - \bar{\chi}) \otimes 1$. In particular, if $v \circ \mu_A$ is characteristic for $\lambda$ (see Definition \ref{def:char}), then $\chi \equiv \bar{\chi}$ mod $2X$, showing that, if $v = v_P \colon A = SP \rightarrow \Z_2$ for some form parameter $P$, then the image of the composition $W_0(P) \xra{W_0(\es)} W_0(Q^+ \oplus A) \rightarrow W_0^{Q^+}(A)$ is in $\Ker(u_v)$ (see Lemma \ref{l:gmtrc}). 
\end{remark}

\subsection{Anti-symmetric form parameters}
%%%%%%%%%%%%%%%%%%%%%%%%%%%%%%%%%%%%%%%%%%%%%%%%%%
Next we consider the restriction of the Witt group functor to the category of anti-symmetric form parameters, $W_0 \colon \FP_- \to \AbF$.

\begin{lem} \label{lem:indec-antis-surj}
%%%%%%%%%%%%%%%%%%%%
If $Q$ is an indecomposable anti-symmetric form parameter, then the induced homomorphism
$W_0(\eql) \colon W_0(Q_- \oplus Q_e) \to W_0(Q)$ is surjective.
\end{lem}

\begin{proof}
If $Q=Q_-$, then the composition $Q_- \rightarrow Q_- \oplus \Z_2 \rightarrow Q_-$ is the identity, and if $Q \neq Q_-$, then $W_0(Q) \cong 0$ by Theorem \ref{t:wf-ind}.
\end{proof}

\begin{lemma} \label{l:gen-antis-surj}
%%%%%%%%%%%%%%%%%%%%%%%
If $P$ is an anti-symmetric form parameter,
then $W_0(\eql) \colon W_0(Q_- \oplus P_e) \to W_0(P)$ is surjective.
\end{lemma}

\begin{proof}
Again we prove this using a maximal splitting $P = Q \oplus G$.
We have $P_e = Q_e \oplus G$, and the induced homomorphism 
\[ W_0(Q_- \oplus Q_e \oplus G) \cong W_0(Q_- \oplus Q_e) \oplus W_0^{Q_- \oplus Q_e}(G) \rightarrow W_0(P) = W_0(Q) \oplus W_0^Q(G)\]
splits as the sum of $W_0(Q_- \oplus Q_e) \rightarrow W_0(Q)$ and 
\[
W_0^{Q_- \oplus Q_e}(G) \cong G \otimes_{\Z} (Q_- \oplus Q_e) \cong (G \otimes_{\Z} Q_-) \oplus (G \otimes Q_e) \cong \Lambda_1(G) \oplus (G \otimes Q_e) \rightarrow W_0^Q(G) \cong G \otimes_{\Z} Q .
\]
The former is surjective by Lemma \ref{lem:indec-antis-surj}. By Corollary \ref{cor:qt-pushout} b) there is a pushout diagram
\[
\xymatrix{
G \otimes \Z_2 \ar[d] \ar[r] & \Lambda_1(G) \ar[d] \\
G \otimes Q_e \ar[r] & G \otimes_{\Z} Q. 
}
\]
The map $W_0^{Q_- \oplus Q_e}(G) \cong \Lambda_1(G) \oplus (G \otimes Q_e) \rightarrow W_0^Q(G) \cong G \otimes_{\Z} Q$ is the sum of the maps $\Lambda_1(G) \rightarrow G \otimes_{\Z} Q$ and $G \otimes Q_e \rightarrow G \otimes_{\Z} Q$ in this diagram, hence it is surjective too. 
\end{proof}

\begin{remark}
Lemma~\ref{l:gen-antis-surj} has the following alternative proof:
Since every nonsingular anti-symmetric bilinear form is isomorphic to a sum of hyperbolic forms, it is easy to see that the induced homomorphism $\HomZns(X, \eql) \colon \HomZns(X, Q_- \oplus P_e) \to \HomZns(X, P)$ is surjective for any torsion free abelian group $X$, and hence $W_0(\eql)$ induces a surjection on Witt groups.
\end{remark}

Let $A$ be an abelian group. Recall that 
$\Lambda_1(A) = (A \otimes A) / \An{x \otimes y + y \otimes x \mid x, y \in A}$.  We will denote the image of $x \otimes y \in A \otimes A$ in the quotient $\Lambda_1(A)$ by $x \ohat y$. Recall from Example \ref{ex:qtp-funct} that $A \otimes_{\Z} Q_- \cong \Lambda_1(A)$, under this identification $a \otimes 1$ and $[a,b] \otimes 1$ correspond to $a \ohat a$ and $a \ohat b$ respectively. By Proposition \ref{prop:qfunctors-ses-antisymm} there is a short exact sequence
\[
\xymatrix{
0 \ar[r] & A \otimes \Z_2 \ar[r] & \Lambda_1(A) \ar[r] & \Lambda(A) \ar[r] & 0
} 
\]
where the first map sends $a \otimes 1$ to $a \ohat a$, and the second map sends $a \ohat b$ to $a \wedge b$. By decomposing $A$ as a direct sum of cyclic groups, we can see that this short exact sequence splits, but the splitting is not natural.

\begin{defin} \label{d:Lambda(v')}
%%%%%%%%%%%%%%%%%%
Let $v' \colon \Z_2 \to A$ be a homomorphism. We define the subgroup 
\[ 
K(v') = \An{ (1, v'(1) \ohat v'(1)), (0, x \ohat x + x \ohat v'(1)) \mid x \in A } \leq \Z_2 \oplus \Lambda_1(A) 
\]
and the quotient group
\[ 
\Lambda(v') = (\Z_2 \oplus \Lambda_1(A)) / K(v') .
\]
\end{defin}

If $\alpha \colon v'_1 \rightarrow v'_2$ is a morphism in $\Z_2 / \AbF$, then the induced map $\Id_{\Z_2} \oplus \Lambda_1(\alpha) \colon \Z_2 \oplus \Lambda_1(A_1) \rightarrow \Z_2 \oplus \Lambda_1(A_2)$ sends the generators of $K(v'_1)$ to generators of $K(v'_2)$. This shows that $K$, and hence $\Lambda$, are functors $\Z_2 / \AbF \to \AbF$. Also note that if $v' = 0 \colon \Z_2 \rightarrow A$ is the trivial homomorphism, then $\Lambda(v') \cong \Lambda(A)$.

\begin{lem} \label{lem:ker1}
Let $v' \colon \Z_2 \to A$ be a homomorphism and suppose that $v' \cong v'_0 \oplus C \colon \Z_2 \rightarrow B \oplus C$ is a maximal splitting of $v'$. Then
\[ 
K(v') \cong 
\begin{cases}
\Z_2 \oplus (C \otimes \Z_2) & \text{if $v'_0 \cong \iota_1$, } \\
\Z_2 \oplus (A \otimes \Z_2) & \text{otherwise. }
\end{cases}
\]
\end{lem}

\begin{proof}
If $v' = 0$, then $K(v')$ is  generated by $(1,0)$ and the elements $(0, x \ohat x)$ for every $x \in A$, therefore it is the direct sum of the subgroups $\Z_2 \oplus 0$ and $0 \oplus A \otimes \Z_2$. 

If $v'_0 \cong \iota_1$, then $\Lambda_1(A) \cong \Lambda_1(\Z_2 \oplus C) \cong \Z_2 \oplus \Z_2 \otimes C \oplus \Lambda_1(C)$, and $v'(1) \ohat v'(1)$ is the generator of the first component $\Lambda_1(\Z_2) \cong \Z_2$. So the subgroup $K(v') \leq \Z_2 \oplus (\Z_2 \oplus \Z_2 \otimes C \oplus \Lambda_1(C))$ is generated by $(1,(1,0,0))$ and the elements $(0,(0,1 \otimes x, x \ohat x))$ for every $x \in C$. Since the map $\Z_2 \otimes C \rightarrow \Z_2 \oplus (\Z_2 \oplus \Z_2 \otimes C \oplus \Lambda_1(C))$, $1 \otimes x \mapsto (0,(0,1 \otimes x, x \ohat x))$ is an isomorphism onto its image, we have $K(v') \cong \Z_2 \oplus \Z_2 \otimes C$. 

If $v'_0 \cong \iota_k$, $k \geq 2$, then $\Lambda_1(A) \cong \Lambda_1(\Z_{2^k} \oplus C) \cong \Z_2 \oplus \Z_{2^k} \otimes C \oplus \Lambda_1(C)$. We have $v'(1) = 2^{k-1} \in \Z_{2^k}$, so $v'(1) \ohat v'(1) = 0 \in \Lambda_1(\Z_{2^k}) \cong \Z_2$. For the generator $1 \in \Z_{2^k}$ we have that $1 \ohat 1 + 1 \ohat v'(1) = 1 \ohat 1$ is the generator of $\Lambda_1(\Z_{2^k}) \cong \Z_2$. Therefore $K(v')$ is generated by $(1,(0,0,0))$, $(0,(1,0,0))$ and the elements $(0,(0, 2^{k-1} \otimes x, x \ohat x))$ for every $x \in C$. Moreover, the map $\Z_2 \otimes A \rightarrow \Z_2 \oplus \Lambda_1(A)$, $1 \otimes x \mapsto (0, x \ohat x + x \ohat v'(1))$ is an isomorphism onto the subgroup generated by $(0,(1,0,0))$ and the elements $(0,(0, 2^{k-1} \otimes x, x \ohat x))$ for every $x \in C$, showing that 
$K(v') \cong \Z_2 \oplus (\Z_2 \otimes A)$. 
\end{proof}

\begin{lem} \label{lem:ker2}
Let $P$ be an anti-symmetric form parameter and suppose that $P \cong Q \oplus G$ is a maximal splitting of $P$. Then
\[ 
\Ker\bigl( W_0(\eql) \colon W_0(Q_- \oplus P_e) \to W_0(P) \bigr) \cong 
\begin{cases}
\Z_2 \oplus (G \otimes \Z_2) & \text{if $Q \cong Q_-$, } \\
\Z_2 \oplus (P_e \otimes \Z_2) & \text{otherwise. }
\end{cases}
\]
\end{lem}

\begin{proof}
Recall from Lemma \ref{l:gen-antis-surj} that $W_0(\eql) \colon W_0(Q_- \oplus P_e) \to W_0(P)$ is the sum of $W_0(Q_- \oplus Q_e) \rightarrow W_0(Q)$ and $W_0^{Q_- \oplus Q_e}(G) \rightarrow W_0^Q(G)$, where the second map is the sum of the two maps in the pushout square of Corollary \ref{cor:qt-pushout} b). Since this is also a pullback square (see Remark \ref{rem:p+p}), we get that the kernel of the second map is isomorphic to $G \otimes \Z_2$. For the first map we have $W_0(Q_- \oplus Q_e) \cong W_0(Q_-) \oplus W_0^{Q_-}(Q_e) \cong \Z_2 \oplus \Lambda_1(Q_e)$ and $W_0(Q) \cong \Z_2$ if $Q \cong Q_-$ and $W_0(Q) \cong 0$ otherwise.

If $Q = Q^-$, then the first map is $\Z_2 \rightarrow 0$, so $\Ker(W_0(\eql)) \cong \Z_2 \oplus G \otimes \Z_2 = \Z_2 \oplus P_e \otimes \Z_2$.

If $Q = Q_-$, then the first map is $\Z_2 \oplus \Z_2 \rightarrow \Z_2$, which is the identity on the first component, which implies that $\Ker(W_0(\eql)) \cong \Z_2 \oplus G \otimes \Z_2$. We note that the second component is generated by $1 \ohat 1 \in \Lambda_1(Q_e)$, and by Definition \ref{d:gamma'} the corresponding element of $W_0^{Q_-}(Q_e)$ is represented by the $(Q_- \oplus Q_e)$-form
\[
\left(
\Z^2, 
\left[ \begin{array}{cc} 0 & 1 \\ -1 & 0 \end{array} \right],
\left( \begin{array}{c} (0,1) \\ (0,1) \end{array} \right) 
\right)
\]
Its image under $W_0(\eql)$ has Arf-invariant $1$, so it is the generator of $W_0(Q)$. Therefore the map $\Z_2 \oplus \Z_2 \rightarrow \Z_2$ restricts to the identity on both components.

Finally, if $Q = \Z^{\Lambda}_k$, $k \geq 2$, then the first map is $\Z_2 \oplus \Z_2 \rightarrow 0$, so $\Ker(W_0(\eql)) \cong \Z_2 \oplus \Z_2 \oplus (G \otimes \Z_2) \cong \Z_2 \oplus (P_e \otimes \Z_2)$.
\end{proof}

\begin{thm} \label{t:W0-nat-anti}
Let $P$ be an anti-symmetric form parameter. 
Then under the identification $W_0(Q_- \oplus P_e) \cong \Z_2 \oplus \Lambda_1(P_e)$ we have $\Ker(W_0(\eql)) = K(v'_P)$. 
\end{thm}

\begin{proof}
First we consider $(1, v'_P(1) \ohat v'_P(1)) \in \Z_2 \oplus \Lambda_1(P_e)$. If we take a maximal splitting $P \cong Q \oplus G$ of $P$ as in Lemma \ref{lem:ker2}, then this element is in the first component of $W_0(Q_- \oplus P_e) \cong W_0(Q_- \oplus Q_e) \oplus W_0^{Q_- \oplus Q_e}(G)$, so it is mapped to $W_0(Q)$ by $W_0(\eql)$. If $Q \not\cong Q_-$, then $W_0(Q) \cong 0$.  If $Q \cong Q_-$, then $(1, v'_P(1) \ohat v'_P(1))$ corresponds to $(1,1) \in \Z_2 \oplus \Z_2 \cong W_0(Q_-) \oplus W_0^{Q_-}(Q_e)$, which generates the kernel of the map to $\Z_2 \cong W_0(Q)$ by the observations in the proof of Lemma \ref{lem:ker2}. Therefore $(1, v'_P(1) \ohat v'_P(1)) \in \Ker(W_0(\eql))$ in both cases.

Next consider $x \ohat x + x \ohat v'_P(1) = x \ohat (x+v'_P(1)) \in \Lambda_1(P_e)$ for some $x \in P_e$. By Definition \ref{d:gamma'} the corresponding element of $W_0^{Q_-}(P_e)$ is represented by the $(Q_- \oplus P_e)$-form
\[
\left(
\Z^2 \left< e,f \right>, 
\left[ \begin{array}{cc} 0 & 1 \\ -1 & 0 \end{array} \right],
\left( \begin{array}{c} (0,x) \\ (0,x+v'_P(1)) \end{array} \right) 
\right).
\]
This is mapped by $W_0(\eql)$ to the element of $W_0(P)$ represented by the $P$-form
\[
\left(
\Z^2 \left< e,f \right>, 
\left[ \begin{array}{cc} 0 & 1 \\ -1 & 0 \end{array} \right],
\left( \begin{array}{c} x \\ x+v'_P(1) \end{array} \right) 
\right)
\]
which is metabolic, because $\left< e-f \right>$ is a lagrangian. This shows that $W_0(\eql)$ sends every generator of $K(v'_P)$ to $0$, hence $K(v'_P) \leq \Ker(W_0(\eql))$. Lemmas \ref{lem:ker1} and \ref{lem:ker2} show that $K(v'_P)$ and $\Ker(W_0(\eql_P))$ are isomorphic finite groups, therefore $K(v'_P) = \Ker(W_0(\eql))$.
\end{proof}

Using Lemma \ref{l:gen-antis-surj}, Theorem \ref{t:W0-nat-anti} and the fact that $\eql$ is a natural transformation, we get the following: 

\begin{corollary} \label{cor:nat-}
The homomorphisms $W_0(\eql)$ induce a natural isomorphism $\Lambda \circ v'_{(-)} \cong W_0 : \FP_- \to \AbF$.
\end{corollary}

We end this section with Theorem \ref{t:lambda-comp} below, which is the dual of Theorem \ref{t:sigma-comp}, and gives an explicit description of $\Lambda(v')$.

\begin{defin}
%%%%%%%%%%%%%%%%%%
Let $v' \colon \Z_2 \to A$ be a homomorphism. We define the subgroup 
\[ 
L(v') = \An{ x \ohat x + x \ohat v'(1) \mid x \in A } \leq \Lambda_1(A) 
\]
and the quotient group
\[ 
\Xi(v') = \Lambda_1(A) / L(v') .
\]
\end{defin}

Note that since $0 \oplus L(v') \leq K(v')$, the quotient map $\Z_2 \oplus \Lambda_1(A) \rightarrow \Lambda(v')$ factors as a composition $\Z_2 \oplus \Lambda_1(A) \rightarrow \Z_2 \oplus \Xi(v') \rightarrow \Lambda(v')$. 

\begin{thm} \label{t:lambda-comp}
%%%%%%%%%%%%%%%%
Let $v' \colon \Z_2 \to A$ be a homomorphism. There is a commutative diagram with short exact rows and columns
\[
\xymatrix{
 & 0 \ar[d] & 0 \ar[d] & & \\
 & \Coker(v'_2) \ar@{=}[r] \ar[d]_-{u'_{v'}} & \Coker(v'_2) \ar[d]_-{u'_{v'}} & & \\
0 \ar[r] & K(v') \ar[d]_-{r} \ar[r] & \Z_2 \oplus \Lambda_1(A) \ar[d] \ar[r] & \Lambda(v') \ar@{=}[d] \ar[r] & 0 \\
0 \ar[r] & \Z_2 \ar[d] \ar[r]^-{\iota} & \Z_2 \oplus \Xi(v') \ar[d] \ar[r] &\Lambda(v') \ar[r] & 0 \\
 & 0 & 0 & & 
}
\]
where 
\begin{compactitem}[$\bullet$]
\item $v'_2 = v' \otimes \Id_{\Z_2} \colon \Z_2 \otimes \Z_2 = \Z_2 \rightarrow A \otimes \Z_2$, hence $v'_2(1)=v'(1) \otimes 1$;
\item The homomorphism $\iota \colon \Z_2 \rightarrow \Z_2 \oplus \Xi(v')$ is given by $\iota(1) = (1,[v'(1) \ohat v'(1)])$;
\item The homomorphism $r \colon K(v') \rightarrow \Z_2$ is the composition $K(v') \rightarrow \Z_2 \oplus \Lambda_1(A) \rightarrow \Z_2$ of the natural inclusion  and projection maps;
\item The homomorphism $u'_{v'} \colon \Coker(v'_2) \rightarrow K(v') \leq \Z_2 \oplus \Lambda_1(A)$ is induced by the map $A \otimes \Z_2 \rightarrow K(v')$, $x \otimes 1 \mapsto (0,x \ohat x + x \ohat v'(1))$;
\item The unlabelled maps are the obvious inclusions/quotient maps. 
\end{compactitem}
\end{thm}

\begin{proof}
First consider the vertical sequence on the left. The map $A \otimes \Z_2 \rightarrow K(v')$ sends $v'(1) \otimes 1$ to $(0,v'(1) \ohat v'(1) + v'(1) \ohat v'(1)) = 0$, hence $u'_{v'}$ is well-defined. Suppose that $v' \cong v'_0 \oplus C \colon \Z_2 \rightarrow B \oplus C$ is a maximal splitting of $v'$. If $v'_0 \cong \iota_1$, then $A \otimes \Z_2 \cong \Z_2 \oplus C \otimes \Z_2$ and $\Coker(v'_2) \cong C \otimes \Z_2$. If $v' = 0$ or $v'_0 \cong \iota_k$, then $v'_2(1)=v'(1) \otimes 1 = 0$, so $\Coker(v'_2) \cong A \otimes \Z_2$. Moreover, the proof of Lemma \ref{lem:ker1} shows that $u'_{v'}$ maps $\Coker(v'_2)$ isomorphically onto $0 \oplus L(v') \leq K(v')$. The quotient $K(v') / (0 \oplus L(v'))$ is generated by $[1,v'(1) \ohat v'(1)]$, and $r$ induces an isomorphism onto $\Z_2$, showing that the sequence is exact. The map $\Z_2 \rightarrow K(v')$, $1 \mapsto (1,v'(1) \ohat v'(1))$ splits $r$. 

The vertical sequence in the middle is exact, because the kernel of $\Z_2 \oplus \Lambda_1(A) \rightarrow \Z_2 \oplus \Xi(v')$ is $0 \oplus L(v')$. The middle horizontal sequence is exact by the definition of $K(v')$ and $\Lambda(v')$. Finally, the bottom horizontal sequence is exact, because $K(v') = (0 \oplus L(v')) \oplus \An{ (1,v'(1) \ohat v'(1)) }$. 
\end{proof}

The projection $\Z_2 \oplus \Xi(v') \rightarrow \Z_2$ splits the bottom exact sequence, hence $\Lambda(v')$ is naturally isomorphic to $\Xi(v')$. 

If $v' = 0$, then $L(v')$ is the image of the natural embedding $A \otimes\Z_2 \rightarrow \Lambda_1(A)$, so $\Xi(v') = \Lambda(A)$. 

If $v' \cong \iota_1 \oplus C$, then, by the proof of Lemma \ref{lem:ker1}, $L(v')$ is the image of the map $\Z_2 \otimes C \rightarrow \Lambda_1(A) \cong \Z_2 \oplus \Z_2 \otimes C \oplus \Lambda_1(C)$, $1 \otimes x \mapsto (0,1 \otimes x, x \ohat x)$. This map can be extended to an automorphism of $\Lambda_1(A)$ by taking the identity on $\Z_2$ and $\Lambda_1(C)$. This automorphism shows that $\Xi(v') = \Lambda_1(A) / L(v') \cong \Lambda_1(A) / \Z_2 \otimes C \cong \Z_2 \oplus \Lambda_1(C)$. Note that $\Lambda_1(C) \cong 0 \oplus C \otimes \Z_2 \oplus \Lambda(C) \cong \Lambda(\Z_2 \oplus C) \cong \Lambda(A)$, hence $\Xi(v') \cong \Z_2 \oplus \Lambda(A)$, but this isomorphism is not natural.

If $v' \cong \iota_k \oplus C$, $k \geq 2$, then $L(v')$ is the direct sum of the image of the map $\Z_2 \otimes C \rightarrow \Lambda_1(A) \cong \Z_2 \oplus \Z_{2^k} \otimes C \oplus \Lambda_1(C)$, $1 \otimes x \mapsto (0, 2^{k-1} \otimes x, x \ohat x)$ and $\Z_2 \oplus 0 \oplus 0$. This time we regard $\Z_2 \otimes C$ as a subgroup of $\Lambda_1(C)$ via a splitting $\Lambda_1(C) \cong \Z_2 \otimes C \oplus \Lambda(C)$, and extend the above map to an automorphism of $\Lambda_1(A)$ by taking the identity on $\Z_2$, $\Z_{2^k} \otimes C$ and $\Lambda(C)$. We get that $\Xi(v') = \Lambda_1(A) / L(v') \cong \Lambda_1(A) / (\Z_2 \oplus \Z_2 \otimes C) \cong \Z_{2^k} \otimes C \oplus \Lambda(C) \cong \Lambda(\Z_{2^k} \oplus C) \cong \Lambda(A)$, but, again, this isomorphism is not natural.

%%%%%%%%%%%%%%%%%%%%%%%%%%%%%%%%%%%%%%%%%%%%%%%%%

\section{Computing Grothendieck-Witt groups and characterising absorbing forms} \label{s:GW}
%%%%%%%%%%%%%%%%%%%%%%%%%%%%%%%%%%%%%%%%%%%%%%%%%%%
In this section we recall the definition of the Grothendieck-Witt group of a form parameter $Q$ and describe its relationship to the Witt group of $Q$.  We also define and characterise absorbing $Q$-forms.

\subsection{The groups $GW_0(Q)$ and $W_0(Q)$}
%%%%%%%%%%%%%%%%%%%%%%%%%%%%%%%%%%%%%%%%%%%%%%%
Let $Q$ be a form parameter, $\ul \mu$ a nonsingular $Q$-form and $\{ \ul \mu \}$ the isometry class of $\ul \mu$.
Recall from Section~\ref{ss:GW-Intro} that $\ICHomZns(Q) = \big\{ \{ \ul \mu \} \big\}$ is the set of isomorphism classes of nonsingular $Q$-forms, which is a monoid under orthogonal sum,
and that the {\em Grothendieck-Witt group of $Q$} is the Grothendieck group of this monoid: 
\[ 
\GW_0(Q) := {\mcal Gr}\big( \ICHomZns(Q), \oplus \big).
\]
There is a natural surjection
\[ \pi \colon GW_0(Q) \to W_0(Q),
\quad \bigl( \{\ul \mu_0\}, \{\ul \mu_1\} \bigr) \mapsto 
[\ul \mu_0] - [\ul \mu_1],
\]
and by definition, there is a short exact sequence of abelian groups
\begin{equation} \label{eq:GW}
0 \to \Ker(\pi) \xra{} GW_0(Q) \xra{\pi} W_0(Q) \to 0.
\end{equation}
\noindent
This exact sequence allows us to compute $GW_0(Q)$ in terms of $W_0(Q)$, using the following definition and lemma.

\begin{defin}
For every form parameter $Q$, there is a well-defined rank homomorphism
\[
\rank \colon \GW_0(Q) \rightarrow \Z \text{\,,} \quad 
\big( \{\ul{\mu_0}\}, \{\ul{\mu_1}\}\big) \mapsto \rank(\ul{\mu_0}) - \rank(\ul{\mu_1}). 
\]
If $Q$ is symmetric, then there is also a signature homomorphism
\[
\sigma \colon \GW_0(Q) \rightarrow \Z \text{\,,} \quad \big( \{\ul{\mu_0}\}, \{\ul{\mu_1}\} \big) \mapsto \sigma(\ul{\mu_0}) - \sigma(\ul{\mu_1}).
\]
\end{defin}

It is well known that $\Ker(\pi) \cong \Z$ is generated by $\big(\{H_\epsilon(\Z)\}, \{ \ul 0\}\big)$, where we recall that $\ul 0$ is the $Q$-form on the trivial group;
see e.g.\ \cite[Lemma 2.31]{S}.
Indeed Schlichting \cite[Definition 2.14]{S} defines $W_0(Q)$ to be the quotient of $\GW_0(Q)$ by the subgroup of formal differences of hyperbolic forms.
Hence we have

\begin{lem} \label{lem:rk-isom}
For every form parameter $Q$ the restriction of the rank homomorphism $\rank \colon \GW_0(Q) \rightarrow \Z$ yields an isomorphism 
\[
\pushQED{\qed} 
\rk \big| _{\Ker(\pi)} \colon \Ker(\pi) \xra{\cong} 2\Z. \qedhere
\popQED
\]
\end{lem}

\noindent
{\em Remark.} We present an {\em a priori} alternative proof of 
Lemma~\ref{lem:rk-isom} below.

\begin{thm} \label{t:GW_splits}
a) For every form parameter $Q$ the map $\rank \oplus \pi \colon \GW_0(Q) \to \Z \oplus W_0(Q)$ is injective. 

b) If $Q$ is symmetric, then $(n, [\ul{\mu}]) \in \Im(\rank \oplus \pi)$ if and only if $\sigma([\ul{\mu}]) \equiv n~\mod~2$. If $Q$ is anti-symmetric, then $\Im(\rank \oplus \pi) = 2\Z \oplus W_0(Q)$.

c) The short exact sequence \eqref{eq:GW} splits.
\end{thm}

\begin{proof}
a) The restriction of $\rank \colon \GW_0(Q) \to \Z$ to $\Ker(\pi)$ is injective by Lemma \ref{lem:rk-isom}. 

b) If $Q$ is symmetric, then we have $\sigma(\ul{\mu})\equiv \rank(\ul{\mu})~\mod~2$ for every $Q$-form $\ul{\mu}$, therefore $\Im(\rank \oplus \pi) \subseteq \left\{ (n, [\ul{\mu}]) \mid \sigma([\ul{\mu}]) \equiv n~\mod~2 \right\}$. If $Q$ is anti-symmetric, then every nonsingular form over $Q$ has even rank, so $\Im(\rank \oplus \pi) \subseteq 2\Z \oplus W_0(Q)$. 
In the other direction, for any form parameter $Q$, every nonsingular $Q$-form $\ul \mu$
and $k \geq 0$, we have $(\rank \oplus \pi)([\ul{\mu} \oplus H_{\epsilon_Q}(\Z^k),\ul 0]) = (\rank(\ul \mu)+2k, [\ul{\mu}])$, and so equality holds for both symmetries.

c) If $Q$ is symmetric, then the signature map $\sigma \colon \GW_0(Q) \to \Im(\sigma)$ has a splitting $s \colon \Im(\sigma) \rightarrow \GW_0(Q)$ (because $\Im(\sigma) \subseteq \Z$), and $\GW_0(Q) = \Ker(\sigma) \oplus \Im(s)$. Since $\Ker(\pi) \subseteq \Ker(\sigma)$ and every form with zero signature has even rank, we can define a splitting map $f \colon \GW_0(Q) \rightarrow \Ker(\pi) \cong 2\Z$ by $f \big| _{\Ker(\sigma)} = \rank\big| _{\Ker(\sigma)}$ and $f \big| _{\Im(s)} = 0$.

If $Q$ is anti-symmetric and $\Ker(\pi)$ is identified with $2\Z$ using Lemma \ref{lem:rk-isom}, then the rank homomorphism $\rank \colon \GW_0(Q) \to 2\Z$ is a canonical splitting of the inclusion $\Ker(\pi) \to \GW_0(Q)$.
\end{proof}

We conclude this subsection with an alternative proof of Lemma~\ref{lem:rk-isom}.  Our main motivation
for doing this is that we believe Lemma~\ref{l:metabolic} is of independent interest in the study
of $Q$-forms; see Remark~\ref{r:Q-form-Class}.
 
\begin{lemma} \label{l:metabolic}
%%%%%%%%%%%%%%%%%%
Let $\ul \mu = (X, \lambda, \mu)$ and $\ul \mu' = (X', \lambda', \mu')$ be full metabolic $Q$-forms on isomorphic free abelian groups $X \cong X'$. Then for some integer $m \geq 0$, $\ul \mu \oplus H_{\epsilon_Q}(\Z^m)$ and $\ul \mu' \oplus H_{\epsilon_Q}(\Z^m)$ are isometric $Q$-forms.
\end{lemma}

\begin{proof}
If $Q$ is symmetric, then by Lemma \ref{l:gmtrc} we can regard a $Q$-form $(X, \lambda, \mu)$ as a $(Q^+ \oplus SQ)$-form $(X, \lambda, (\mu_0,\mu_1))$ such that $v(\mu_1(x)) \equiv \lambda(x, x)$ mod $2$. The analogous statement for such $(Q^+ \oplus SQ)$-forms is proved in \cite[Theorem 2.22 a)]{N2}.

Here we consider the anti-symmetric case. As $\ul \mu$ and $\ul \mu'$ are metabolic, they have lagrangians $L < X$ and $L' < X'$. Let $e_1, \ldots , e_k, f_1, \ldots , f_k$ and $e'_1, \ldots , e'_k, f'_1, \ldots , f'_k$ be symplectic bases for the bilinear forms $(X, \lambda)$ and $(X', \lambda')$ respectively such that $L = \left< e_1, \ldots , e_k \right>$ and $L' = \left< e_1', \ldots , e_k' \right>$, and let $N = \left< f_1, \ldots , f_k \right>$ and $N' = \left< f'_1, \ldots , f'_k \right>$. Since $L$ and $L'$ are lagrangians, we have $S\mu \big| _L = 0$ and $S\mu' \big| _{L'} = 0$, so $S\mu \big| _N \colon N \rightarrow SQ$ and $S\mu' \big| _{N'} \colon N' \rightarrow SQ$ are surjective (because $\ul \mu$ and $\ul \mu'$ are full). By \cite[Lemma 2.21 a)]{N2} there is an $m \geq 0$ and an isomorphism $h \colon N \oplus \Z^m \rightarrow N' \oplus \Z^m$ such that $(S\mu \big| _N + 0) = (S\mu' \big| _{N'} + 0) \circ h \colon N \oplus \Z^m \rightarrow SQ$. 

Let $\tilde{X} = X \oplus \Z^{2m}$, $\tilde{L} = L \oplus \Z^m$ and $\tilde{N} = N \oplus \Z^m$, so that $\tilde{X} = \tilde{L} \oplus \tilde{N}$ and if $H_-(\Z^m) = (\Z^m \oplus \Z^m, \lambda_H, \mu_H)$ (where each copy of $\Z^m$ is a lagrangian), then $\ul \mu \oplus H_-(\Z^m) = (\tilde{X}, \lambda \oplus \lambda_H, \mu \oplus \mu_H)$ and $\tilde{L}$ is a lagrangian in $\ul \mu \oplus H_-(\Z^m)$. Let $l = k+m$ and suppose that $e_{k+1}, \ldots , e_l, f_{k+1}, \ldots , f_l$ is a symplectic basis of $H_-(\Z^m)$ (with $\left< e_{k+1}, \ldots , e_l \right>$ and $\left< f_{k+1}, \ldots , f_l \right>$ being the standard lagrangians). Then $\tilde{L} = \left< e_1, \ldots , e_l \right>$, $\tilde{N} = \left< f_1, \ldots , f_l \right>$ and $e_1, \ldots , e_l, f_1, \ldots , f_l$ is a symplectic basis of $(\tilde{X}, \lambda \oplus \lambda_H)$. Similarly let $\tilde{X}' = X' \oplus \Z^{2m}$, $\tilde{L}' = L' \oplus \Z^m$ and $\tilde{N}' = N' \oplus \Z^m$ so that $\ul \mu' \oplus H_-(\Z^m) = (\tilde{X}', \lambda' \oplus \lambda_H, \mu' \oplus \mu_H)$. 

The map $h$ can be regarded as an isomorphism $\tilde{N} \rightarrow \tilde{N}'$ with $S(\mu \oplus \mu_H \big| _{\tilde{N}}) = S(\mu' \oplus \mu_H \big| _{\tilde{N}'}) \circ h \colon \tilde{N} \rightarrow SQ$. For every $1 \leq i \leq l$ let $f''_i = h(f_i)$, then $f''_1, \ldots , f''_l$ is a basis of $\tilde{N}'$. Since the adjoint of $\lambda' \oplus \lambda_H$ determines an isomorphism $\tilde{L}' \cong (\tilde{N}')^*$, we get a basis $e''_1, \ldots , e''_l$ of $\tilde{L}'$ such that $e''_1, \ldots , e''_l, f''_1, \ldots , f''_l$ is a symplectic basis of $(\tilde{X}', \lambda' \oplus \lambda_H)$. For every $1 \leq i \leq l$ we have $S(\mu \oplus \mu_H)(f_i) = S(\mu' \oplus \mu_H)(f''_i) \in SQ = Q_e / \Im({\rm p})$, so $(\mu' \oplus \mu_H)(f''_i) \in Q_e$ is either $(\mu \oplus \mu_H)(f_i)$ or $(\mu \oplus \mu_H)(f_i) + {\rm p}(1)$. We define $f'''_i := f''_i$ in the case where $(\mu' \oplus \mu_H)(f''_i) = (\mu \oplus \mu_H)(f_i)$ and $f'''_i := f''_i + e''_i$ otherwise. 

Then $e''_1, \ldots , e''_l, f'''_1, \ldots , f'''_l$ is again a symplectic basis of $(\tilde{X}', \lambda' \oplus \lambda_H)$. Moreover, we have the equalities $(\mu \oplus \mu_H)(e_i) = 0 = (\mu' \oplus \mu_H)(e''_i)$ and $(\mu \oplus \mu_H)(f_i) = (\mu' \oplus \mu_H)(f'''_i)$ by the definition of $f'''_i$. Therefore the isomorphism $I \colon \tilde{X} \rightarrow \tilde{X}'$ given by $I(e_i)=e''_i$ and $I(f_i) = f'''_i$ is an isometry between $\ul \mu \oplus H_-(\Z^m)$ and $\ul \mu' \oplus H_-(\Z^m)$. 
\end{proof}

Recall from Definition \ref{def:full} that a $Q$-form 
$(X, \lambda, \mu)$ is called full if its linearisation $S\mu \colon X \to SQ$ is surjective.

\begin{lem} \label{lem:full}
For any form parameter $Q$ there exists a full metabolic $Q$-form. 
\end{lem}

\begin{proof}
Let $Q = (Q_e, {\rm h}, {\rm p})$ and suppose that $q_1, q_2, \ldots , q_k$ is a generating system for $Q_e$. Consider the triple 
\[
\left( \Z^k \oplus \Z^k, 
\left[ \begin{array}{cc} 0 & I_k \\ \epsilon_Q I_k & D \end{array} \right], 
( 0, \ldots, 0, q_1, \ldots, q_k)^{\rm t}
\right),
\]
where $I_k$ is the $k \times k$ identity matrix, $D$ is a $k \times k$ diagonal matrix with entries ${\rm h}(q_1), \ldots, {\rm h}(q_k)$ on the diagonal and $(\cdot)^{\rm t}$ denotes the transpose. By Remark \ref{rem:eqf-not} this data defines a $Q$-form,
which is metabolic because the matrix is nonsingular and 
$\Z^k \oplus 0$ a lagrangian. It is also full, because if $\pi\colon Q_e \rightarrow SQ$ denotes the quotient homomorphism, then $\pi(q_1), \ldots, \pi(q_k)$ is a generating system for $SQ$.
\end{proof}

\begin{proof}[Alternative proof of Lemma~\ref{lem:rk-isom}]
Metabolic forms have even rank, and so $\rank(\Ker(\pi)) \subseteq 2\Z$. Moreover, we have $\rank([H_{\epsilon_Q}(\Z^m),\ul 0]) = 2m$, 
and so $\rank \colon \Ker(\pi) \to 2\Z$ is surjective. Finally assume that 
$\ul{\mu_0}$ and $\ul{\mu_1}$ are metabolic $Q$-forms of equal rank, and consider $\big(\{\ul{\mu_0}\},\{\ul{\mu_1}\}\big) \in \Ker(\rank \big| _{\Ker(\pi)})$.
By Lemma \ref{lem:full} there is a full metabolic $Q$-form $\ul{\mu'}$, and by applying Lemma \ref{l:metabolic} to $\ul{\mu_0} \oplus \ul{\mu'}$ and $\ul{\mu_1} \oplus \ul{\mu'}$ we get that $\ul{\mu_0} \oplus \ul{\mu'} \oplus H_{\epsilon_Q}(\Z^m) \cong \ul{\mu_1} \oplus \ul{\mu'} \oplus H_{\epsilon_Q}(\Z^m)$ for some $m$. This shows that 
$\big(\{\ul{\mu_0}\},\{\ul{\mu_1}\}\big)  = 0$, 
and therefore $\rank \big| _{\Ker(\pi)} \colon \Ker(\pi) \to 2\Z$ is an isomorphism. 
\end{proof}

\subsection{Absorbing $Q$-forms}
%%%%%%%%%%%%%%%%%%%%%%%%%%%%%%%%%%%%%%%%%%%%%
In this subsection we characterise absorbing $Q$-forms.

\begin{defin}
\quad
\begin{compactenum}[a)]
\item Let $\ul \mu$ and $\ul \eta$ be nonsingular $Q$-forms. We say that $\ul \mu$ \textit{absorbs} $\ul \eta$ if $\ul \eta$ can be embedded in $k\ul \mu = \ul \mu \oplus \ldots \oplus \ul \mu$ for some $k \geq 1$. 

\item A nonsingular $Q$-form $\ul \mu$ is called \textit{absorbing}, if it absorbs all nonsingular $Q$-forms.
\end{compactenum}
\end{defin}

\begin{thm}
A nonsingular $Q$-form is absorbing if and only if it is indefinite and full.
\end{thm}

\begin{proof}
A definite form cannot absorb any indefinite form, and a non-full form cannot absorb any full form, so the conditions are necessary. 

For the other direction first note that $\ul \mu$ absorbs $\ul \eta \oplus \ul \eta'$ if and only if it absorbs both $\ul \eta$ and $\ul \eta'$. As we saw before Definition \ref{def:WG}, $\ul \eta \oplus (-\ul \eta)$ is metabolic for every nonsingular $\ul \eta$.
Moreover, it elementary to prove that every metabolic form can be written as a direct sum of rank-$2$ metabolic forms,
and therefore a $Q$-form is absorbing if and only if it absorbs every rank-$2$ metabolic form.

Now let $\ul \mu = (X, \lambda, \mu)$ be a full indefinite nonsingular form over $Q = (Q_e, {\rm h}, {\rm p})$, and let $\pi \colon Q_e \rightarrow SQ$ denote the projection. Let $\ul \eta$ be a rank-$2$ metabolic form, it is isomorphic to 
\[
\left(
\Z^2\left< e,f \right>, 
\left[ \begin{array}{cc} 0 & 1 \\ \epsilon_Q & a \end{array} \right],
\left( \begin{array}{c} 0 \\ q \end{array} \right) 
\right)
\]
for some $q \in Q_e$ and $a = {\rm h}(q) \in \Z$. We will show that $\ul \eta$ embeds into $3\ul \mu$.

Since $\ul \mu$ is indefinite, there is a primitive $x \in X$ such that $\lambda(x,x) = 0$, and since it is nonsingular, there is a $y \in X$ such that $\lambda(x,y) = 1$. Finally, $\ul \mu$ is full, so there is a $z \in X$ such that $S\mu(z) = \pi(q)-S\mu(y)$. If $Q$ is symmetric, then let $\delta = \frac{1}{2}{\rm h}(q-\mu(y)-\mu(z))$, which is an integer, because $\pi(q-\mu(y)-\mu(z)) = \pi(q)-S\mu(y)-S\mu(z)=0$.
It follows that $q-\mu(y)-\mu(z) \in \Ker(\pi)= \Im({\rm p})$, and ${\rm h} \circ {\rm p} = 2 \Id_{\Z}$. If $Q$ is anti-symmetric, then let $\delta = 0$ if $\mu(y)+\mu(z) = q$ and $\delta = 1$ otherwise.

We will use the notation $(X^3, \bar{\lambda}, \bar{\mu})$ for $3\ul \mu$. We define the homomorphism $i \colon \Z^2 \rightarrow X^3$ by 
\[
i(e) = (x,-x,0) \quad \text{and} \quad i(f) = (y + \delta x,-\delta x,z).
\] 
Then $\bar{\lambda}(i(e),i(e)) = 0$ and $\bar{\lambda}(i(e),i(f)) = 1$. 
We also have 
\[ \bar{\mu}(i(e)) =  \mu(x) + \mu(-x) = \mu(x) + \mu(-x) + {\rm p}(\lambda(x,-x)) = \mu(x-x) =0\]
and $S\bar{\mu}(i(f)) = S\mu(y)+S\mu(z) = \pi(q)$. 
If $Q$ is symmetric, then 
\[ \bar{\lambda}(i(f),i(f)) = \lambda(y,y) + 2\delta + \lambda(z,z) = {\rm h}(\mu(y)) + {\rm h}(q-\mu(y)-\mu(z)) + {\rm h}(\mu(z)) = {\rm h}(q) = a.\] 
This means that ${\rm h}(\bar{\mu}(i(f))) = {\rm h}(q)$ and $\pi(\bar{\mu}(i(f))) = \pi(q)$, so by Lemma \ref{l:Wu-pb} $\bar{\mu}(i(f)) = q$. 
If $Q$ is anti-symmetric, then $\bar{\lambda}(i(f),i(f)) = 0 = a$. Moreover,
\[ \bar{\mu}(i(f)) = \mu(y+\delta x) + \mu(-\delta x) + \mu(z) = \mu(y)+\mu(\delta x) + \delta{\rm p}(1) + \mu(-\delta x) + \mu(z) = \mu(y) + \mu(z) + \delta{\rm p}(1).\] 
The previously established equality $S\mu(y)+S\mu(z) = \pi(q)$ implies that $\mu(y)+\mu(z)-q \in \Ker(\pi) = \Im({\rm p})$, so either $\mu(y) + \mu(z)=q$ or $\mu(y) + \mu(z)=q + {\rm p}(1)$. By the definition of $\delta$ we have $\bar{\mu}(i(f)) = q$ in both cases. 

Therefore $i$ is a morphism of $Q$-forms. Since $\ul \eta$ is nonsingular, this implies that $i$ is injective. Therefore $\ul \eta$ embeds into $3\ul \mu$, and by our earlier observations this shows that $\ul \mu$ is absorbing. 
\end{proof}

%%%%%%%%%%%%%%%%%%%%%%%%%%%%%%%%%%%%%%%%%%%%%%%%%%%%%%%%%%%%

\section{Witt groups and algebraic Poincar\'e cobordism} \label{s:APC}
%%%%%%%%%%%%%%%%%%%%%%%%%%%%%%%%%%%%%%%%%%%%%%%%%%

Extended quadratic forms and their Witt groups were previously
defined by Ranicki in \cite{R1}, where they are located in the larger context
of the bordism groups of algebraic Poincar\'e complexes with structure maps to chain bundles.
In this section we briefly recall this larger setting and relate Ranicki's algebraic 
Wu classes to the quasi-Wu classes of Definition~\ref{d:q-Wu}.

As the objects and definitions of algebraic surgery are involved, we will not go into 
them here in detail.  Rather, we assume that the reader is familiar with the theory of chain bundles $(B_*, \beta)$
and algebraic Poincar\'e complexes with $(B_*, \beta)$-structure as developed by Weiss and Ranicki \cite{We, R1};
we use \cite{R1} as a reference.  
We fix the ring $\Z$ with trivial involution and let $B_*$ be a chain complex of finitely generated 
free $\Z$-modules with $B_i = 0$ for all $i < 0$.
A chain bundle $(B_*, \beta)$ is an algebraic model for a stable vector bundle,
and more generally a stable spherical fibration.
Indeed, if $\xi$ is a spherical fibration over a simply connected $CW$-complex $X$,
then $\xi$ defines a chain bundle $(C_*(X), \beta(\xi))$ over the cellular chain complex of $X$.
An $n$-dimensional $(B_*, \beta)$-Poincar\'e complex
is an algebraic analogue of closed $n$-dimensional normal $\xi$-manifold and 
the $\beta$-twisted $L$-groups, $L^n(B_*, \beta)$, which are the algebraic bordism classes of
$n$-dimensional $(B_*,\beta)$-Poincar\'e complexes, model 
the bordism group of closed $n$-dimensional normal $\xi$-manifolds.

The $\beta$-twisted $L$-groups are natural for maps of chain bundles,
and by definition, $L_n(\Z)$, the usual quadratic $L$-group of the ring $\Z$ (see e.g.\ \cite[Defns.\ 11.57 \& 12.23]{R2}),
is $L_n(0, 0)$, where $(0, 0)$ denotes the unique chain bundle over the zero chain complex.
It follows 
that there is a natural map $L_n(\Z) \to L^n(B_*, \beta)$, which 
Weiss \cite[2.21(iii) \& (iv)]{We} places into a long exact sequence of abelian groups
\begin{equation} \label{eq:LBbeta}
\dots \xra{} Q_{n+1}(B_*, \beta) \xra{} L_n(\Z) \xra{} L^n(B_*, \beta) \xra{} Q_n(B_*, \beta) \xra{} L_{n-1}(\Z) \xra{} \dots ~ ,
\end{equation}
where $Q_n(B_*, \beta)$ is the group of $n$-dimensional $\beta$-twisted quadratic structures on $B_*$ defined in \cite[\S 6]{R1}.

The $\beta$-twisted $L$-groups receive homomorphisms from certain Witt groups, as we now explain.
For each non-negative integer $q$, Ranicki defines a $(-1)^q$-symmetric form parameter $Q_{(B_*, \beta)}(q)$ \cite[\S 10]{R1}, which for brevity we denote by $Q_\beta(q)$.
Moreover, given a nonsingular $Q_\beta(q)$-form $(X, \lambda, \mu)$,
Ranicki constructs a $2q$-dimensional
$(B_*, \beta)$-Poincar\'e complex concentrated in dimension $q$ \cite[p.\ 250]{R1}, and proves that this construction induces a natural homomorphism
\[ \eta_\beta \colon W_0(Q_\beta(q)) \to L^{2q}(B_*, \beta).\]

Below we formulate a conjecture which describes conditions under which the homomorphism $\eta_\beta$ may be a split injection or an isomorphism.  At all events, to use $\eta_\beta$ it is important to identify the form parameter $Q_\beta(q)$ by finding the quasi-Wu class of $Q_\beta(q)$, and we do this now.
By construction (see \cite[\S 10]{R1}) 
the linearisation of the form parameter $Q_\beta(q)$ is given by
\[ SQ_\beta(q) = H_q(B_*).\]
The chain bundle $\beta$ also defines algebraic Wu classes $\wh v_i(\beta) \in H^i(B_*; \Z_2)$ \cite[\S 4]{R1} 
such that if $\beta = \beta(\xi)$ is the chain bundle defined by a spherical fibration then,
$\wh v_i(\beta(\xi)) = v_i(\xi)$, where $v_i(\xi)$ are the usual Wu classes of $\xi$.
By the Universal Coefficient Theorem, we have a split short exact sequence
\[ 0 \to \Ext_1^\Z(H_{q-1}(B_*), \Z_2) \to H^q(B_*; \Z_2) \xra{\ev} \Hom(H_q(B_*), \Z_2) \to 0.\]
Recalling Definition~\ref{d:q-Wu}, we note that the quasi-Wu classes of $Q_\beta(q)$ are such that if $q = 2k$ is even then $v_{Q_\beta(2k)} \in \Hom(H_q(B_*), \Z_2)$
and if $q = 2k{+}1$ is odd, then either $v'_{Q_\beta(2k+1)} = 0$, or $v'_{Q_\beta(2k+1)} \neq 0$ 
is determined by the extension
\[ [v'_{Q_\beta(2k+1)} ] := \bigl[ 0 \to \Z_2 \to (Q_\beta(2k{+}1))_e \to H_{2k+1}(B_*) \to 0 \bigr] .\]

\begin{proposition} \label{p:Q_beta}
The algebraic Wu classes $\wh v_*(\beta)$ of a chain bundle $(B_*, \beta)$ determine the quasi-Wu classes 
of the form parameter $Q_\beta(q)$ as follows:
\begin{enumerate}
\item[a)]
If $q = 2k$, then $v_{Q_\beta(2k)} = \ev(\wh v_{2k}(\beta))$.
\item[b)]
If $q = 2k{+}1$ and $\ev(\wh v_{2k+2}(\beta)) \neq 0$ then then $v'_{Q_\beta(2k+1)} = 0$.
\item[c)]
If $q=2k{+}1$ and $\ev(\wh v_{2k+2}(\beta)) = 0$, then $\wh v_{2k+2}(\beta) \in \Ext^\Z_1(H_{2k+1}(B_*), \Z_2)$
and $[v'_{Q_\beta(2k+1)}] = \wh v_{2k+2}(\beta)$.
\end{enumerate}
\end{proposition}

\begin{remark}
%%%%%%%%%%%%%%
For $q = 2k{+}1$ odd in Proposition~\ref{p:Q_beta}, we observe that $v'_{Q_\beta(q)} = 0$ if and only if $\hat v_{2k+2}(\beta) = 0$,
and the classification of anti-symmetric form parameters $Q$ shows that $Q_-$ is a summand of $Q$ if and only if 
$v'_Q = 0$.  By Theorem \ref{t:wf-ind}
the Arf invariant is defined for $Q$-forms if and only if $v'_Q = 0$.
Hence we obtain an algebraic analogue of results of Browder \cite[\S 3]{Browd1} and Brown \cite[Corollary 1.22]{Brown}, stating that $v_{2q+2}(\xi) = 0$ is sufficient to define an Arf-invariant for closed $(4q{+}2)$-dimensional normal $\xi$-manifolds, where normal $\xi$-manifolds are defined in Section~\ref{ss:scm}.
Much deeper algebraic analogues of the work of Browder and Brown on quadratic refinements can be found in \cite{We}.
\end{remark}

\begin{proof}[Proof of Proposition~\ref{p:Q_beta}]
%%%%%%%%%%%%%
a) This follows directly by inspecting the formula for ${\rm h}$ in \cite[p.\,249]{R1} (it is denoted there by $s$), 
and the formula for 
$\widehat v_{2k}$ on \cite[p.\,232]{R1}.

b) The braid \cite[p.\,249]{R1} gives 
an exact sequence of abelian groups
\[ H_{2k+2}(B) \xra{\ev(\wh v_{2k+2}(\beta))} (Q_{-})_e \xra{} Q_\beta(2k{+}1)_e \to H_{2k+1}(B) \to 0, \]
where $(Q_{-})_e  \to Q_\beta(2k{+}1)$ 
is $v'_{Q_\beta(2k{+}1)}$.  The exactness of this sequence proves 
Part b).

c) This follows by inspecting the formula on \cite[p.\,249]{R1} for addition in $Q_\beta(2k{+}1)$.
\end{proof}

With Proposition~\ref{p:Q_beta} and Theorem~\ref{t:ES+EQL}
we can compute the Witt group $W_0(Q_\beta(q))$ from knowledge of $H_q(B)$ and the algebraic Wu classes of $\beta$.  To apply these computations to the $\beta$-twisted $L$-groups and $\beta$-twisted $Q$-groups, we need more information about
the homomorphism $\eta_\beta$.

\begin{conjecture}
\label{c:cbsurg}
%%%%%%%%%%%%%%%%%%%%%%%%%%
Let $(B_*, \beta)$ be a chain bundle and $q \geq 0$.
If $H_{q-1}(B_*)$ is torsion free, then the natural homomorphism 
$\eta_\beta \colon W_0(Q_\beta(q)) \to  L^{2q}(B_*, \beta)$ 
is split injective.
If, in addition, $H_i(B_*) = 0$ for all $i < q$, then $\eta_\beta$ is an isomorphism.
\end{conjecture}

\begin{remark}
One idea for constructing an inverse to $\eta_\beta$ goes as follows.
In \cite[\S 10]{R1} Ranicki proves that the $q$-th homology group of a $2q$-dimensional $(B_*, \beta)$-Poincar\'e complex supports a $Q_\beta(q)$-form.  If the map from $(B_*, \beta)$-Poincar\'e complex to $B_*$ is $q$-connected, then since $H_{q-1}(B_*)$ is torsion free,
Poincar\'e duality ensures that this $Q_\beta(q)$-form is nonsingular,
and we can take its Witt class a candidate for the value of $\eta_\beta^{-1}$.
If the map from the $(B_*, \beta)$-Poincar\'e complex to $B_*$ is not $q$-connected, we expect that algebraic surgery below the middle dimension, e.g.\  \cite[Proposition 5.4]{We}, can be used to show both that the $(B_*, \beta)$-Poincar\'e complex is $(B_*, \beta)$-cobordant to a $(B_*, \beta)$-Poincar\'e complex-connected complex with a $q$-connected map to $B_*$, and that the Witt class of the $Q_\beta(q)$-forms obtained is well-defined.
\end{remark}

To apply Conjecture~\ref{c:cbsurg} we need a further definition in the symmetric case. This arises since Weiss's sequence \eqref{eq:LBbeta} involves the quadratic $L$-group $L_{4k}(\Z) = W_0(Q_{+})$, whereas our methods compare $W_0(Q_\beta(2k))$ with the symmetric $L$-group $W_0(Q^{+})= L^{4k}(\Z)$.
For $v \colon A \to \Z_2$ a homomorphism, we recall $\Sigma(v) \subset \Z \oplus \Gamma(A)$ 
from Definition~\ref{d:Sigma(v)}. 
We make the following mod~$8$ modification of $\Sigma(v)$.
\begin{definition} \label{d:Sigma8}
%%%%%%%%%%%%%%%%%%
Let $\Sigma_8(v) \leq \Z/8 \oplus \Gamma(A)$ denote the subgroup generated by the elements
\[
(1,x \otimes 1),
\quad \text{and} \quad
(0,[k_1, k_2] \otimes 1),
\]
where $x \in v^{-1}(1)$ and $k_1, k_2 \in \Ker(v)$.
\end{definition}
The next conjecture would be a direct consequence of Conjecture \ref{c:cbsurg}, Proposition~\ref{p:Q_beta} and Theorem~\ref{t:ES+EQL}. Recall that $\iota_1 \colon \Z_2 \to \Z_2$ is the identity.  For any homomorphism $v' \colon \Z_2 \to A$, there is a canonical map of homomorphisms $\iota_1 \to v'$, which induces a homomorphism $\Lambda(\iota_1) \to \Lambda(v')$.

\begin{conjecture}  \label{c:LBbeta}
%%%%%%%%%
Let $(B_*, \beta)$ be a chain bundle with 
algebraic Wu class $\wh v_{2k}(\beta)$. \\
If $(B_*, \beta)$ is $(2k{-}1)$-connected then
\begin{enumerate}
\item[a)]
$L^{4k}(B_*, \beta) \cong W_0(Q_{\beta}(2k)) \cong \Sigma
\bigl(\ev(\wh v_{2k}(\beta))\bigr)$;
\item[b)]
$Q_{4k}(B_*, \beta) \cong  \Sigma_8
\bigl(\ev(\wh v_{2k}(\beta))\bigr)$.
\end{enumerate}
If $(B_*, \beta)$ is $(2k{-}2)$-connected 
and $\ev(\wh v_{2k}(\beta)) = 0$, then $[v'_{Q_\beta(2k{-}1)}] = \wh v_{2k}(\beta)$ and
\begin{enumerate}
\item[c)]
$L^{4k-2}(B_*, \beta) \cong W_0(Q_{\beta}(2k{-}1)) \cong \Lambda(v'_{Q_\beta(2k{-}1)})$;
\item[d)]
$Q_{4k-2}(B_*, \beta) \cong 
\Coker \bigl(\Lambda(\iota_1) \to \Lambda(v'_{Q_\beta(2k{-1})}) \bigr)$.
\end{enumerate}
If $(B_*, \beta)$ is $(2k{-}2)$-connected and 
$\ev(\wh v_{2k}(\beta)) \neq 0$ then $Q_\beta(2k{-1}) \cong Q^- \oplus H_{2k-1}(B_*)$ and
\begin{enumerate}
\item[e)]
$L^{4k-2}(B_*, \beta) \cong Q_{4k-2}(B_*, \beta) \cong    W_0(Q_{\beta}(2k{-}1)) \cong \Lambda(H_{2k-1}(B_*))$.
\end{enumerate}
\end{conjecture}

\begin{remark}
The conjectural computation of $\beta$-twisted $Q$-groups in Conjecture~\ref{c:LBbeta}
is consistent with the computations of Banagl and Ranicki \cite[Proposition 52]{B-R}.
\end{remark}

%%%%%%%%%%%%%%%%%%%%%%%%%%%%%%%%%%%%%%%%%%%%%%%%
%%%%%%%%%%%%%%%%%%%%%%%%%%%%%%%%%%%%%%%%%%%%%%%%

%%%%%%%%%%%%%%%%%%%%%%%%%%%%%%%%%%%%%%%%%%%%%%%%%

\end{document}